\documentclass[oneside,english]{amsart}
\usepackage[T1]{fontenc}
\usepackage[latin9]{inputenc}
\usepackage{babel}
\usepackage{amstext}
\usepackage{amsthm}
\usepackage{amssymb}
\usepackage{graphicx}
\usepackage{mathtools} 
\usepackage{graphicx}
\usepackage{tikz}
\usetikzlibrary{arrows.meta}
\usepackage[unicode=true,
 bookmarks=false,
 breaklinks=false,pdfborder={0 0 1},backref=section,colorlinks=false]
 {hyperref}

\makeatletter
\numberwithin{equation}{section}
\numberwithin{figure}{section}
\theoremstyle{plain}
\newtheorem{thm}{\protect\theoremname}
\theoremstyle{plain}
\newtheorem{cor}[thm]{\protect\corollaryname}
\theoremstyle{plain}
\newtheorem{lem}[thm]{\protect\lemmaname}
\theoremstyle{plain}
\newtheorem{prop}[thm]{\protect\propositionname}
\theoremstyle{remark}
\newtheorem{rem}[thm]{\protect\remarkname}
\theoremstyle{definition}
\newtheorem{defn}[thm]{\protect\definitionname}
\theoremstyle{definition}
\newtheorem{example}[thm]{\protect\examplename}

\usepackage{babel}
\numberwithin{thm}{section}

\providecommand{\corollaryname}{Corollary}
\providecommand{\definitionname}{Definition}
\providecommand{\examplename}{Example}
\providecommand{\lemmaname}{Lemma}
\providecommand{\propositionname}{Proposition}
\providecommand{\remarkname}{Remark}
\providecommand{\theoremname}{Theorem}

\makeatother

\providecommand{\corollaryname}{Corollary}
\providecommand{\definitionname}{Definition}
\providecommand{\examplename}{Example}
\providecommand{\lemmaname}{Lemma}
\providecommand{\propositionname}{Proposition}
\providecommand{\remarkname}{Remark}
\providecommand{\theoremname}{Theorem}

\newcommand{\CC}{\mathbb{C}}
\newcommand{\DD}{\mathbb{D}}

\newcommand{\di}[1]{\mathop{\mathrm{d}#1}}
\newcommand{\iu}{i}
\newcommand{\dist}{\mathrm{dist}}

\usepackage{tikz}
\usepackage{subcaption}
\usetikzlibrary{calc,arrows.meta,decorations.markings}

\title{The Brown measure of a sum of two free nonselfadjoint random variables,
one of which is $R$-diagonal}
\author{Hari Bercovici, Yu Kitagawa, and Ping Zhong}
\address{Mathematics Department, Indiana University, Bloomington, IN 47405,
USA}
\email{bercovic@indiana.edu }
\address{Department of Mathematics, Faculty of Science, Kyoto University, Japan}
\email{kitagawa.yu.57z@st.kyoto-u.ac.jp}
\address{Department of Mathematics, University of Houston, Houston,
TX 77204, USA}
\email{pzhong@central.uh.edu}
\subjclass[2000]{Primary: 46L54. Secondary: 30D05}

\begin{document}

\begin{abstract}
Suppose that $X_1$ and $X_2$ are two $*$-free, possibly unbounded,
random variables with Brown measures $\mu_{X_1}$ and $\mu_{X_2}$ respectively. Using properties of classic free additive convolution, we
develop a method for calculating $\mu_{X_1+X_2}$ when $X_2$ is
$R$-diagonal. This method determines a density relative to Lebesgue
measure on an open set whose closure equals the support of $\mu_{X_{1}+X_{2}}$. 
We prove that $\mu_{X_1+X_2}$ has no singular continuous part,
that its density is strictly positive on this open set, and
we determine the masses of its atoms. Explicit calculations
are possible in several important cases. As an application,
we combine our Fuglede--Kadison determinant formula with known subordination
convergence results to calculate the limiting logarithmic
potentials of free multiplicative random walks. We identify
the limit with the Fuglede--Kadison determinant formula for
free multiplicative Brownian motion previously proved by some PDE methods.  
\end{abstract}

\keywords{Free probability, Brown measure, $R$-diagonal operators, free convolution, nonselfadjoint random variables}

\maketitle

\tableofcontents

\section{Introduction}

Consider a tracial $W^{*}$-probability space $(\mathcal{A},\tau)$.
Thus, $\mathcal{A}$ is a von Neumann algebra, and $\tau$ is a normal,
faithful, tracial state on $\mathcal{A}$. In his extension of Lidski\u{\i}'s
theorem \cite[Theorem 1]{L-brown} to this context (and, indeed, to
the larger context of semifinite von Neumann algebras), L. Brown \cite{L-brown}
associated a kind of spectral distribution measure $\mu_{X}$ to possibly
unbounded operators $X$ affiliated with $\mathcal{A}$. We use the
notation $\widetilde{\mathcal{A}}$ for the algebra consisting of
all such operators, and we denote by ${\rm Log}^{+}(\tau)$ the collection
of those $X\in\widetilde{A}$ for which 
\[
\tau(\log^{+}|X|)<+\infty.
\]
Given $X\in{\rm Log}^{+}(\tau)$, Brown shows that the function 
\[
\lambda\mapsto\tau(\log|X-\lambda|),\quad\lambda\in\mathbb{C},
\]
is subharmonic and, in fact, it is the logarithmic potential of a
Borel probability measure $\mu_{X}$ on $\mathbb{C}$ such that 
\[
\int_{\mathbb{C}}\log^{+}|z|\,d\mu_{X}(z)<+\infty.
\]
In other words, 
\[
\tau(\log|X-\lambda|)=\int_{\mathbb{C}}\log|z-\lambda|\,d\mu_{X}(z),\quad z\in\mathbb{C}.
\]
The measure $\mu_{X}$ is obtained as 
\[
d\mu_{X}(\lambda)=\frac{1}{2\pi}\Delta_{\lambda}\tau(\log|X-\lambda|),
\]
where $\Delta_{\lambda}=4\overline{\partial}_{\lambda}\partial_{\lambda}$
denotes the Laplace operator, applied in the sense of distributions
to the locally integrable function $\tau(\log|X-\lambda|)$. The concept
of the Brown measure overlaps the classical one of the scalar spectral
measure for normal operators. To clarify, suppose that $X\in\widetilde{\mathcal{A}}$
is a normal operator (not necessarily in ${\rm Log}^{+}(\mu)$) with
spectral measure $E_{X}$, and set 
\[
\mu_{X}(\sigma)=\tau(E_{X}(\sigma))
\]
for every Borel set $\sigma\subset\mathbb{C}$. If $X\in{\rm Log^{+}(\mu)}$, this is precisely the Brown
measure of $X$, so the notation is consistent.

Suppose now that $X_{1},X_{2}\in{\rm Log}^{+}(\tau)$ are $*$-free
(relative to $\tau$; see, for instance \cite{vdn}) and $X=X_{1}+X_{2}$.
We wish to calculate the Brown measure $\mu_{X}$ in terms of data
coming from $X_{1}$ and $X_{2}$. When $X_{1}$ and $X_{2}$ are
selfadjoint, the resulting measure $\mu_{X}$ depends only on the
measures $\mu_{X_{1}}$ and $\mu_{X_{2}}$, namely 
\[
\mu_{X}=\mu_{X_{1}}\boxplus\mu_{X_{2}},
\]
where $\boxplus$ denotes the additive free convolution \cite{v-add}.
The nonselfadjoint situation is more difficult to deal with, and progress
has only been made in special situations. After one important case
calculated by Haagerup and Schulz \cite{HA-SC-brown}, Biane and Lehner
\cite{biane-lehner} developed techniques that work in greater generality.
In this work, we consider the case in which $X_{2}$ is $R$-diagonal
(see \cite{Ni-Sp-waterloo} and \cite{HA-SC-brown}) but no further
hypothesis is made about $X_{1}$. Previously, the case of a selfadjoint
$X_{1}$ (and circular $X_{2}$) was studied in \cite{Ho-Pi-jems}
and subsequently \cite{pzh-ajm} using different approaches. We further
develop the methods of \cite{pzh-ajm} without assuming that $X_{1}$
is selfadjoint.

\begin{thm}[Informal statement of Theorems \ref{thm:positivity} and  \ref{thm:atom-main}]
 There exists an open set
$\Omega\subset\mathbb{C}$, whose closure equals to the support of
$\mu_{X}$, and such that $\mu_{X}$ has a real-analytic and strictly positive density on
$\Omega$ relative to planar Lebesgue measure. Moreover, if $X_1$ is normal and $\lambda_0$ is an atom of $\mu_{X_1}$, then 
\[
  \mu_{X_1+X_2}(\{\lambda_0\})=(\mu_{X_1}(\{\lambda_0\})+\mu_{X_2}(\{0\})-1)_+.
\]
\end{thm}

The main ideas of our approach is based on Hermitian reduction method and on the fact that the calculation of the Brown measure for the sum of two $R$-diagonal operators can be reduced to calculating an ordinary free convolution of two symmetric probability measures on the real line. The latter fact was proved in
\cite{NSbook} for bounded operators, reformulated in that case in
\cite{ha-larsen}, and extended in \cite{HA-SC-brown} to unbounded
operators. More precisely, if $X_{1}$
and $X_{2}$ are $R$-diagonal operators then the symmetrization
$\widetilde{\mu}_{|X|}$ of $\mu_{|X|}$ is the usual free convolution
of $\widetilde{\mu}_{|X_{1}|}$ with $\widetilde{\mu}_{|X_{2}|}$. This allows us to use the well understood mechanism
of free convolution to gather information about $\mu_{X}$. 

Brown measures have been useful in several ways. First, of course,
is Brown's theorem \cite{L-brown} extending Lidski\u{\i}'s theorem
as well as \cite[Theorem 2]{FuKa-det}. Subsequently, Brown measures
were instrumental in parametrizing an important family of hyperinvariant
subspaces for arbitrary operators in type II$_1$ factors \cite{HA-SC-inv} (see also \cite{DH-circular,DH-DT,SS-invariant}).
Knowledge of $\mu_{X}$ is often predictive of the asymptotic behavior
of the empirical eigenvalue distribution of certain random matrix
ensembles as their size tends to $+\infty$. For instance, the Brown measure of an $R$-diagonal operator is the limiting distribution of single ring theorem \cite{GKZ}. The calculation
of $\mu_{X}$ for selfadjoint $X_{1}$ (\cite{Ho-Pi-jems} and \cite{pzh-ajm})
has a random matrix counterpart \cite{bordenave-caputo-chafai}. The
results of the present paper have their own random matrix counterpart
that is the subject of subsequent works \cite{Ho-Pi, HYZ}.

A further application concerns Biane's free multiplicative Brownian
motion $\{b_t\}_{t\geq0}$ \cite{BianeSB,BianeBM}.
Driver--Hall--Kemp calculated its Brown measure using the
Hamilton--Jacobi method \cite{DHK}. Ho and the third author extended
this approach to $ub_t$, where $u$ is a unitary operator free from
$b_t$ \cite{Ho-Pi-jems}.
Related matrix models were studied in the physics literature:
Gudowska-Nowak, Janik, Jurkiewicz, and Nowak considered
multiplicative diffusion \cite{GNJJN03,GNJJN05}, while Burda,
Janik, and Nowak developed a multiplication law and an
$S$-transform for non-Hermitian random matrices \cite{BJN}.
After the first version of our paper appeared on arXiv,
Driver--Hall--Ho--Kemp--Nemish--Nikitopoulos--Parraud
\cite{Seven} proved eigenvalue convergence for matrix random walks
and, using our results on the additive model, established
convergence of the Brown measures of free multiplicative random
walks $B_k(t)$ to that of $b_t$.
Combining their subordination convergence results with our
Fuglede--Kadison determinant formula, we prove locally uniform
convergence of the logarithmic potentials on compact subsets of the interior of the Brown measure support. We identify the limit with the potential
calculated in \cite{DHK}.
These calculations also suggest a completely new approach to the Brown measure
of $ub_t$ through additive subordination using free probability method, which we intend to develop in subsequent work.

The paper is organized as follows. In Section 2, we establish
identities for logarithmic integrals using the subordination
functions for free additive convolution. In Section 3, we apply
these identities to calculate the Brown measure of $X_1+X_2$,
where $X_1$ and $X_2$ are $*$-free and $X_2$ is $R$-diagonal.
In Section 4.1, we prove that the continuous part of this measure
is absolutely continuous with respect to planar Lebesgue measure.
In Section 4.2, we characterize its atoms and determine their
masses from the local behavior of the logarithmic potential.
We also give examples in which the formulas can be evaluated
explicitly in Section 5. Finally, we apply the determinant formula to free
multiplicative random walks in Section 6. Using the subordination convergence
results of \cite{Seven}, we prove locally uniform convergence
of their logarithmic potentials in the interior of the support and identify the limit with the potential of free multiplicative
Brownian motion calculated previously by some PDE methods in \cite{DHK}.

\section{Logarithmic integrals and subordination\label{sec:Logarithmic-integrals-and}}

As in the introduction, we work in the context of a tracial $W^{*}$-probability
space $(\mathcal{A},\tau)$. Given a selfadjoint operator $T\in\widetilde{\mathcal{A}}$,
we consider the functions 
\[
G_{T}(z)=\tau((z-T)^{-1}),\quad H_{T}(z)=\frac{1}{G_{T}(z)}-z,
\]
defined for $z$ in the complex upper half-plane $\mathbb{C}^{+}$.
We also employ the notation 
\[
G_{\mu}(z)=\int_{\mathbb{R}}\frac{d\mu(t)}{z-t},
\]
for an arbitrary Borel probability measure $\mu$ on $\mathbb{R},$
so $G_{T}=G_{\mu_{T}}$. The function $G_{T}$ maps $\mathbb{C}^{+}$
to $-\mathbb{C}^{+}$. Unless $\mu_{T}$ is a point mass, $H_{T}$
maps $\mathbb{C}^{+}$ to itself. The following result, containing
facts established in \cite{v-fish1,bi-mz} and \cite{be08,be14},
describes the subordination property of free convolution. The notation
$\overline{\mathbb{C}^{+}}$ indicates the closure of $\mathbb{C}^{+}\cup\mathbb{R}$
of $\mathbb{C}^{+}$ in the complex plane, while $\overline{\mathbb{C}^{+}}\cup\{\infty\}$
is the closure of $\mathbb{C}^{+}$ in the Riemann sphere. 
\begin{thm}
\label{thm:subordination}Suppose that $T_{1},T_{2}\in\widetilde{\mathcal{A}}$
are selfadjoint operators free relative to $\tau$, and $T=T_{1}+T_{2}$.
If neither $T_{1}$ nor $T_{2}$ is a constant multiple of the identity
operator, then there exist unique continuous functions $\omega_{1},\omega_{2}:\overline{\mathbb{C}^{+}}\to\overline{\mathbb{C}^{+}}\cup\{\infty\}$
that are analytic on $\mathbb{C}^{+}$ such that 
\begin{equation}
G_{T}(z)=G_{T_{1}}(\omega_{1}(z))=G_{T_{2}}(\omega_{2}(z)),\quad z\in\mathbb{C}^{+},\label{eq:G_jcompose with omega_j_}
\end{equation}
and 
\begin{equation}
\omega_{1}(z)+\omega_{2}(z)=z+\frac{1}{G_{T}(z)},\quad z\in\mathbb{C}^{+}.\label{eq:omega1+omega2}
\end{equation}
\end{thm}

The values of the functions $\omega_1$ and $\omega_2$ in Theorem \ref{thm:subordination} can be obtained as Denjoy-Wolff points of certain functions (\cite{denjoy} and \cite{wolff};
see also \cite{shapiro}). For an arbitrary analytic function
$\varphi:\mathbb{C}^{+}\to\mathbb{C}^{+}$, other than a conformal
automorphism: the iterates 
\[
\varphi^{\circ n}=\underbrace{\varphi\circ\cdots\circ\varphi}_{n\text{ times}}
\]
converge pointwise to a constant function, whose value (possibly in
$\mathbb{R}\cup\{\infty\}$) is now called the Denjoy-Wolff point
of $\varphi$. If $\varphi$ has a fixed point $\alpha\in\mathbb{C}^{+}$,
then the Denjoy-Wolff point of $\varphi$ is $\alpha$, and 
\[
|\varphi'(\alpha)|<1.
\]
If $\varphi$ does not have a fixed point in $\mathbb{C}^{+}$, then
the Denjoy-Wolff point $\alpha$ belongs to $\mathbb{R}\cup\{\infty\}$
and it is fixed in the sense that the nontangential limit of $\varphi$
at $\alpha$ exists and is equal to $\alpha$. In this case we have
$0<\varphi'(\alpha)\le1$, where $\varphi'(\alpha)$ denotes the Julia-Carath\'{e}odory
derivative. The number $\varphi'(\alpha)$ can be understood either
as the nontangential derivative of $\varphi$ at $\alpha$, or as
the nontangential limit of $\varphi'$ at $\alpha$. 

\begin{thm}[{\cite{bbh}}]
\label{thm:dw-subordination}
Let $T_1,T_2,T$ and $\omega_1,\omega_2$ be as in
Theorem~\ref{thm:subordination}. For
$z\in\overline{\mathbb{C}^{+}}$, define
\[
\varphi_z(w)=z+H_{T_2}(z+H_{T_1}(w)),\qquad
\psi_z(w)=z+H_{T_1}(z+H_{T_2}(w)),\quad w\in\mathbb{C}^{+}.
\]
Both maps are analytic self-maps of $\mathbb{C}^{+}$. If
$\varphi_z$ is not the identity, its Denjoy--Wolff point is
$\omega_1(z)$. The analogous statement
holds for $\psi_z$ and $\omega_2(z)$.
If neither map is a conformal automorphism and
$\omega_1(z),\omega_2(z)\in\mathbb{C}^{+}$, then
\[
|\varphi_z'(\omega_1(z))|
=|\psi_z'(\omega_2(z))|
=|H_{T_1}'(\omega_1(z))H_{T_2}'(\omega_2(z))|<1.
\]
\end{thm}
A continuity result for Denjoy--Wolff points is instrumental, and we recall it here for later reference.

\begin{thm}[{\cite{he-ajm, bbh}}]
\label{thm:dw-continuity}
    Let $f_n,f\colon\DD\to\overline{\DD}$ be analytic maps such that $f_n$ converges pointwise to $f$. If $f$ is not the identity map, then the Denjoy--Wolff points of $f_n$ converge to the Denjoy--Wolff point of $f$.
\end{thm}

When both $T_{1}$ and $T_{2}$ have two-point spectra, $\omega_{1}(z)$
is still a fixed point of $\varphi_{z}$, which is a conformal automorphism
in this case. The map $\varphi_{z}$ has a unique fixed point unless
it is the identity map, a situation that arises for at most one value
of $z$ for which $\varphi_{z}$ is the identity map. The spectrum
hypothesis in the statement above simply eliminates this rare occurrence. 

We note a simple consequence. 
\begin{cor}
\label{cor:when omega1(0) is in H}With the notation of Theorem \emph{\ref{thm:subordination}},
suppose that $z\in\mathbb{R}.$ Then $\omega_{1}(z)$ belongs to $\mathbb{C}^{+}$
if and only if $\omega_{2}(z)$ belongs to $\mathbb{C}^{+}.$ When
this occurs, $G_{T}$ extends to $z$ by continuity and the equations
\emph{(\ref{eq:G_jcompose with omega_j_})} and \emph{(\ref{eq:omega1+omega2})}
are satisfied. In particular, we have $H_{T_{1}}(\omega_{1}(0))=\omega_{2}(0)$
and $H_{T_{2}}(\omega_{2}(0))=\omega_{1}(0)$ whenever $\omega_{1}(0)\in\mathbb{C}^{+}$. 
\end{cor}

\begin{proof}
If $f,g:\mathbb{C}^{+}\to\mathbb{C}^{+}$ are analytic functions and
$f\circ g$ has a fixed point $\alpha\in\mathbb{C}^{+}$, then clearly
$g\circ f$ has the fixed point $g(\alpha).$ The first claim follows
from this observation applied to $f(w)=z+H_{T_{1}}(w)$ and $g(w)=z+H_{T_{2}}(w)$.
The second claim follows from the continuity of $\omega_{j}$ on $\overline{\mathbb{C}^{+}}$
for $j=1,2$. 
\end{proof}
We focus next on selfadjoint operators $T$ that belong to ${\rm Log}^{+}(\tau)$,
so the function 
\begin{align*}
L_{T}(z) & =\tau(\log(z-T))=\int_{\mathbb{R}}\log(z-t)\,d\mu_{T}(t),\quad z\in\mathbb{C}^{+},
\end{align*}
is defined and analytic. The principal value of the logarithm should
be used in these formulas, that is, $\Im\log(z-t)\in(0,\pi)$ for
$z\in\mathbb{C}^{+}$and $t\in\mathbb{R}$. Eventually, we will only
be interested in the logarithmic potential 
\[
\Re L_{T}(z)=\tau(\log|z-T|)=\int_{\mathbb{R}}\log|z-t|\,d\mu_{T}(t).
\]
Observe that the complex derivative of $L_{T}$ is 
\[
L_{T}'(z)=G_{T}(z),\quad z\in\mathbb{C}^{+}.
\]

\begin{lem}
\label{lem:logarithmic potential using omega}Let $T_{1},T_{2}\in{\rm Log}^{+}(\tau)$
be free selfadjoint operators, let $T=T_{1}+T_{2}$, and let $\omega_{1},\omega_{2}$
be the subordination functions provided by Theorem\emph{ \ref{thm:subordination}}.
Then we have 
\begin{equation}
L_{T_{1}}(\omega_{1}(z))+L_{T_{2}}(\omega_{2}(z))=L_{T}(z)+\log(\omega_{1}(z)+\omega_{2}(z)-z),\quad z\in\mathbb{C}^{+}.\label{eq:L_T in terms of L_Tj}
\end{equation}
\end{lem}

\begin{proof}
Using (\ref{eq:omega1+omega2}) we can rewrite the identity above
as 
\begin{equation}
L_{T_{1}}(\omega_{1}(z))+L_{T_{2}}(\omega_{2}(z))=L_{T}(z)-\log(G_{T}(z)).\label{eq:messenger}
\end{equation}
We show first that the difference between the two sides of (\ref{eq:messenger})
is a constant by verifying that the derivatives are the same: 
\[
G_{T_{1}}(\omega_{1}(z))\omega_{1}'(z)+G_{T_{2}}(\omega_{2}(z))\omega_{2}'(z)=G_{T}(z)-\frac{G_{T}'(z)}{G_{T}(z)}.
\]
Using (\ref{eq:G_jcompose with omega_j_}), this can be rewritten
as 
\[
\omega_{1}'(z)+\omega_{2}'(z)=1-\frac{G_{T}'(z)}{G_{T}(z)^{2}},
\]
which is recognized as the identity obtained by differentiating (\ref{eq:omega1+omega2}).
Thus, there exists a constant $c$ such that 
\[
L_{T_{1}}(\omega_{1}(z))+L_{T_{2}}(\omega_{2}(z))=L_{T}(z)-\log(G_{T}(z))+c,\quad z\in\mathbb{C}^{+}.
\]
To identify the constant $c$, we consider $z=iy$ with $y>0$ and
subtract $2\log(iy)$ from both sides to see that 
\begin{align*}
\tau\left(\log\frac{\omega_{1}(iy)-T_{1}}{iy}\right)+\tau\left(\log\frac{\omega_{2}(iy)-T_{2}}{iy}\right) & =\tau\left(\log\frac{iy-T}{iy}\right)-\log(iyG_{T}(iy))+c,\\
\end{align*}
for $y>0$. Now, let $y\to+\infty$ and note that each of the four
terms in this identity that involve $y$ tends to zero. For instance,
\[
\tau\left(\log\frac{iy-T}{iy}\right)=\int_{\mathbb{R}}\log\left(1-\frac{t}{iy}\right)\,d\mu_{T}(t),
\]
where the integrands are dominated by the integrable function $\log(1+|t|)$
for $y\ge1$, and tend to zero pointwise as $y\to+\infty$. For the other terms,
use the well-known limits 
\[
\lim_{y\uparrow\infty}\frac{\omega_{1}(iy)}{iy}=\lim_{y\uparrow\infty}\frac{\omega_{2}(iy)}{iy}=\lim_{y\uparrow\infty}iyG_{T}(iy)=1.\qedhere
\]
\end{proof}
For the remainder of this section, we restrict ourselves to selfadjoint
operators $T$ that have a symmetric distribution on $\mathbb{R}$.
This amounts to saying that $T$ and $-T$ have the same distribution,
or that $d\mu_{T}(t)=d\mu_{T}(-t)$. In terms of $G_{T}$, symmetry
is expressed by the identity 
\[
-\overline{G_{T}(-\overline{z})}=G_{T}(z),\quad z\in\mathbb{C}^{+},
\]
where the bar indicates complex conjugation. Symmetry also yields
\begin{equation}
G_{T}(z)=\frac{1}{2}[\tau((z-T)^{-1})+\tau((z+T)^{-1})]=z\tau[(z^{2}-T^{2})^{-1})],\quad z\in\mathbb{C}^{+}.\label{eq:G_T-symmetric}
\end{equation}
In particular, when $z=iy$ is purely imaginary, 
\begin{align*}
G_{T}(iy) & =\int_{\mathbb{R}}\frac{-iy\,d\mu_{T}(t)}{t^{2}+y^{2}},
\end{align*}
\begin{equation}
H_{T}(iy)=\frac{\tau(T^{2}(T^{2}+y^{2})^{-1})}{-iy\tau((T^{2}+y^{2})^{-1})}=\frac{{\displaystyle \int_{\mathbb{R}}\frac{t^{2}\,d\mu_{T}(t)}{t^{2}+y^{2}}}}{{\displaystyle \int_{\mathbb{R}}\frac{-iy\,d\mu_{T}(t)}{t^{2}+y^{2}}}},\label{eq:explicit HT}
\end{equation}
and 
\begin{align}
\Re L_{T}(iy) & =\frac{1}{2}[\tau(\log|iy-T|)+\tau(\log|iy+T|)]\nonumber \\
 & =\frac{1}{2}\tau(\log(T^{2}+y^{2}))=\frac{1}{2}\int_{\mathbb{R}}\log(t^{2}+y^{2})\,d\mu_{T}(t).\label{eq:log integral for symmetric T}
\end{align}
The symmetry property is inherited by the subordination functions. 
\begin{lem}
\label{lem:omegas are symmetric also}Suppose that $T_{1},T_{2}\in{\rm Log}^{+}(\tau)$
are selfadjoint, have symmetric distributions, and are free. Let $T=T_{1}+T_{2}$,
and let $\omega_{1},\omega_{2}$ be the subordination functions provided
by Theorem\emph{ \ref{thm:subordination}}. Then we have 
\[
-\overline{\omega_{j}(-\overline{z})}=\omega_{j}(z),\quad z\in\overline{\mathbb{C}^{+}},j=1,2.
\]
\end{lem}

\begin{proof}
Since $-T_{1}$ and $-T_{2}$ are free and have the same distributions
as $T_{1}$ and $T_{2}$, it follows that $-T=(-T_{1})+(-T_{2})$
also has the same distribution as $T$. If we set now 
\[
\check{\omega}_{j}(z)=-\overline{\omega_{j}(-\overline{z})},\quad z\in\overline{\mathbb{C}^{+}},j=1,2,
\]
we conclude that the equations (\ref{eq:G_jcompose with omega_j_})
and (\ref{eq:omega1+omega2}) are also satisfied with $\check{\omega}_{j}$
in place of $\omega_{j}$. The lemma follows from the uniqueness of
$\omega_{j}$. 
\end{proof}
Given $T\in\widetilde{\mathcal{A}}$ and $p\in\mathbb{R}$, we use
the notation 
\[
m_{p}(T)=\tau(|T|^{p})\in[0,+\infty]
\]
for the $p$-th absolute moment of $T$. If $T$ is not injective,
we set $m_{p}(T)=\infty$ for $p<0$. (Of course, it is possible that
$m_{p}(T)=\infty$ even when $\ker T=0$.) We use a similar notation
\[
m_{p}(\mu)=\int_{\mathbb{C}}|z|^{p}\,d\mu
\]
for an arbitrary probability measure $\mu$ on $\mathbb{C},$ and
observe that, when $T$ is normal, 
\[
m_{p}(T)=m_{p}(\mu_{T}),\quad p\in\mathbb{R}.
\]
Unless $m_{2}(\mu)=0,$ the Schwarz inequality shows that 
\begin{equation}
m_{-2}(\mu)m_{2}(\mu)\ge1,\label{eq:Shw}
\end{equation}
and equality occurs precisely when $\mu$ is supported on a circle
centered at the origin. When $\mu$ is a symmetric measure supported
on the real line, equality occurs when 
\[
\mu=\frac{1}{2}(\delta_{a}+\delta_{-a})
\]
for some $a\in(0,+\infty)$. In particular, the inequality (\ref{eq:Shw})
is strict if the support of such a measure contains at least three
points.

The following result is \cite[Lemma 4.8]{HA-SC-brown}. We provide
a somewhat more transparent argument below. 
\begin{lem}
\label{lem:increasing function}Suppose that $T\in\widetilde{A}$
is a symetric selfadjoint operator whose distribution is not of the
form $(\delta_{a}+\delta_{-a})/2$ for any $a\in[0,+\infty)$. Then
the function 
\[
y\mapsto-iyH_{T}(iy),\quad y\in(0,+\infty),
\]
is strictly increasing and its range is the open interval $(1/m_{-2}(T),m_{2}(T))$.
If $\mu=(\delta_{a}+\delta_{-a})/2$, we have $-iyH_{T}(iy)=a^2$ for
every $y>0$. 
\end{lem}

\begin{proof}
The limits at $0$ and at $+\infty$ are easily calculated from (\ref{eq:explicit HT}),
so we only verify that the function is increasing. The function $H_{T}$
maps $\mathbb{C}^{+}$ to itself, 
\[
\lim_{y\uparrow\infty}H_{T}(iy)/iy=0,
\]
and 
\[
-\overline{H_{T}(-\overline{z})}=H_{T}(z),\quad z\in\mathbb{C}^{+}.
\]
It follows that $H_{T}$ has a Nevanlinna representation \cite{akh}
of the form 
\begin{equation}\label{eqn:Nevanlinna}
 H_{T}(z)=\int_{\mathbb{R}}\frac{1+tz}{t-z}\,d\rho(t)
\end{equation}
for some positive, finite, symmetric measure $\rho$ on $\mathbb{R}$
with support different from $\{0\}$. We can then calculate 
\[
-iyH_{T}(iy)=\int_{\mathbb{R}}\frac{(1+t^{2})y^{2}}{y^{2}+t^{2}}\,d\rho(t),\quad y>0.
\]
Clearly, the integrand is a strictly increasing function of $y$ at
points $t\ne0$. The lemma follows. 
\end{proof}

\begin{rem}\label{rem:about-rho}
Suppose that $T\in\widetilde{A}$
is a symetric selfadjoint operator whose distribution is not of the
form $(\delta_{a}+\delta_{-a})/2$ for any $a\in[0,+\infty)$. Represent the function $H_{T}$ as in \eqref{eqn:Nevanlinna}. Then, Lemma \ref{lem:increasing function} also implies 
\begin{enumerate}
  \item $\rho(\{0\})=\lim_{y\rightarrow {0}}-iyH_T(iy)=1/m_{-2}(T)$, and 
  \item $\int_{\mathbb{R}}(1+t^2)d\rho(t)=\lim_{y\rightarrow {\infty}}-iyH_T(iy)=m_2(T)$.
\end{enumerate}
\end{rem}

\begin{prop}\label{prop:negative-derivative-HT}
Under the assumption of Lemma {\ref{lem:increasing function}}, we have 
\[
 y^2H_T'(iy)\leq \frac{m_2(T)-1/m_{-2}(T)}{2}+iyH_T(iy), \qquad \forall y>0.
\]
In particular, if
$m_2(T) m_{-2}(T)\leq 3$, then $H_T'(iy)\leq 0$ for any $y>0$. 
\end{prop}
\begin{proof}
For every $z\in\mathbb{C}^+$, 
\[
  H_T'(z)=\int_{\mathbb{R}}\frac{1+t^2}{(t-z)^2}d\rho(t).
\]
Hence, for $z=iy$, we have 
\begin{equation}\label{eqn:deritative-HT}
    H_T'(iy)=\int_{\mathbb{R}}\frac{1+t^2}{(t-iy)^2}d\rho(t)=\int_\mathbb{R}\frac{t^2-y^2}{(t^2+y^2)^2}(1+t^2)d\rho(t),
\end{equation}
because $\rho$ is a symmetric measure. 
Using Remark \ref{rem:about-rho}, we see that 
\begin{align*}
  y^2H_T'(iy)&=\int_{\mathbb{R}}\frac{y^2(t^2-y^2)}{(t^2+y^2)^2}(1+t^2)d\rho(t)\\
    &=2\int_{\mathbb{R}}\frac{y^2t^2}{(t^2+y^2)^2}(1+t^2)d\rho(t)-\int_{\mathbb{R}}\frac{y^2}{t^2+y^2}(1+t^2)d\rho(t)\\
    &\leq \frac{1}{2}\int_{\mathbb{R}\backslash\{0\}}(1+t^2)d\rho(t)+iyH_T(iy)\\
    &=\frac{m_2(T)-1/m_{-2}(T)}{2}+iyH_T(iy). 
\end{align*}
Since $iyH_T(iy)\leq -1/m_{-2}(T)$, it follows that 
\[
   y^2H_T'(iy)\leq \frac{m_2(T)-3/m_{-2}(T)}{2}\leq 0
\]
provided that $m_2(T)m_{-2}(T)\leq 3$. 
\end{proof}

The following result is needed in Section \ref{sec:The-calculation-of}. 
\begin{prop}
\label{prop:omega(0)}Let $T_{1},T_{2}\in\widetilde{\mathcal{A}}$
be nonzero, free, selfadjoint elements that have symmetric distributions,
and let $\omega_{1},\omega_{2}$ be the subordination functions provided
by Theorem \emph{\ref{thm:subordination}}. 
\begin{enumerate}
\item If $\omega_{1}(0)\in\mathbb{C}^{+}$, then $\omega_{2}(0)\in\mathbb{C}^{+}$.
When this occurs, the closed intervals $[1/m_{-2}(T_{1}),m_{2}(T_{1})]$
and $[1/m_{-2}(T_{2}),m_{2}(T_{2})]$ have a common point. Unless
$1/m_{-2}(T_{1})=m_{2}(T_{1})$, the common point is in $(1/m_{-2}(T_{1}),m_{2}(T_{1}))$. 
\item We cannot have $\omega_{1}(0)=\omega_{2}(0)=\infty$. 
\item We have $\omega_{1}(0)=\omega_{2}(0)=0$ if and only if $\mu_{T_{1}}(\{0\})+\mu_{T_{2}}(\{0\})\ge1$.
When this occurs, we have $G_{T_{1}+T_{2}}(0)=\infty$. 
\item If $\omega_{1}(0)=0$ and $\omega_{2}(0)=\infty$, we have 
\[
m_{2}(T_{2})\le\frac{1}{m_{-2}(T_{1})}.
\]
In particular, $m_{2}(T_{2})<+\infty$ and $m_{-2}(T_{1})<+\infty$. 
\end{enumerate}
\end{prop}

\begin{proof}
The first assertion in (1), and the equality $\omega_{2}(0)=H_{T_{1}}(\omega_{1}(0))$,
follow from Corollary \ref{cor:when omega1(0) is in H} and its proof.
For the second assertion, we multiply this equality by $-\omega_{1}(0)$
and apply Lemma \ref{lem:increasing function} to conclude that 
\begin{align*}
-\omega_{1}(0)\omega_{2}(0) & =-\omega_{1}(0)H_{T_{1}}(\omega_{1}(0))\\
 & =-\omega_{2}(0)H_{T_{2}}(\omega_{2}(0))
\end{align*}
belongs to the the intervals $[1/m_{-2}(T_{j}),m_{2}(T_{j})]$ since
$\omega_{j}(0)$ is purely imaginary.

We prove (2) by contradiction. Suppose that $\omega_{1}(0)=\omega_{2}(0)=\infty$,
and rewrite (\ref{eq:omega1+omega2}) as 
\[
\frac{\omega_{2}(iy)}{\omega_{1}(iy)}=\frac{iy}{\omega_{1}(iy)}-1+\frac{1}{\omega_{1}(iy)G_{T_{1}}(\omega_{1}(iy))},\quad y>0,
\]
and recall the fact that 
\[
\lim_{t\uparrow\infty}\frac{1}{itG_{T_{1}}(it)}=1.
\]
Writing $\omega_1(iy)=i\eta_1(iy)$, the assumping that $\omega_1(0)=\infty$ amounts to
\[
\lim_{y\to \infty}\eta_1(iy)=+\infty. 
\]
and 
letting $y\downarrow0$ in the above inequality, we see that $\lim_{y\downarrow0}\omega_{2}(iy)/\omega_{1}(iy)=0$.
By symmetry, we also have $\lim_{y\downarrow0}\omega_{1}(iy)/\omega_{2}(iy)=0$,
and these last two equalities cannot hold simultaneously.

The equality $\omega_{1}(0)=\omega_{2}(0)=0$ means that $0$ is the
Denjoy-Wolff point of the maps $H_{T_{1}}\circ H_{T_{2}}$ and $H_{T_{2}}\circ H_{T_{1}}.$
Since these maps preserve the positive imaginary line, we
deduce that $H_{T_{j}}$ extends to $0$ such that $H_{T_{j}}(0)=0$.
Moreover, the Julia-Carath\'{e}odory derivative of $H_{T_{j}}$ at zero
is 
\[
\lim_{y\downarrow0}\frac{H_{T_{j}}(iy)}{iy}=\frac{1}{\mu_{j}(\{0\})}-1.
\]
Thus, applying the chain rule for functions of one real variable,
we see that $0$ is the Denjoy-Wolff point of $H_{T_{1}}\circ H_{T_{2}}$
precisely when 
\[
\left(\frac{1}{\mu_{1}(\{0\})}-1\right)\left(\frac{1}{\mu_{2}(\{0\})}-1\right)\le1,
\]
and this is equivalent to the inequality stated in (3). The fact that
$G_{T_{1}+T_{2}}(0)=\infty$ follows from (\ref{eq:omega1+omega2}).

Finally, suppose that $\omega_{1}(0)=0$ and $\omega_{2}(0)=\infty$.
Multiplying (\ref{eq:omega1+omega2}) by $\omega_{1}$ or $\omega_{2},$
we obtain the identities 
\begin{align*}
\omega_{1}(iy)\omega_{2}(iy) & =iy\omega_{1}(iy)+\omega_{1}(iy)H_{T_{1}}(\omega_{1}(iy))\\
 & =iy\omega_{2}(iy)+\omega_{2}(iy)H_{T_{2}}(\omega_{2}(iy)),\quad y>0.
\end{align*}
Observe that $iy\omega_{2}(iy)\le0$ to conclude that 
\[
iy\omega_{1}(iy)+\omega_{1}(iy)H_{T_{1}}(\omega_{1}(iy))\le\omega_{2}(iy)H_{T_{2}}(\omega_{2}(iy)),
\]
and now let $y\to0$ and use Lemma \ref{lem:increasing function}
to conclude that 
\[
-1/m_{-2}(T_{1})\le-m_{2}(T_{2}),
\]
thus establishing (4). 
\end{proof}
\begin{rem}
\label{rem:special cases T_1 T_2}In the preceding result, it is useful
to separate the cases in which one or both of the measures $\mu_{T_{j}}$
are of the form $(\delta_{a}+\delta_{-a})/2$ for some $a>0$. Denote
by $J_{T}$ or $J_{\mu_{T}}$ the interval $[1/m_{-2}(T),m_{2}(T)]$. 
\end{rem}

\begin{enumerate}
\item Suppose that the intervals $J_{T_{1}}$ and $J_{T_{2}}$ have interior
points. Then one of the following mutually exclusive possibilities
arises: 
\begin{enumerate}
\item $\mu_{T_{1}}(\{0\})+\mu_{T_{2}}(\{0\})\ge1$, in which case $\omega_{1}(0)=\omega_{2}(0)=0$. 
\item $\mu_{T_{1}}(\{0\})+\mu_{T_{2}}(\{0\})<1$ but $J_{T_{1}}$ and $J_{T_{2}}$
have common interior points. In this case, $\omega_{1}(0)$ and $\omega_{2}(0)$
are finite and different from zero. 
\item $J_{T_{1}}$ is on the right of $J_{T_{2}},$ that is, $0<m_{2}(T_{2})\le1/m_{-2}(T_{1})<+\infty$,
in which case $\omega_{1}(0)=0$ and $\omega_{2}(0)=\infty$. 
\item $0<m_{2}(T_{1})\le1/m_{-2}(T_{2})<+\infty$, in which case $\omega_{1}(0)=\infty$
and $\omega_{2}(0)=0$. 
\end{enumerate}
\item Suppose that $J_{T_{1}}$ has interior points but $J_{T_{2}}$ does
not. In this case, $J_{T_{2}}=\{a_{2}^2\}$, $\mu_{T_{2}}=(\delta_{a_{2}}+\delta_{-a_{2}})/2$,
and $a_{2}>0$ because $T_{2}$ was assumed to be nonzero. Then there
are three mutually exclusive possibilities: 
\begin{enumerate}
\item $a_{2}^2$ is an interior point of $J_{T_{1}},$ in which case $\omega_{1}(0)$
and $\omega_{2}(0)$ are finite and nonzero. 
\item $J_{T_{1}}$ is on the right of $a_{2}^2,$ that is, $0<m_{2}(T_{2})=a_{2}^2\le1/m_{-2}(T_{1})<+\infty$,
in which case $\omega_{1}(0)=0$ and $\omega_{2}(0)=\infty$. 
\item $0<m_{2}(T_{1})\le a_{2}^2=1/m_{-2}(T_{2})<+\infty$, in which case
$\omega_{1}(0)=\infty$ and $\omega_{2}(0)=0$. 
\item A symmetric statement holds if $J_{T_{1}}$ is a singleton and
$J_{T_{2}}$ has interior points. 
\end{enumerate}
\item Neither $J_{T_{1}}$ nor $J_{T_{2}}$ has interior points, so $J_{T_{1}}=\{a_{1}^2\}$
and $J_{T_{2}}=\{a_{2}^2\}$ with $a_{1},a_{2}>0$. The alternatives,
which are also easily verified by explicit calculation, are as follows: 
\begin{enumerate}
\item $a_{1}=a_{2}$, in which case $\omega_{1}(0)=\omega_{2}(0)=ia_{1}$. 
\item $a_{1}>a_{2},$ in which case $\omega_{1}(0)=0$ and $\omega_{2}(0)=\infty$. 
\item $a_{1}<a_{2}$, in which case $\omega_{1}(0)=\infty$ and $\omega_{2}(0)=0$. 
\end{enumerate}
\end{enumerate}
\begin{rem}
\label{rem:JC derivatives of omega at 0 in (0,infty) case} Suppose
$\omega_{1}(0)=0$ and $\omega_{2}(0)=\infty$. The last calculation
in the proof of Proposition \ref{prop:omega(0)} also implies the
identity 
\[
\lim_{y\downarrow0}(iy\omega_{2}(iy))=m_{2}(T_{2})-\frac{1}{m_{-2}(T_{1})},
\]
or, equivalently, 
\[
\lim_{y\downarrow0}\frac{1}{-iy\omega_{2}(iy)}=\frac{m_{-2}(T_{1})}{1-m_{2}(T_{2})m_{-2}(T_{1})}.
\]
This is the Julia-Carath\'{e}odory derivative of $-1/\omega_{2}$ at zero.
The Julia-Carath\'{e}odory derivative of $\omega_{1}$ at zero can also
be calculated by noting that 
\begin{align*}
\frac{\omega_{1}(iy)}{iy} & =1+\frac{H_{T_{2}}(iy)}{iy}\\
 & =1+\frac{1}{iy\omega_{2}(iy)}\omega_{2}(iy)H_{T_{2}}(\omega_{2}(iy)).
\end{align*}
Let now $y\to0$ and recall that $\omega_{2}(iy)\to\infty$ to obtain
\[
\lim_{y\downarrow0}\frac{\omega_{1}(iy)}{iy}=1+\frac{m_{2}(T_{2})m_{-2}(T_{1})}{1-m_{2}(T_{2})m_{-2}(T_{1})}=\frac{1}{1-m_{2}(T_{2})m_{-2}(T_{1})},
\]
where we use once again Lemma \ref{lem:increasing function}. 
\end{rem}

\section{The calculation of $\mu_{X}$\label{sec:The-calculation-of}}

With $(\mathcal{A},\tau)$ as before, we consider $*$-free operators
$X_{1},X_{2}\in{\rm Log}^{+}(\tau)$ such that $X_{2}$ is $R$-diagonal.
We exclude the trivial case in which either $X_{1}$ or $X_{2}$ is
a scalar multiple of the identity. We also assume that there exists
a Haar unitary operator $U\in\mathcal{A}$ that is $*$-free from
$\{X_{1},X_{2}\}$. This is not a true restriction as we can always
enlarge the algebra $\mathcal{A}$ without affecting the measure $\mu_{X}$,
where $X=X_{1}+X_{2}$. For every $\lambda\in\mathbb{C}$, the selfadjoint
operators 
\[
|X-\lambda|=|(X_{1}-\lambda)+X_{2}|
\]
and 
\[
T^{(\lambda)}=|X_{1}-\lambda+U^{*}X_{2}|=|U(X_{1}-\lambda)+X_{2}|
\]
have the same ditribution. Now, $U(X_{1}-\lambda)$ is $R$-diagonal
and so is $X_2$ by hypothesis, according to \cite{Ni-Sp-waterloo,ha-larsen,HA-SC-brown},
\[
\widetilde{\mu}_{T^{(\lambda)}}=\mu_{1}^{(\lambda)}\boxplus\mu_{2},
\]
where 
\begin{equation}
\mu_{1}^{(\lambda)}=\widetilde{\mu}_{|X_{1}-\lambda|}\text{ and }\mu_{2}=\widetilde{\mu}_{|X_{2}|}.\label{eq:definition of mu1lambda and mu2}
\end{equation}
Denote by $\omega_{1}^{(\lambda)}$ and $\omega_{2}^{(\lambda)}$
the subordination functions associated with this free convolution.
Thus, 
\[
G_{\mu_{1}^{(\lambda)}}(\omega_{1}^{(\lambda)}(z))=G_{\mu_{2}}(\omega_{2}^{(\lambda)}(z))=G_{\widetilde{\mu}_{|X-\lambda|}}(z),
\]
and 
\begin{equation}
  \label{eqn:subordination-lambda-index}
    \omega_{1}^{(\lambda)}(z)+\omega_{2}^{(\lambda)}(z)=z+\frac{1}{G_{\widetilde{\mu}_{|X-\lambda|}}(z)},\quad z\in\mathbb{C}^{+}.
\end{equation}

We now introduce several subsets of the complex plane $\mathbb{C}$
that are related with the the values of $\omega_{1}^{(\lambda)}(0)$
and $\omega_{2}^{(\lambda)}(0)$. 
\begin{defn}
\label{def:DSF1F2F} Let $X_{1},X_{2}\in{\rm Log}^{+}(\tau)\backslash\mathbb{C}$
be free, let $X=X_{1}+X_{2}$, suppose that $X_{2}$ is $R$-diagonal,
and let $\omega_{1}^{(\lambda)},\omega_{2}^{(\lambda)}$ be the subordonation
functions arising from the free convolution $\widetilde{\mu}_{|X-\lambda|}=\widetilde{\mu}_{|X_{1}-\lambda|}\boxplus\widetilde{\mu}_{|X_{2}|}$.
We define subsets of $\mathbb{C}$ as follows: 
\begin{align*}
S & =\{\lambda:\omega_{1}^{(\lambda)}(0)=\omega_{2}^{(\lambda)}(0)=0\},\\
F_{1} & =\left\{ \lambda:m_{2}(|X_{2}|)\le\frac{1}{m_{-2}(|X_{1}-\lambda|)}\right\} ,\\
F_{2} & =\left\{ \lambda:m_{2}(|X_{1}-\lambda|)\le\frac{1}{m_{-2}(|X_{2}|)}\right\} ,\\
F & =F_{1}\cap F_{2},\quad\text{ and }\\
\Omega & =\mathbb{C}\backslash(S\cup F_{1}\cup F_{2}).
\end{align*}
\end{defn}

Thus, the sets $S,F_{1}\backslash F,F_{2}\backslash F,F,$ and $\Omega$
form a partition of $\mathbb{C}$. Remark \ref{rem:special cases T_1 T_2}
allows for a different description of these sets in terms of the values
of $\omega_{1}^{(\lambda)}(0)$ and $\omega_{2}^{(\lambda)}(0)$.
For the following statement, $\ker(Y)\in\mathcal{A}$ denotes the
orthogonal projection onto the null space of $Y$. Also, note that
$\widetilde{\mu}_{|Y|}$ is of the form $\text{\ensuremath{(\delta_{a}+\delta_{-a})/2}, \ensuremath{a>0},}$
precisely when $|Y|=a$ or, equivalently, when $Y/a$ is unitary.
If $Y$ is also $R$-diagonal, the unitary operator $Y/a$ is a Haar
unitary. 
\begin{lem}
\label{lem:desciption of four parts}Suppose that $X_{1},X_{2}\in{\rm Log}^{+}(\tau)\backslash\mathbb{C}$
are $*$-free and $X_{2}$ is $R$-diagonal. 
\begin{enumerate}
\item We have $S=\{\lambda:\tau(\ker(X_{1}-\lambda))+\tau(\ker(X_{2}))\ge1\}.$
In particular, $S$ is a finite set, disjoint from $F_{1}\cup F_{2}$.
If $S$ is not empty then $F_{2}=\varnothing$. 
\item 
$F$ consists of those $\lambda\in\mathbb{C}$ such that $|X_{1}-\lambda|=|X_{2}|=\|X_{2}\|$. In particular, $\lambda\in F$ if and only if $H_{\widetilde\mu_{|X_2|}}\circ H_{\widetilde\mu_{|X_1-\lambda|}}$ is the identity map on $\mathbb{C}^+$. Moreover, $F$ contains at most two points and $F$ is contained in the boundary of $\Omega$. For $\lambda\in F$, we have $\omega_{1}^{(\lambda)}(0)\in\mathbb{C}^{+}$ and $\omega_{2}^{(\lambda)}(0)\in\mathbb{C}^{+}$.

\item The set $F_{1}\backslash F$ consists of those $\lambda\in\mathbb{C}$
for which $\omega_{1}^{(\lambda)}(0)=0$ and $\omega_{2}^{(\lambda)}(0)=\infty$.
In particular, $F_{1}$ is empty if $m_{2}(|X_{2}|)=+\infty$. If
$F_{1}\ne\varnothing$ then it is a closed set. 
\item The set $F_{2}\backslash F$ consists of those $\lambda\in\mathbb{C}$
for which $\omega_{1}^{(\lambda)}(0)=\infty$ and $\omega_{2}^{(\lambda)}(0)=0$.
In particular, $F_{2}$ is empty if $m_{2}(|X_{1}|)=+\infty$. In
general, $F_{2}$ is $\varnothing$, $\{\tau(X_{1})\}$, or a closed
disk centered at $\tau(X_{1})$. 
\item We have $\lambda\in\Omega$ if and only if $\lambda\notin F$ and
neither $\omega_{1}^{(\lambda)}(0)$ nor $\omega_{2}^{(\lambda)}(0)$
belongs to $\{0,\infty\}$. The sets $\Omega$ and $\Omega\cup S$
are open. 
\end{enumerate}
\end{lem}

\begin{proof}
The first assertion in (1) follows from Proposition \ref{prop:omega(0)}(3).
If $S$ is not empty, then $\ker(X_{2})\ne0$. This implies $m_{-2}(|X_{2}|)=+\infty$,
and hence $m_{2}(|X_{1}-\lambda|)=0$ for $\lambda\in F_{2}$. Since
$X_{1}$ is not a scalar multiple of the identity, there is no such
$\lambda$, in other words, $F_{2}=\varnothing$.

Suppose now that $F\ne\varnothing$ and $\lambda_{0}\in F$. The definition
of the sets $F_{1}$and $F_{2}$, combined with (\ref{eq:Shw}), shows
that
\begin{align*}
\frac{1}{m_{-2}(|X_{2}|)} & \le m_{2}(|X_{2}|)\le\frac{1}{m_{-2}(|X_{1}-\lambda|)},\text{ and}\\
\frac{1}{m_{-2}(|X_{1}-\lambda|)} & \le m_{2}(|X_{1}-\lambda|)\le\frac{1}{m_{-2}(|X_{2}|)}.
\end{align*}
These inequalities imply
\begin{equation}
m_{2}(|X_{2}|)=\frac{1}{m_{-2}(|X_{2}|)}=m_{2}(|X_{1}-\lambda|)=\frac{1}{m_{-2}(|X_{1}-\lambda|)},\label{eq:some identity}
\end{equation}
and thus $|X_{2}|=|X_{1}-\lambda|$ is a scalar multiple of the identity. Thus, $F$ contains at most two points. 

If $H_{\widetilde\mu_{|X_2|}}\circ H_{\widetilde\mu_{|X_1-\lambda|}}$ is the identity map, the Schwarz--Pick inequality implies that both maps are automorphisms of $\CC^+$. For a symmetric probability measure $\nu$, an automorphism $H_\nu$ must have the form $H_\nu(z)=-b^2/z$, $b>0$. Indeed, it preserves the positive imaginary axis and satisfies $H_\nu(\iu y)/(\iu y)\to0$ as $y\to\infty$. This is equivalent to $\widetilde\mu_{|X_1-\lambda|}=(\delta_{b_1}+\delta_{-b_1})/2$ and $\widetilde\mu_{|X_2|}=(\delta_{b_2}+\delta_{-b_2})/2$. For these laws, the composition is $(b_2^2/b_1^2)z$, which is the identity precisely when $b_1=b_2$. This proves first two assertions in (2).
To see that $F$ is contained in the boundary of $\Omega$, suppose
that $\lambda_{0}\in F$, so $|X_{1}-\lambda_{0}|=|X_{2}|=\|X_{2}\|$.
Thus, $\ker X_{2}=0$ and so $S$ is empty. Also, $X_{1}-\lambda_{0}=aU,$where
$a=\|X_{2}\|$ and $U$ is unitary. For $\lambda$ close to $\lambda_{0}$,
we have 
\begin{align*}
m_{-2}(|X_1-\lambda|)=\int_{|\zeta|=1}\frac{d\mu_U(\zeta)}{|a\zeta+\lambda_0-\lambda|^2}.
\end{align*}
Near $\lambda_0$, this function of $\lambda$ and $\lambda\mapsto m_2(|X_1-\lambda|)$ are continuous and have strictly positive Laplacians. In particular, these functions do not have any local extreme near $\lambda_0$. Since $m_{-2}(|X_{1}-\lambda_0|)=a^{-2}$ and $ m_2(|X_1-\lambda_0|)=a^2$, neither
$F_1=\{m_{-2}(|X_1-\lambda|)\le a^{-2}\}$ nor $F_2=\{ m_2(|X_1-\lambda_0|)\le a^2\}$ contains a neighborhood of $\lambda_0$. This completes the proof of (2).

The first two assertions in (3) follow from Remark \ref{rem:special cases T_1 T_2}.
To show that $F_{1}$ is closed, it suffices to verify that the function
$\lambda\mapsto m_{-2}(|X_{1}-\lambda|)$ is lower semicontinuous.
If this function is not identically $+\infty$, lower semicontinuity
follows because 
\[
m_{-2}(|X_{1}-\lambda|)=\sup_{\varepsilon>0}\tau((|X_{1}-\lambda|^{2}+\varepsilon)^{-1}),
\]
and the function $\lambda\mapsto\tau((|X_{1}-\lambda|^{2}+\varepsilon)^{-1})$
is continuous for every $\varepsilon>0$.

The first two assertions in (4) also follow from Remark \ref{rem:special cases T_1 T_2}.
The fact that $F_{2}$ is a closed disk follows because 
\begin{align*}
m_{2}(|X_{1}-\lambda|) & =\|X_{1}-\lambda\|_{2}^{2}\\
 & =\|(X_{1}-\tau(X_{1}))-(\lambda-\tau(X_{1}))\|_{2}^{2}\\
 & =\|X_{1}-\tau(X_{1})\|_{2}^{2}+|\lambda-\tau(X_{1})|^{2}.
\end{align*}
Thus, provided $m_{2}(X_{2})<+\infty$, $F_{2}$ is nonempty precisely
when $m_{2}(|X_{1}-\tau(X_{1})|)\le1/m_{-2}(X_{2})$, in which case
$F_{2}$ is a closed disk (possibly of zero radius).

Finally, (5) follows from Remark \ref{rem:special cases T_1 T_2}
and from the fact that $F_{1}\cup F_{2}$ and $F_{1}\cup F_{2}\cup S$
are closed. 
\end{proof}

\begin{cor}
\label{cor:omega-boundary}
For $\varepsilon>0$, write $\omega_j^{(\lambda)}(\iu\varepsilon)=\iu\eta_j(\lambda,\varepsilon)$, $j=1,2$. Then, as $(\lambda,\varepsilon)\to(\lambda_0,0)$ with $\varepsilon>0$,
\begin{align*}
(\eta_1(\lambda,\varepsilon),\eta_2(\lambda,\varepsilon))&\to(0,\infty),\qquad \lambda_0\in F_1\setminus F,\\
(\eta_1(\lambda,\varepsilon),\eta_2(\lambda,\varepsilon))&\to(\infty,0),\qquad \lambda_0\in F_2\setminus F.
\end{align*}
The convergence is uniform as $\dist(\lambda,K)+\varepsilon\to0$, for every compact subset $K$ of $F_1\setminus F$ or $F_2\setminus F$, respectively.
\end{cor}

\begin{proof}
Set $\nu_\lambda=\widetilde\mu_{|X_1-\lambda|}$ and $\nu_2=\widetilde\mu_{|X_2|}$. Suppose that $(\lambda_n,\varepsilon_n)\to(\lambda_0,0)$ with $\varepsilon_n>0$. Observe that $H_{\nu_{\lambda_n}}\to H_{\nu_{\lambda_0}}$ locally uniformly on $\CC^+$. Hence
\begin{align*}
\iu\varepsilon_n+H_{\nu_2}\bigl(\iu\varepsilon_n+H_{\nu_{\lambda_n}}(z)\bigr)\to H_{\nu_2}(H_{\nu_{\lambda_0}}(z))
\end{align*}
locally uniformly on $\CC^+$. If $\lambda_0\in(F_1\cup F_2)\setminus F$, Lemma~\ref{lem:desciption of four parts} shows that the limiting map is not the identity. Theorem~\ref{thm:dw-continuity} therefore implies that $\omega_1^{(\lambda_n)}(\iu\varepsilon_n)\to\omega_1^{(\lambda_0)}(0)$. Interchanging the two measures gives the corresponding convergence for $\omega_2$. The description of $F_1$ and $F_2$ in Lemma~\ref{lem:desciption of four parts} gives the first two assertions. The uniform convergence near compact subsets follows from the preceding statement and compactness.
\end{proof}

\begin{rem}
\label{rem:two points in F}The argument above shows that $F$ is
empty unless $X_{2}$ is a scalar multiple of a Haar unitary, and
$X_{1}$ is normal with spectrum contained in a circle of radius $\|X_{2}\|$.
Consider the case in which $X_{2}$ is a Haar unitary and $X_{1}$
is normal. Then $F$ consists of the centers of those circles of radius
one which contain the spectrum of $X_{1}$. It the spectrum of $X_{1}$
contains three or more points, or if it consists of two points at
distance $2$, there is exactly one such point. We illustrate in Figure
\ref{fig:one point in F} the support of the Brown measure when the
$X_{1}$ is distributed like the $3\times3$ diagonal matrix with
spectrum $\{1,i,-1\}$. In this case, $F$ is a singleton and $F_{2}$
is a disk. In Figure \ref{fig:lemniscate}, $X_{1}$ has spectrum
$\{1,-1\},$so $F=\{0\}$ in both cases. For the first image, $X_{1}$
is distributed like the $3\times3$ diagonal matrix with entries $1,1,-1$
and $F_{2}$ is a disk. For the second one (also calculated and illustrated
in \cite[Figure 1]{biane-lehner}), $X_{1}$ is distributed like the
$2\times2$ diagonal matrix with entries $1,-1$, and $F=F_{2}=\{0\}$.
If $X_{1}$ is normal and its spectrum consists of just two points
with distance in $(0,2)$, the set $F$ contains two points. $\lambda_{2}\ne\lambda_{1}$
and the center of the disk $F_{2}$ is $(\lambda_{1}+\lambda_{2})/2$.
Figure \ref{fig:F has two points} illustrates the Brown measure when
$X_{1}$ is distributed like the $2\times2$ diagonal matrix with
spectrum $\{0,\sqrt{2}\}$. The pictures are obtained from simulations
of random matrices that approximate $X_{1}+X_{2}$. If $X_{2}$ is
Haar unitary and $X_{1}$ is unitary, we always have $0\in F$, and
$\tau(X_{1})$ is the center of $F_{2}$; see also \cite[Figure 6]{biane-lehner}
for another illustration.

\begin{figure}
\includegraphics[scale=0.4]{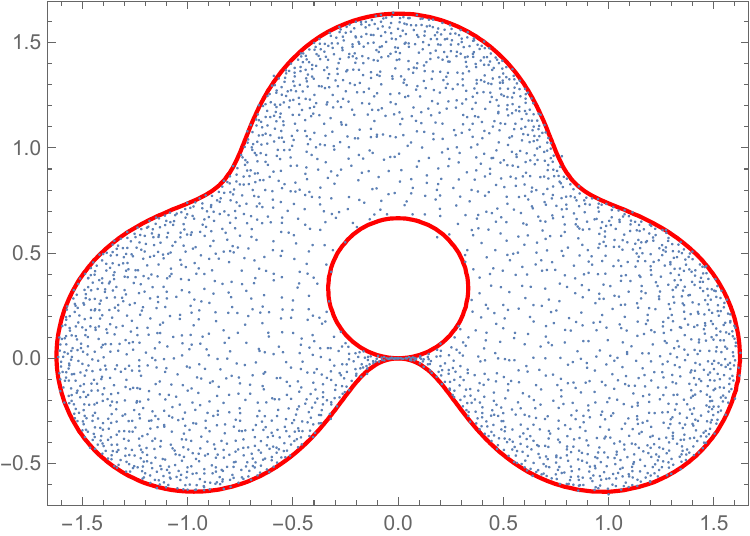}\caption{\label{fig:one point in F}$F=\{0\}$ }
\end{figure}

\begin{figure}
\includegraphics[scale=0.4]{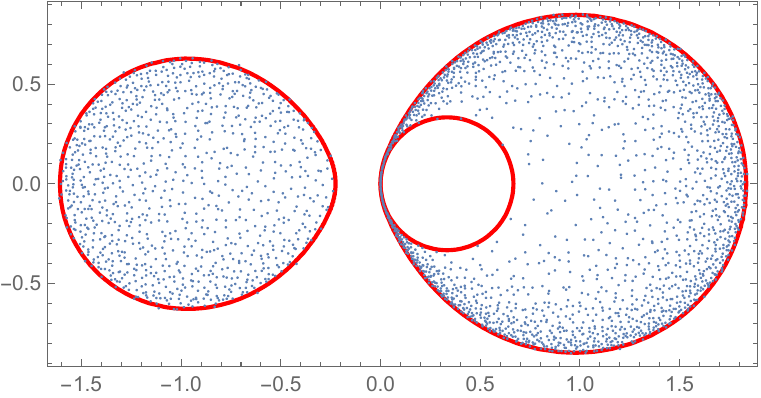}\includegraphics[scale=0.4]{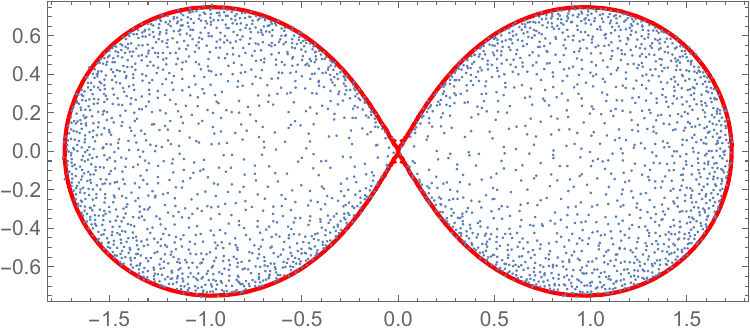}

\caption{\label{fig:lemniscate}$X_{1}$ unitary with spectrum $\{\pm1\}$}

\end{figure}

\begin{figure}
\includegraphics[scale=0.4]{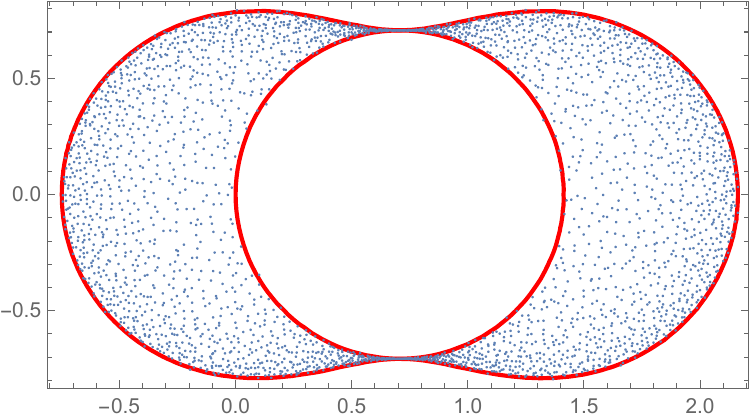}

\caption{\label{fig:F has two points}$F$ contains the points $\sqrt{2},\pm i\sqrt{2})$
at distance $1$ from $0$ and $\sqrt{2}$}

\end{figure}
\end{rem}

We can now calculate the logarithmic integral of $\mu_{X}$ in terms
of $\mu_{X_{1}},\mu_{X_{2}},$ and the values of $\omega_{j}^{(\lambda)}(0)$
for $\lambda\in\Omega$.
\begin{prop}
\label{prop:calculating the logarithmic integral}Suppose that $\lambda\in\mathbb{C}\backslash S$.
Then\emph{:} 
\begin{enumerate}
\item If $F_{1}$ is not empty, we have $m_{-2}(|X_{1}-\lambda|)<+\infty$
and $\tau(\log(|X-\lambda|)=\tau(\log|X_{1}-\lambda|)$ for every
$\lambda\in F_{1}$. 
\item If $F_{2}$ is not empty, we have $m_{-2}(|X_{2}|)<+\infty$, and
$\tau(\log(|X-\lambda|)=\tau(\log|X_{2}|)$ for every $\lambda\in F_{2}$. 
\item If $\lambda\in\Omega$, we have 
\begin{align*}
\tau(\log(|X-\lambda|) & =\frac{1}{2}\tau(\log(|X_{1}-\lambda|^{2}-\omega_{1}^{(\lambda)}(0)^{2})+\frac{1}{2}\tau(\log(|X_{2}|^{2}-\omega_{2}^{(\lambda)}(0)^{2})\\
 & -\log(-i\omega_{1}^{(\lambda)}(0)-i\omega_{2}^{(\lambda)}(0)).
\end{align*}
\end{enumerate}
\end{prop}

\begin{proof}
Taking real parts in (\ref{eq:L_T in terms of L_Tj}) and setting $z=iy$ for $y>0$,
we obtain 
\begin{align*}
\int_{\mathbb{C}}\log|t+iy|\,d\mu_{T}(t) & =\int_{\mathbb{C}}\log|t+\omega_{1}^{(\lambda)}(iy)|\,d\mu_{T_{1}}(t)L_{T_{1}}+\int_{\mathbb{C}}\log|t+\omega_{2}^{(\lambda)}(iy)|\,d\mu_{T_{2}}(t)\\
 & -\log|\omega_{1}^{(\lambda)}(z)+\omega_{2}^{(\lambda)}(z)-z|.
\end{align*}
In the case of symmetric measures, we can use (\ref{eq:log integral for symmetric T})
to rewrite this as 
\begin{align*}
\frac{1}{2}\int_{\mathbb{R}}\log(t^{2}+y^{2})\,d\mu_{T}(t) & =\frac{1}{2}\int_{\mathbb{R}}\log(t^{2}-\omega_{1}^{(\lambda)}(iy)^{2})\,d\mu_{T_{1}}(t)\\
 & +\frac{1}{2}\int_{\mathbb{R}}\log(t^{2}-\omega_{2}^{(\lambda)}(iy)^{2})\,d\mu_{T_{2}}(t)\\
 & -\log(-i\omega_{1}^{(\lambda)}(iy)-i\omega_{2}^{(\lambda)}(iy)-y),
\end{align*}
where we used that fact that $-i\omega_{j}^{(\lambda)}(iy)=|\omega_{j}^{(\lambda)}(iy)|\ge y$,
hence $\omega_{j}^{(\lambda)}(iy)^{2}=-|\omega_{j}^{(\lambda)}(iy)|^{2}$.
We apply this to $\mu=\widetilde{\mu}_{|X-\lambda|}$, $\mu_{1}=\mu_{1}^{(\lambda)}=\widetilde{\mu}_{|X_{1}-\lambda|}$,
and $\mu_{2}=\widetilde{\mu}_{X_{2}},$ and use the fact that 
\[
\int_{\mathbb{R}}\log(t^{2}+y^{2})\,d\widetilde{\mu}(t)=\int_{\mathbb{R}}\log(t^{2}+y^{2})d\mu(t)
\]
for an arbitrary Borel measure on $\mathbb{R},$ to obtain 
\begin{align*}
\frac{1}{2}\int_{\mathbb{R}}\log(t^{2}+y^{2})\,d\mu_{|X-\lambda|}(t) & =\frac{1}{2}\int_{\mathbb{R}}\log(t^{2}-\omega_{1}^{(\lambda)}(iy)^{2})\,d\mu_{|X_{1}-\lambda|}(t)\\
 & +\frac{1}{2}\int_{\mathbb{R}}\log(t^{2}-\omega_{2}^{(\lambda)}(iy)^{2})\,d\mu_{|X_{2}|}(t)\\
 & -\log(-i\omega_{1}^{(\lambda)}(iy)-i\omega_{2}^{(\lambda)}(iy)-y).
\end{align*}
Equivalently,

\begin{equation}\label{eq:very local}
\begin{aligned}
\frac{1}{2}\tau\bigl(\log(|X-\lambda|^{2}+y^{2})\bigr)
&=\frac{1}{2}\tau\bigl(
  \log(|X_{1}-\lambda|^{2}
  -\omega_{1}^{(\lambda)}(iy)^{2})\bigr)\\
&\quad+\frac{1}{2}\tau\bigl(
  \log(|X_{2}|^{2}
  -\omega_{2}^{(\lambda)}(iy)^{2})\bigr)\\
&\quad-\log\bigl(
  -i\omega_{1}^{(\lambda)}(iy)
  -i\omega_{2}^{(\lambda)}(iy)-y\bigr).
\end{aligned}
\end{equation}
Part (3) of the statement follows now by letting $y\downarrow0$ in
(\ref{eq:very local}). Next, suppose that $\lambda\in F_{1}$. Subtract
$\log(-i\omega_{2}^{(\lambda)}(iy))$ from the last two lines of (\ref{eq:very local})
to see that 
\begin{align*}
\frac{1}{2}\tau(\log(|X-\lambda|^{2}+y^{2}) & =\frac{1}{2}\tau(\log(|X_{1}-\lambda|^{2}-\omega_{1}^{(\lambda)}(iy)^{2})\\
 & +\frac{1}{2}\tau\left(\log\left(\frac{|X_{2}|^{2}}{-\omega_{2}^{(\lambda)}(iy)^{2}}+1\right)\right)\\
 & -\log\left(\frac{-i\omega_{1}^{(\lambda)}(iy)}{-i\omega_{2}^{(\lambda)}(iy)}+1-\frac{y}{-i\omega_{2}^{(\lambda)}(iy)}\right).
\end{align*}
Now, let $y\downarrow0$ and observe that the last two terms above
tend to $0$; note that $m_{2}(X_{2})<+\infty$ because $\lambda\in F_{1}$.
The resulting equation yields (1). To verify (most of) (2), simply
switch the roles of $|X_{1}-\lambda|$ and $|X_{2}|$ in the preceding
argument. 
\end{proof}
\begin{rem}
\label{rem:the strange case of dr. S}For points $\lambda\in S$,
there are two situations to consider. If $\tau(\ker(X_{1}-\lambda))+\tau(\ker X_{2})>1,$
then (see \cite{bv-otaa}) 
\[
\tau(\ker(X-\lambda))=\tau(\ker(X_{1}-\lambda))+\tau(\ker X_{2})-1>0,
\]
and it follows that $\tau(\log(|X-\lambda|))=-\infty$. According
to \cite{bv-otaa}, the only points $\lambda$ such that $\tau(\ker(X-\lambda))>0$
belong to $S$. If $\tau(\ker(X_{1}-\lambda))+\tau(\ker X_{2})=1,$
it may happen that $\tau(\log(|X-\lambda|))$ is finite. In any case,
the values of $\tau(\log(|X-\lambda|))$ for $\lambda$ in a finite
set do not affect the calculation of $\mu_{X}$. 
\end{rem}

We are now ready to prove the first main result of this section. 
\begin{thm}
\label{thm:main-sec3}Let $X_{1},X_{2}\in{\rm Log}^{+}(\tau)\backslash\mathbb{C}$
be two free variables such that $X_{2}$ is $R$-diagonal, set $X=X_{1}+X_{2}$,
and let $\Omega$ be the set introduced in Definition \emph{\ref{def:DSF1F2F}}.
Then the support of the Brown measure $\mu_{X}$ is contained in the
closure of $\Omega$. Moreover, $\mu_{X}|\Omega$ is absolutely continuous
relative to area measure and it has a real-analytic density in that
open set. 
\end{thm}

\begin{proof}
Since the set $S$ and the boundaries of $F_{1}$ and $F_{2}$ are
contained in the closure of $\Omega$, we need to show that $\mu_{X}$
is equal to zero on the interiors of $F_{1}$ and $F_{2}$. We prove
the equivalent statement that the logarithmic integral $\tau(\log(|X-\lambda|)$
is a harmonic function on these interiors. For the set $F_{2}$, this
follows from Proposition \ref{prop:calculating the logarithmic integral}(2).
Indeed, $\tau(\log(|X-\lambda|)$ is constant on $F_{2}$. Now, Proposition
\ref{prop:calculating the logarithmic integral}(1) states that $m_{-2}(|X_{1}-\lambda|)\leq\frac{1}{m_{2}(|X_{2}|)}<+\infty$
for every $\lambda\in F_{1}$.
According to \cite[Theorem 4.5]{pzh-ajm},
this implies that $\mu_{X_{1}}(F_{1})=0,$ and hence $\tau(\log(|X_{1}-\lambda|)$
is a harmonic function for $\lambda$ in the interior of $F_{1}$.
Since, according to the same proposition, $\tau(\log(|X-\lambda|)=\tau(\log(|X_{1}-\lambda|)$
for $\lambda\in F_{1}$, the desired conclusion about the support
of $\mu_{X}$ follows.

To conclude the proof, we show that $\mu_{X}$ has a real-analytic
density on the set $\Omega$. The formula established in Proposition
\ref{prop:calculating the logarithmic integral}(3) shows that it
suffices to prove that $\omega_{j}^{(\lambda)}(0)$ are real-analytic
functions of $\lambda\in\Omega$. 
For $\lambda\in\Omega$, $\omega_{1}^{(\lambda)}(0)$
is the unique fixed point of the map $\psi_{\lambda}=H_{\mu_{2}}\circ H_{\mu_{|X_{1}-\lambda|}}:\mathbb{C}^{+}\to\mathbb{C}^{+}$
and this map is not a conformal automorphism. 
Therefore, Theorem~\ref{thm:dw-subordination} gives
$|\psi_\lambda'(\omega_1^{(\lambda)}(0))|<1$ 
for every $\lambda\in\Omega$.
We observe next that the map $\lambda\mapsto H_{\widetilde{\mu}_{|X_{1}-\lambda|}}(iy)$
is real-analytic for every $y>0$. Indeed, (\ref{eq:explicit HT})
shows that 
\begin{equation}
H_{\widetilde{\mu}_{|X_{1}-\lambda|}}(iy)=\frac{\tau((\lambda-X_{1})^{*}(\lambda-X_{1})((\lambda-X_{1})^{*}(\lambda-X_{1})+y^{2})^{-1})}{-iy\tau(((\lambda-X_{1})^{*}(\lambda-X_{1})+y^{2})^{-1})}.\label{eq:explicit HmutildeX1-lambda}
\end{equation}
It follows that the map $(\lambda,w)\mapsto\psi_{\lambda}(w)$ is
a real-analytic map on $\Omega\times i\mathbb{R}_{+}$. The real analyticity
of $\omega_{1}^{(\lambda)}$ follows now from the implicit function
theorem because $\partial_{y}(iy-\psi_{\lambda}(iy))=i(1-\psi'_{\lambda}(iy))\ne0$
for $iy=\omega_{1}^{(\lambda)}(0)$. The fact that $\omega_{2}^{(\lambda)}(0)$
is real-analytic follows the same way if we replace $\psi_{\lambda}$
by $H_{\widetilde{\mu}_{|X_{1}-\lambda|}}\circ H_{\widetilde{\mu}_{|X_{2}|}}$. 
\end{proof}
Effective calculation of the density of $\mu_{X}$ depends on finding
sufficiently explicit expressions for the subordination functions
$\omega_{j}^{(\lambda)}(0)$ and for their derivatives. Recall that
$\mu_{X}$ is equal to 
\[
\frac{1}{2\pi}\Delta_{\lambda}\tau(\log(|X-\lambda|))=\frac{2}{\pi}\overline{\partial}_{\lambda}\partial_{\lambda}\tau(\log(|X-\lambda|)).
\]
Inside $\Omega$, these differential operators can be applied in the
classical sense to yield the density of $\mu_{X}$ relative to area
measure. We begin with the calculation of $2\partial_{\lambda}\tau(\log(|X-\lambda|))$,
using Proposition \ref{prop:calculating the logarithmic integral}(3). 
\begin{lem}
\label{lem:magic simplification}With the notation of Proposition
\emph{\ref{prop:calculating the logarithmic integral}}, we have 
\[
2\partial_{\lambda}\tau(\log(|X-\lambda|))=\tau((\lambda-X_{1})^{*}(|\lambda-X_{1}|^{2}-\omega_{1}^{(\lambda)}(0)^{2})^{-1}),\quad\lambda\in\Omega.
\]
\end{lem}

\begin{proof}
Differentiating the identity in Proposition \ref{prop:calculating the logarithmic integral}(3),
we obtain 
\begin{align*}
2\partial_{\lambda}\tau(\log(|X-\lambda|) & =\partial_{\lambda}\tau(\log(|X_{1}-\lambda|^{2}-\omega_{1}^{(\lambda)}(0)^{2})+\partial_{\lambda}\tau(\log(|X_{2}|^{2}-\omega_{2}^{(\lambda)}(0)^{2})\\
 & -2\partial_{\lambda}\log(-i\omega_{1}^{(\lambda)}(0)-i\omega_{2}^{(\lambda)}(0)).
\end{align*}
The three terms on the right can be written as
\begin{equation}
    \begin{aligned}\label{eq:(3.6)}
        &\partial_{\lambda}\tau(\log(|X_{1}-\lambda|^{2}-\omega_{1}^{(\lambda)}(0)^{2})\\ &=\partial_{\lambda}\tau(\log((\lambda-X_{1})^{*}(\lambda-X_{1})-\omega_{1}^{(\lambda)}(0)^{2})\\
 & =\tau((\lambda-X_{1})^{*}(|\lambda-X_{1}|^{2}-\omega_{1}^{(\lambda)}(0)^{2})^{-1})\\
 & \qquad\qquad-2\omega_{1}^{(\lambda)}(0)(\partial_{\lambda}\omega_{1}^{(\lambda)}(0))\tau((|\lambda-X_{1}|^{2}-\omega_{1}^{(\lambda)}(0)^{2})^{-1}),
    \end{aligned}
\end{equation}
\begin{equation}
\partial_{\lambda}\tau(\log(|X_{2}|^{2}-\omega_{2}^{(\lambda)}(0)^{2})=-2\omega_{2}^{(\lambda)}(0)(\partial_{\lambda}\omega_{2}^{(\lambda)}(0))\tau((|X_{2}|^{2}-\omega_{2}^{(\lambda)}(0)^{2})^{-1}),\label{eq:two}
\end{equation}
and 
\begin{equation}
-2\partial_{\lambda}\log(-i\omega_{1}^{(\lambda)}(0)-i\omega_{2}^{(\lambda)}(0))=-\frac{2}{\omega_{1}^{(\lambda)}(0)+\omega_{2}^{(\lambda)}(0)}(\partial_{\lambda}\omega_{1}^{(\lambda)}(0)+\partial_{\lambda}\omega_{2}^{(\lambda)}(0)),\label{eq:three}
\end{equation}
where we have used the differentiation rule under the trace pointed
out in \cite{FuKa-det} and \cite{HA-SC-brown}. In order to simplify
the result, we note that 
\begin{align*}
\omega_{1}^{(\lambda)}(0)+\omega_{2}^{(\lambda)}(0) & =\frac{1}{G_{\mu_{j}^{(\lambda)}}(\omega_{j}^{(\lambda)}(0))},
\end{align*}
by Corollary \ref{cor:when omega1(0) is in H}. Thus (\ref{eq:three})
yields 
\begin{align}
   \label{eqn:3.10.G-identity}
    \begin{split}
        -\frac{2}{\omega_{1}^{(\lambda)}(0)+\omega_{2}^{(\lambda)}(0)} & =2\omega_{1}^{(\lambda)}(0)\tau(|\lambda-X_{1}|^{2}-\omega_{1}^{(\lambda)}(0)^{2})^{-1})\\
 & =2\omega_{2}^{(\lambda)}(0)\tau(|X_{2}|^{2}-\omega_{2}^{(\lambda)}(0)^{2})^{-1})
    \end{split}
\end{align}
for $j=1,2$. Combining this with (\ref{eq:(3.6)}) and (\ref{eq:two}).
\end{proof}
We proceesd next to the second derivative. 
\begin{lem}
\label{lemma:lasciate ogne speranza voi ch'entrate}For every $\lambda\in\Omega$,
we have 
\begin{align}
&\frac{1}{2\pi}\Delta_{\lambda}\tau(\log|X-\lambda|) \nonumber\\
& =-\frac{1}{\pi}\omega_{1}^{(\lambda)}(0)^{2}\tau((|\lambda-X_{1}|^{2}-\omega_{1}^{(\lambda)}(0)^{2})^{-1}(|(\lambda-X_{1})^{*}|^{2}-\omega_{1}^{(\lambda)}(0)^{2})^{-1})\label{eq:speranza}\\
 &\qquad +\frac{2}{\pi}\omega_{1}^{(\lambda)}(0)(\overline{\partial}_{\lambda}\omega_{1}^{(\lambda)}(0))\tau((\lambda-X_{1})^{*}(|\lambda-X_{1}|^{2}-\omega_{1}^{(\lambda)}(0)^{2})^{-2}).\nonumber 
\end{align}
\end{lem}

\begin{proof}
We write the result of Lemma \ref{lem:magic simplification} as 
\[
2\partial_{\lambda}\tau(\log|X-\lambda|)=\tau(Z_{\lambda}^{*}(Z_{\lambda}^{*}Z_{\lambda}-\omega_{1}^{(\lambda)}(0)^{2})^{-1}),
\]
where $Z_{\lambda}=\lambda-X_{1}$ satisfies $\overline{\partial}_{\lambda}Z_{\lambda}^{*}=1$
and $\overline{\partial}_{\lambda}Z_{\lambda}=0$. We apply $\overline{\partial}_{\lambda}$
to this identity and obtain 
\begin{align*}
 & \quad2\overline{\partial}_{\lambda}\partial_{\lambda}\tau(\log(|X-\lambda|))\\
 & =\tau((Z_{\lambda}^{*}Z_{\lambda}-\omega_{1}^{(\lambda)}(0)^{2})^{-1})\\
 & -\tau(Z_{\lambda}^{*}(Z_{\lambda}^{*}Z_{\lambda}-\omega_{1}^{(\lambda)}(0)^{2})^{-1}(Z_{\lambda}-2\omega_{1}^{(\lambda)}(0)\overline{\partial}_{\lambda}\omega_{1}^{(\lambda)}(0))(Z_{\lambda}^{*}Z_{\lambda}-\omega_{1}^{(\lambda)}(0)^{2})^{-1})\\
 & =\tau((Z_{\lambda}^{*}Z_{\lambda}-\omega_{1}^{(\lambda)}(0)^{2})^{-1})-\tau(Z_{\lambda}^{*}(Z_{\lambda}^{*}Z_{\lambda}-\omega_{1}^{(\lambda)}(0)^{2})^{-1}Z_{\lambda}(Z_{\lambda}^{*}Z_{\lambda}-\omega_{1}^{(\lambda)}(0)^{2})^{-1})\\
 & +2\omega_{1}^{(\lambda)}(0)(\overline{\partial}_{\lambda}\omega_{1}^{(\lambda)}(0))\tau(Z_{\lambda}^{*}(Z_{\lambda}^{*}Z_{\lambda}-\omega_{1}^{(\lambda)}(0)^{2})^{-2}).
\end{align*}
The equality $Z_{\lambda}(Z_{\lambda}^{*}Z_{\lambda}-\omega_{1}^{(\lambda)}(0)^{2})^{-1}=(Z_{\lambda}Z_{\lambda}^{*}-\omega_{1}^{(\lambda)}(0)^{2})^{-1}Z_{\lambda}$,
and the trace identity, allow us to replace the second term in the
last expression by 
\[
\tau((Z_{\lambda}^{*}Z_{\lambda}-\omega_{1}^{(\lambda)}(0)^{2})^{-1}(Z_{\lambda}Z_{\lambda}^{*}-\omega_{1}^{(\lambda)}(0)^{2})^{-1}Z_{\lambda}Z_{\lambda}^{*}),
\]
and thus 
\begin{align*}
 & \quad\overline{\partial}_{\lambda}\partial_{\lambda}\tau(\log(|X-\lambda|))\\
 & =\tau((Z_{\lambda}^{*}Z_{\lambda}-\omega_{1}^{(\lambda)}(0)^{2})^{-1})\\
 & -\tau((Z_{\lambda}^{*}Z_{\lambda}-\omega_{1}^{(\lambda)}(0)^{2})^{-1}(Z_{\lambda}Z_{\lambda}^{*}-\omega_{1}^{(\lambda)}(0)^{2})^{-1}Z_{\lambda}Z_{\lambda}^{*})\\
 & +2\omega_{1}^{(\lambda)}(0)(\overline{\partial}_{\lambda}\omega_{1}^{(\lambda)}(0))\tau(Z_{\lambda}^{*}(Z_{\lambda}^{*}Z_{\lambda}-\omega_{1}^{(\lambda)}(0)^{2})^{-2})\\
 & =\tau((Z_{\lambda}^{*}Z_{\lambda}-\omega_{1}^{(\lambda)}(0)^{2})^{-1})\\
 & -\tau((Z_{\lambda}^{*}Z_{\lambda}-\omega_{1}^{(\lambda)}(0)^{2})^{-1}Z_{\lambda}Z_{\lambda}^{*}(Z_{\lambda}Z_{\lambda}^{*}-\omega_{1}^{(\lambda)}(0)^{2})^{-1})\\
 & +2\omega_{1}^{(\lambda)}(0)(\overline{\partial}_{\lambda}\omega_{1}^{(\lambda)}(0))\tau(Z_{\lambda}^{*}(Z_{\lambda}^{*}Z_{\lambda}-\omega_{1}^{(\lambda)}(0)^{2})^{-2})\\
 & =\tau((Z_{\lambda}^{*}Z_{\lambda}-\omega_{1}^{(\lambda)}(0)^{2})^{-1}[(Z_{\lambda}Z_{\lambda}^{*}-\omega_{1}^{(\lambda)}(0)^{2})-Z_{\lambda}Z_{\lambda}^{*}](Z_{\lambda}Z_{\lambda}^{*}-\omega_{1}^{(\lambda)}(0)^{2})^{-1})\\
 & +2\omega_{1}^{(\lambda)}(0)(\overline{\partial}_{\lambda}\omega_{1}^{(\lambda)}(0))\tau(Z_{\lambda}^{*}(Z_{\lambda}^{*}Z_{\lambda}-\omega_{1}^{(\lambda)}(0)^{2})^{-2}).
\end{align*}
Canceling $Z_{\lambda}Z_{\lambda}^{*}$ yields the stated formula. 
\end{proof}
The formula in Lemma \ref{lemma:lasciate ogne speranza voi ch'entrate}
involves $\overline{\partial}_{\lambda}\omega_{1}^{(\lambda)}(0)$,
which we proceed to calculate next. As in the proof of Theorem \ref{thm:main-sec3},
$\psi'_{\lambda}$ denotes the derivative of $\psi$ for fixed $\lambda$,
while $\overline{\partial}_{\lambda}\psi_{\lambda}(iy)$ denotes the
$\overline{\partial}$ derivative for fixed $y$. 

\begin{lem}
\label{lem:delbaromega} With the notation $\psi_{\lambda}=H_{\widetilde{\mu}_{|X_{2}|}}\circ H_{\widetilde{\mu}_{|X_{1}-\lambda|}}$,
we have 
\[
\overline{\partial}_{\lambda}\omega_{1}^{(\lambda)}(0)=\frac{(\overline{\partial}_{\lambda}\psi_{\lambda})(\omega_{1}^{(\lambda)}(0))}{1-\psi'_{\lambda}(\omega_{1}^{(\lambda)}(0))},\quad\lambda\in\Omega.
\]
\end{lem}

\begin{proof}
Apply $\overline{\partial}_{\lambda}$ to the equation $\omega_{1}^{(\lambda)}(0)=\psi_{\lambda}(\omega_{1}^{(\lambda)}(0))$. 
\end{proof}
The derivative $(\overline{\partial}_{\lambda}\psi_{\lambda})(\omega_{1}^{(\lambda)}(0))$
can be made more explicit.
\begin{lem}
\label{lem:delbar psi_lambda} For $\lambda\in\Omega$, we have 
\begin{align*}
(\overline{\partial}_{\lambda}\psi_{\lambda})(\omega_{1}^{(\lambda)}(0))  =  -H'_{\mu_{2}}(\omega_{2}^{(\lambda)}(0))\frac{\tau((\lambda-X_{1})(\omega_{1}^{(\lambda)}(0)^{2}-|\lambda-X_{1}|^{2})^{-2})}{\omega_{1}^{(\lambda)}(0)\tau((\omega_{1}^{(\lambda)}(0)^{2}-|\lambda-X_{1}|^{2})^{-1})^{2}},
\end{align*}
and 
\[
H'_{\mu_{2}}(\omega_{2}^{(\lambda)}(0))=\frac{\tau((\omega_{2}^{(\lambda)}(0)^{2}+|X_{2}|^{2})(\omega_{2}^{(\lambda)}(0)^{2}-|X_{2}|^{2})^{-2})}{\tau(\omega_{2}^{(\lambda)}(0)(\omega_{2}^{(\lambda)}(0)^{2}-|X_{2}|^{2})^{-1})^{2}}-1.
\]
\end{lem}

\begin{proof}
Since the function $z\mapsto H_{\mu_{2}}(z)$ is analytic in $\mathbb{C}^{+}$,
we have 
\[
(\overline{\partial}_{\lambda}\psi_{\lambda})(z)=H'_{\mu_{2}}(H_{\widetilde{\mu}_{|X_{1}-\lambda|}}(z))(\overline{\partial}_{\lambda}H_{\widetilde{\mu}_{|X_{1}-\lambda|}})(z).
\]
We proceed to calculate the two derivatives in the right hand side,
recalling that $H_{\mu}(z)=(1/G_{\mu}(z))-z$ and that $\mu_{2}$
is a symmetric measure: 
\begin{align*}
H'_{\mu_{2}}(z) & =\frac{-G'_{\mu_{2}}(z)}{(G_{\mu_{2}}(z))^{2}}-1=\frac{{\displaystyle \int_{\mathbb{R}}\frac{1}{(z-t)^{2}}d\mu_{2}}}{{\displaystyle \left(\int_{\mathbb{R}}\frac{1}{z-t}d\mu_{2}\right)^{2}}}-1\\
 & =\frac{{\displaystyle \frac{1}{2}\int_{\mathbb{R}}\left[\frac{1}{(z-t)^{2}}+\frac{1}{(z+t)^{2}}\right]d\mu_{2}}}{{\displaystyle \left(\frac{1}{2}\int_{\mathbb{R}}\left[\frac{1}{z-t}+\frac{1}{z+t}\right]d\widetilde{\mu}_{|X_{2}|}\right)^{2}}}-1=\frac{{\displaystyle \int_{\mathbb{R}}\frac{z^{2}+t^{2}}{(z^{2}-t^{2})^{2}}d\mu_{2}}}{{\displaystyle \left(\int_{\mathbb{R}}\frac{z}{z^{2}-t^{2}}d\mu_{2}\right)^{2}}}-1\\
 & =\frac{\tau((z^{2}+|X_{2}|^{2})(z^{2}-|X_{2}|^{2})^{-2})}{\tau(z(z^{2}-|X_{2}|^{2})^{-1})^{2}}-1.
\end{align*}
We substitute now $H_{\mu_{1}^{(\lambda)}}(\omega_{1}^{(\lambda)}(0))=\omega_{2}^{(\lambda)}(0)$
to obtain 
\begin{equation}
H'_{\mu_{2}}(H_{\mu_{1}^{(\lambda)}}(\omega_{1}^{(\lambda)}(0)))=\frac{\tau((\omega_{2}^{(\lambda)}(0)^{2}+|X_{2}|^{2})(\omega_{2}^{(\lambda)}(0)^{2}-|X_{2}|^{2})^{-2})}{\tau(\omega_{2}^{(\lambda)}(0)(\omega_{2}^{(\lambda)}(0)^{2}-|X_{2}|^{2})^{-1})^{2}}-1.\label{eq:3.11}
\end{equation}
Next, 
\begin{align*}
(\overline{\partial}_{\lambda}H_{\mu_{1}^{(\lambda)}})(z) & =\frac{-(\overline{\partial}_{\lambda}G_{\mu_{1}^{(\lambda)}})(z)}{((G_{\mu_{1}^{(\lambda)}})(z))^{2}}=\frac{-\overline{\partial}_{\lambda}\tau(z(z^{2}-|X_{1}-\lambda|^{2})^{-1})}{\tau(z(z^{2}-|X_{1}-\lambda|^{2})^{-1})^{2}}\\
 & =-\frac{\tau((\lambda-X_{1})(z^{2}-|X_{1}-\lambda|^{2})^{-2})}{z\cdot\tau((z^{2}-|X_{1}-\lambda|^{2})^{-1})^{2}},
\end{align*}
and substituting $\omega_{1}^{(\lambda)}(0)$ for of $z$, 
\begin{equation}
(\overline{\partial}_{\lambda}H_{\mu_{1}^{(\lambda)}})(\omega_{1}^{(\lambda)}(0))=-\frac{\tau((\lambda-X_{1})(\omega_{1}^{(\lambda)}(0)^{2}-|\lambda-X_{1}|^{2})^{-2})}{\omega_{1}^{(\lambda)}(0)\tau((\omega_{1}^{(\lambda)}(0)^{2}-|\lambda-X_{1}|^{2})^{-1})^{2}}.\label{eq:3.12}
\end{equation}
The stated formula is obtained by multiplying (\ref{eq:3.11}) and
(\ref{eq:3.12}). 
\end{proof}
Since $\tau(\log|X-\lambda|)$ is subharmonic, the Laplacian $\overline{\partial}_{\lambda}\partial_{\lambda}\tau(\log|X-\lambda|)$
must be a nonnegative function in $\Omega$. It may be worthwhile
to observe that the formula for this function that follows from the
preceding lemmas does obviously yield a real-valued function in $\Omega$.

\begin{lem}
For $\lambda\in\Omega$, we have 
\begin{equation}
 \label{eqn:H_1-derivative}
    H'_{\widetilde{\mu}_{|X_{1}-\lambda|}}(\omega_{1}^{(\lambda)}(0))=\frac{\tau((\omega_{1}^{(\lambda)}(0)^{2}+|\lambda-X_{1}|^{2})(\omega_{1}^{(\lambda)}(0)^{2}-|\lambda-X_{1}|^{2})^{-2})}{\tau(\omega_{1}^{(\lambda)}(0)(\omega_{1}^{(\lambda)}(0)^{2}-|\lambda-X_{1}|^{2})^{-1})^{2}}-1.
\end{equation}
Consequently, we have 
\begin{align} 
     &\frac{2\tau(Z_\lambda^*Z_\lambda(\omega_{1}^{(\lambda)}(0))^2-Z_\lambda^*Z_\lambda)^{-2})}{\tau((\omega_{1}^{(\lambda)}(0))^2-Z_\lambda^*Z_\lambda)^{-1})^2}\nonumber\\
&\quad\quad=\omega_{1}^{(\lambda)}(0)^2(H'_{\widetilde{\mu}_{|X_{1}-\lambda|}}(\omega_{1}^{(\lambda)}(0))+1)
 -\frac{1}{\tau((\omega_{1}^{(\lambda)}(0))^2-Z_\lambda^*Z_\lambda)^{-1})}\label{eqn:3.15.2nd.part.identity}\\
&\quad\quad=\omega_{1}^{(\lambda)}(0)^2\cdot H'_{\widetilde{\mu}_{|X_{1}-\lambda|}}(\omega_{1}^{(\lambda)}(0))-\omega_{1}^{(\lambda)}(0)\cdot \omega_{2}^{(\lambda)}(0),\nonumber
\end{align}
where $Z_\lambda=\lambda-X_1$.
\end{lem}
\begin{proof}
A similar calculation as $H'_{\mu_2}(\omega_{2}^{(\lambda)}(0))$ as in Lemma \ref{lem:delbaromega} yields \eqref{eqn:H_1-derivative}. By multiplying $\omega_{1}^{(\lambda)}(0)$ for both sides of \eqref{eqn:H_1-derivative} and arrange terms we then have 
\begin{align*}
&\omega_{1}^{(\lambda)}(0)^2(H'_{\widetilde{\mu}_{|X_{1}-\lambda|}}(\omega_{1}^{(\lambda)}(0))+1)\\
    =&\frac{\tau((\omega_{1}^{(\lambda)}(0)^{2}+|\lambda-X_{1}|^{2})(\omega_{1}^{(\lambda)}(0)^{2}-|\lambda-X_{1}|^{2})^{-2})}{\tau((\omega_{1}^{(\lambda)}(0)^{2}-|\lambda-X_{1}|^{2})^{-1})^{2}}\\
    =&\frac{\tau((\omega_{1}^{(\lambda)}(0)^{2}-|\lambda-X_{1}|^{2})(\omega_{1}^{(\lambda)}(0)^{2}-|\lambda-X_{1}|^{2})^{-2})}{\tau((\omega_{1}^{(\lambda)}(0)^{2}-|\lambda-X_{1}|^{2})^{-1})^{2}}\\
      &\qquad\qquad\quad+\frac{2\tau(|\lambda-X_{1}|^{2}(\omega_{1}^{(\lambda)}(0)^{2}-|\lambda-X_{1}|^{2})^{-2})}{\tau((\omega_{1}^{(\lambda)}(0)^{2}-|\lambda-X_{1}|^{2})^{-1})^{2}}\\
   =&\frac{1}{\tau((\omega_{1}^{(\lambda)}(0))^2-Z_\lambda^*Z_\lambda)^{-1})}+
    \frac{2\tau(Z_\lambda^*Z_\lambda(\omega_{1}^{(\lambda)}(0))^2-Z_\lambda^*Z_\lambda)^{-2})}{\tau((\omega_{1}^{(\lambda)}(0))^2-Z_\lambda^*Z_\lambda)^{-1})^2}.
\end{align*}
By \eqref{eqn:3.10.G-identity} we have
\begin{align}
  \label{eqn:3.17.sum.identity}
    \frac{1}{\tau((\omega_{1}^{(\lambda)}(0))^2-Z_\lambda^*Z_\lambda)^{-1})}
      =\omega_{1}^{(\lambda)}(0)^2+\omega_{1}^{(\lambda)}(0)\cdot\omega_{2}^{(\lambda)}(0).
\end{align}
This completes the proof for \eqref{eqn:3.15.2nd.part.identity}.
\end{proof}

We are now ready to prove our second main result in this section. 
\begin{thm}
\label{thm:positivity}
For $\lambda\in\Omega$, the absolutely continuous part of $\mu_{X_1+X_2}$ has a density function $\rho$ given by
\begin{align} \label{eqn:density-rho}
    \begin{split}
         \rho(\lambda)&= \frac{1}{2\pi}\Delta_{\lambda}\tau(\log|X-\lambda|)\\
 &=-\frac{1}{\pi}\omega_{1}^{(\lambda)}(0)^2\tau((\omega_{1}^{(\lambda)}(0)^2-Z_\lambda^*Z_\lambda)^{-1}(\omega_{1}^{(\lambda)}(0)^2-Z_\lambda Z_\lambda^*)^{-1}) \\
  &\quad\quad\quad\quad\quad-\frac{2}{\pi}\frac{H_{\mu_2}^\prime(\omega_{2}^{(\lambda)}(0))}{1-\psi_\lambda^\prime(\omega_{1}^{(\lambda)}(0))}\frac{|\tau(Z_\lambda(\omega_{1}^{(\lambda)}(0)^2-Z_\lambda^*Z_\lambda)^{-2})|^2}{\tau((\omega_{1}^{(\lambda)}(0)^2-Z_\lambda^*Z_\lambda)^{-1})^2},
    \end{split}
\end{align}
where $Z_\lambda=\lambda-X_1$. 
Moreover, $\rho(\lambda)>0$ for any $\lambda\in\Omega$. 
\end{thm}
\begin{proof}
The density formula is a direct consequence of Lemma \ref{lemma:lasciate ogne speranza voi ch'entrate},
Lemma \ref{lem:delbaromega} and Lemma \ref{lem:delbar psi_lambda}. It remains to prove that $\rho(\lambda)>0$. 

We observe that 
\[
Z_\lambda(Z_\lambda^*Z_\lambda-\omega_{1}^{(\lambda)}(0)^2)^{-1}
=(Z_\lambda Z_\lambda^*-\omega_{1}^{(\lambda)}(0)^2)^{-1}Z_\lambda.
\]
Denote $h=Z_\lambda^*Z_\lambda-\omega_{1}^{(\lambda)}(0)^2$ and $k=Z_{\lambda}Z_\lambda^*-\omega_{1}^{(\lambda)}(0)^2$. We note that $h>0$ and $k>0$. Following an estimate from the proof of \cite[Lemma 7.11]{byzh-ADV}, the Cauchy-Schwarz inequality yields
\begin{equation}
\label{eq:CS-estimate}
\begin{aligned}
\left|\tau\left(
Z_\lambda
\bigl(\omega_1^{(\lambda)}(0)^2-Z_\lambda^*Z_\lambda\bigr)^{-2}
\right)\right|^2&=|\tau(k^{-1}Z_\lambda h^{-1})|^2\\
&=\left|\tau\left(
(h^{-1/2}k^{-1/2})
(k^{-1/2}Z_\lambda h^{-1/2})
\right)\right|^2\\
&\le
\tau(h^{-1}k^{-1})\,
\tau\left(
k^{-1/2}Z_\lambda h^{-1}Z_\lambda^*k^{-1/2}
\right)\\
&=
\tau(h^{-1}k^{-1})\,
\tau(Z_\lambda^*k^{-1}Z_\lambda h^{-1})\\
&=
\tau(h^{-1}k^{-1})\,
\tau(Z_\lambda^*Z_\lambda h^{-2}).
\end{aligned}
\end{equation}
Hence, we have 
\begin{align}
&|\tau(Z_\lambda(\omega_{1}^{(\lambda)}(0)^2-Z_\lambda^*Z_\lambda)^{-2})|^2=|\tau(Z_\lambda(Z_\lambda^*Z_\lambda-\omega_{1}^{(\lambda)}(0)^2)^{-2})|^2\nonumber\\
&\leq\tau((Z_\lambda Z_\lambda^*-\omega_{1}^{(\lambda)}(0)^2)^{-1}
(Z_\lambda^*Z_\lambda-\omega_{1}^{(\lambda)}(0)^2)^{-1})\tau(Z_\lambda^*Z_\lambda(Z_\lambda^*Z_\lambda-\omega_{1}^{(\lambda)}(0)^2)^{-2}).\nonumber
\end{align}
We note that $-\omega_{1}^{(\lambda)}(0)^2>0$ and hence the first term of \eqref{eqn:density-rho} is strictly positive. Recall that $1-\psi_\lambda^\prime(\omega_{1}^{(\lambda)}(0))>0$. It suffices to prove $\rho(\lambda)>0$ when $H_{\mu_2}^\prime(\omega_{2}^{(\lambda)}(0))>0$. When $H_{\mu_2}^\prime(\omega_{2}^{(\lambda)}(0))>0$, we have 
\begin{align*}
    \rho(\lambda)\geq& \frac{1}{\pi}\tau((\omega_{1}^{(\lambda)}(0)^2-Z_\lambda^*Z_\lambda)^{-1}(\omega_{1}^{(\lambda)}(0)^2-Z_\lambda Z_\lambda^*)^{-1})\\
     &\times \left(-\omega_{1}^{(\lambda)}(0)^2 -\frac{2 H_{\mu_2}^\prime(\omega_{2}^{(\lambda)}(0))}{1-\psi_\lambda^\prime(\omega_{1}^{(\lambda)}(0))}\frac{\tau(Z_\lambda^* Z_\lambda(\omega_{1}^{(\lambda)}(0)^2-Z_\lambda^*Z_\lambda)^{-2})}{\tau((\omega_{1}^{(\lambda)}(0)^2-Z_\lambda^*Z_\lambda)^{-1})^2} \right).
\end{align*}
By applying the identity \eqref{eqn:3.15.2nd.part.identity} we obtain
\begin{align*}
        &-\omega_{1}^{(\lambda)}(0)^2 -\frac{2 H_{\mu_2}^\prime(\omega_{2}^{(\lambda)}(0))}{1-\psi_\lambda^\prime(\omega_{1}^{(\lambda)}(0))}\frac{\tau(Z_\lambda^* Z_\lambda(\omega_{1}^{(\lambda)}(0)^2-Z_\lambda^*Z_\lambda)^{-2})}{\tau((\omega_{1}^{(\lambda)}(0)^2-Z_\lambda^*Z_\lambda)^{-1})^2} \nonumber\\
    =&-\omega_{1}^{(\lambda)}(0)^2 -\frac{H_{\mu_2}^\prime(\omega_{2}^{(\lambda)}(0)) (\omega_{1}^{(\lambda)}(0)^2\cdot H'_{\widetilde{\mu}_{|X_{1}-\lambda|}}(\omega_{1}^{(\lambda)}(0))-\omega_{1}^{(\lambda)}(0)\cdot \omega_{2}^{(\lambda)}(0)) }{1-H'_{\widetilde{\mu}_{|X_{1}-\lambda|}}(\omega_{1}^{(\lambda)}(0))
 H_{\mu_2}^\prime(\omega_{2}^{(\lambda)}(0))}\nonumber\\
  =&\frac{\omega_{1}^{(\lambda)}(0)\cdot (-\omega_{1}^{(\lambda)}(0)+\omega_{2}^{(\lambda)}(0)H_{\mu_2}^\prime(\omega_{2}^{(\lambda)}(0)))}{1-H'_{\widetilde{\mu}_{|X_{1}-\lambda|}}(\omega_{1}^{(\lambda)}(0))
 H_{\mu_2}^\prime(\omega_{2}^{(\lambda)}(0))}
\end{align*}
where we used the identity
\[
\psi_\lambda^\prime(\omega_{1}^{(\lambda)}(0))=H'_{\widetilde{\mu}_{|X_{1}-\lambda|}}(\omega_{1}^{(\lambda)}(0))
 H_{\mu_2}^\prime(\omega_{2}^{(\lambda)}(0).
\]
Since $H_{\mu_2}$ is not a conformal automorphism, the Schwarz--Pick inequality implies that 
\begin{align}
|H_{\mu_2}^\prime(\omega_{2,\lambda})|
<
\frac{\Im H_{\mu_2}(\omega_2^{(\lambda)}(0))}{\Im\omega_2^{(\lambda)}(0)}
=
\frac{\Im\omega_1^{(\lambda)}(0)}{\Im\omega_{2,\lambda}}.
\label{eq:schwarz-pick-H2}
\end{align}
We note that $\omega_j^{(\lambda)}(0)$ is pure imaginary and $\Im \omega_2^{(\lambda)}(0)>0$. Then \eqref{eq:schwarz-pick-H2} implies that
\[
 \omega_{1}^{(\lambda)}(0)\cdot (-\omega_{1}^{(\lambda)}(0)+\omega_{2}^{(\lambda)}(0)H_{\mu_2}^\prime(\omega_{2}^{(\lambda)}(0)))>0.
\]
We therefore prove that $\rho(\lambda)>0$ for $\lambda\in\Omega$. 
\end{proof}

\section{Regularity properties of $\mu_{X_1+X_2}$}
In this section, we prove regularity results for the Brown measure $\mu_X=\mu_{X_1+X_2}$. 

\subsection{No singular continuous part}
The following theorem is the first main result of this section. We first show that $\mu_X$ has no singular continuous part. We denote planar Lebesgue measure by $\di{A}$ and use the density $\rho$ calculated in Theorem~\ref{thm:positivity}. 

\begin{thm}
\label{thm:no-boundary-singular}
We have a nonnegative $\rho_X\in L^1(m_2)$ such that
\begin{equation}\label{eq:main}
 \mu_X=\rho_X\di{A}+
       \sum_{\lambda\in S}\mu_X(\{\lambda\})\delta_\lambda.
\end{equation}
In particular, $\mu_X$ has no singular continuous part.
\end{thm}

By Theorem~\ref{thm:main-sec3}, $\mu_X$ is supported on $\overline\Omega$ and is absolutely continuous on $\Omega$. Since $S$ is finite and $F\subset F_2$, it suffices to prove
absolute continuity in a neighborhood of each point of
$\partial\Omega\cap(F_1\setminus F)$ and to show that $\mu_X(F_2)=0$.

In the circular case, \cite{byzh-ADV} proved absolute
continuity by calculating an upper bound for the density of the regularized Brown measure and
passing to the weak limit. We use this argument locally near the boundary points in $F_1$. For $F_2$, we use Green's theorem to estimate the regularized mass of disks enclosing $F_2$ in terms of the first derivative of the regularized logarithmic potential. For $\varepsilon>0$, set
\[
 \Phi_\varepsilon(\lambda)=\frac{1}{2}\tau(\log(|X-\lambda|^2+\varepsilon^2))
\]
and 
\[
\mu_X^{(\varepsilon)}=\frac{1}{2\pi}\Delta_\lambda\Phi_\varepsilon(\lambda)=\rho_\varepsilon(\lambda)\di{A}
\]
Recall that $\mu_X^{(\varepsilon)}\Rightarrow\mu_X$ as $\varepsilon\downarrow0$ (see \cite[Lemma 2.8]{HA-SC-brown}). 

In the following proof, we use the notation $\omega_j^{(\lambda)}(i\varepsilon)=\iu\eta_j(\lambda,\varepsilon),\ \lambda\in\CC$ for $j=1,2$ when $\varepsilon$ is fixed. We may also denote $\eta_j=\eta_j(\lambda,\varepsilon)$ when there is no ambiguity, and as in the preceding section, set
$Z_\lambda=\lambda-X_1$. The subordination relations give
\begin{align}
\eta_1(\tau({Z_\lambda^*Z_\lambda+\eta_1^2})^{-1})
=\eta_2\tau((|X_2|^2+\eta_2^2)^{-1})
=\frac{1}{\eta_1+\eta_2-\varepsilon}.
\label{eq:regularized-subordination}
\end{align}

We first record the density formula for the regularized Brown measure, which is an analogue of Theorem~\ref{thm:positivity}.

\begin{lem}
For every $\varepsilon>0$,
\begin{align}
\partial_\lambda\Phi_\varepsilon(\lambda)
&=\frac12\tau(h^{-1}Z_\lambda^*),\label{eq:regularized-gradient}\\
\pi\rho_\varepsilon(\lambda)
&=\eta_1^2\tau(h^{-1}k^{-1})
-\frac{2H_{\mu_2}'(\iu\eta_2)|\tau(Z_\lambda h^{-2})|^2}
{\tau(h^{-1})^2\bigl(1-H_{\mu_1^{(\lambda)}}'(\iu\eta_1)H_{\mu_2}'(\iu\eta_2)\bigr)},
\label{eq:regularized-density}
\end{align}
where $h=Z_\lambda^*Z_\lambda+\eta_1^2$ and
$k=Z_\lambda Z_\lambda^*+\eta_1^2$.
Moreover,
$|H_{\mu_1^{(\lambda)}}'(\iu\eta_1)H_{\mu_2}'(\iu\eta_2)|<1$, and $\eta_1$ and $\eta_2$ are real-analytic functions in terms of $\lambda$.
\end{lem}

\begin{proof}
By symmetry, $H_{\mu_1^{(\lambda)}}'(\iu\eta_1)$ and
$H_{\mu_2}'(\iu\eta_2)$ are real. Theorem~\ref{thm:dw-subordination},
applied at $z=\iu\varepsilon$, gives $|H_{\mu_1^{(\lambda)}}'(\iu\eta_1)H_{\mu_2}'(\iu\eta_2)|<1$.
Thus, for fixed $\varepsilon>0$, the implicit function theorem applied to the subordination equations shows that $\eta_1$ and $\eta_2$ depend smoothly on $\lambda$.

Setting $y=\varepsilon$ in \eqref{eq:very local} gives
\begin{align*}
\Phi_\varepsilon(\lambda)
=\frac12\tau(\log (Z_\lambda^*Z_\lambda+\eta_1^2))+\frac12\tau(\log(|X_2|^2+\eta_2^2))
-\log(\eta_1+\eta_2-\varepsilon).
\end{align*}
By \eqref{eq:regularized-subordination} and applying a similar calculation as Lemma \ref{lem:magic simplification}, the derivatives of the right-hand side with respect to $\eta_1$ and $\eta_2$ vanish, which implies 
\[
\partial_\lambda\Phi_\varepsilon(\lambda)
=\frac12\tau(h^{-1}Z_\lambda^*). 
\]

At fixed $\eta_1$, we have
$\overline\partial_\lambda\tau((Z_\lambda^*Z_\lambda+\eta_1^2)^{-1})=-\tau(Z_\lambda (Z_\lambda^*Z_\lambda+\eta_1^2)^{-2})$. Differentiating the subordination equations yields
\begin{align}
\overline\partial_\lambda\eta_1(\lambda)
=\frac{H_{\mu_2}'(\iu\eta_2)\tau(Z_\lambda h^{-2})}
{\eta_1\tau(h^{-1})^2\bigl(1-H_{\mu_1^{(\lambda)}}'(\iu\eta_1)H_{\mu_2}'(\iu\eta_2)\bigr)}.
\label{eq:delbar-eta-regularized}
\end{align}
Differentiating \eqref{eq:regularized-gradient}, as in the proof of Lemma~\ref{lemma:lasciate ogne speranza voi ch'entrate}, gives
\begin{align*}
2\overline\partial_\lambda\partial_\lambda\Phi_\varepsilon(\lambda)
=\eta_1^2\tau(h^{-1}k^{-1})
-2\eta_1(\overline\partial_\lambda\eta_1)\tau(Z_\lambda^*h^{-2}).
\end{align*}
Substitution of \eqref{eq:delbar-eta-regularized} proves \eqref{eq:regularized-density}.
\end{proof}

When $H_{\mu_2}'(\iu\eta_2)\ge0$, the first term in
\eqref{eq:regularized-density} bounds $\pi\rho_\varepsilon(\lambda)$
from above. When $H_{\mu_2}'(\iu\eta_2)<0$, we divide the numerator
and denominator of the second term by $-2H_{\mu_2}'(\iu\eta_2)$.
The following lemma gives a lower bound for the resulting denominator
in terms of two resolvent variances. The variance associated with
$|X_2|$ will be used when $|X_1-\lambda_0|$ is scalar at a point
$\lambda_0\in\partial\Omega\cap(F_1\setminus F)$.

\begin{lem}
\label{lem:denominator}
If $H_{\mu_2}'(\iu\eta_2)<0$, then
\begin{align}
&\frac{\tau((|Z_\lambda|^2+\eta_1^2)^{-1})^2(H_{\mu_1^{(\lambda)}}'(\iu\eta_1)H_{\mu_2}'(\iu\eta_2)-1)}{2H_{\mu_2}'(\iu\eta_2)}
\nonumber\\
&\quad\ge\tau((|Z_\lambda|^2+\eta_1^2)^{-2})-\tau((|Z_\lambda|^2+\eta_1^2)^{-1})^2\nonumber\\
&\qquad\qquad
+\tau((|Z_\lambda|^2+\eta_1^2)^{-1})^2\frac{\eta_2(\eta_2-\varepsilon)}
{\eta_1(\eta_1-\varepsilon)}
\left(\frac{\tau((|X_2|^2+\eta_2^2)^{-2})}
{\tau((|X_2|^2+\eta_2^2)^{-1})^2}-1\right).
\label{eq:regularized-denominator-bound}
\end{align}
In both cases $H_{\mu_2}'(\iu\eta_2)<0$ and
 $H_{\mu_2}'(\iu\eta_2)\ge0$, whenever $\tau((|Z_\lambda|^2+\eta_1^2)^{-2})>\tau((|Z_\lambda|^2+\eta_1^2)^{-1})^2$,
\begin{equation}
\pi\rho_\varepsilon(\lambda)\le\tau((|Z_\lambda|^2+\eta_1^2)^{-1})
\left(1+\frac{\tau((|Z_\lambda|^2+\eta_1^2)^{-1})^2}
{\tau((|Z_\lambda|^2+\eta_1^2)^{-2})-\tau((|Z_\lambda|^2+\eta_1^2)^{-1})^2}\right).
\label{eq:regularized-density-bound}
\end{equation}
\end{lem}
\begin{proof}
For this computation, put
\[
a=\frac{\tau((|Z_\lambda|^2+\eta_1^2)^{-2})}{\tau((|Z_\lambda|^2+\eta_1^2)^{-1})^2}-1,\qquad
b=\frac{\tau((|X_2|^2+\eta_2^2)^{-2})}
{\tau((|X_2|^2+\eta_2^2)^{-1})^2}-1.
\]
Both are nonnegative by the Cauchy--Schwarz inequality. Differentiating the Cauchy
transforms similar to the proof for Lemma \ref{lem:delbaromega} and using \eqref{eq:regularized-subordination} gives
\[
H_{\mu_1^{(\lambda)}}'(\iu\eta_1)
=2a-\frac{\eta_2-\varepsilon}{\eta_1},\qquad
H_{\mu_2}'(\iu\eta_2)
=2b-\frac{\eta_1-\varepsilon}{\eta_2}.
\]
If $H_{\mu_2}'(\iu\eta_2)<0$, then
\begin{align*}
&\frac12\left(H_{\mu_1^{(\lambda)}}'(\iu\eta_1)
-\frac1{H_{\mu_2}'(\iu\eta_2)}\right)\\
&\quad=a-\frac{\eta_2-\varepsilon}{2\eta_1}
+\frac1{2((\eta_1-\varepsilon)/\eta_2-2b)}\\
&\quad=a+\frac{\displaystyle\frac{\eta_2-\varepsilon}{\eta_1}b
+\displaystyle\frac{\varepsilon(\eta_1+\eta_2-\varepsilon)}{2\eta_1\eta_2}}
{\displaystyle\frac{\eta_1-\varepsilon}{\eta_2}-2b}.
\end{align*}
The denominator in the last fraction satisfies
\[
0<\frac{\eta_1-\varepsilon}{\eta_2}-2b
\le\frac{\eta_1-\varepsilon}{\eta_2}.
\]
Dropping the nonnegative numerator term containing $\varepsilon$ yields
\[
\frac12\left(H_{\mu_1^{(\lambda)}}'(\iu\eta_1)
-\frac1{H_{\mu_2}'(\iu\eta_2)}\right)
\ge a+\frac{\eta_2(\eta_2-\varepsilon)}
{\eta_1(\eta_1-\varepsilon)}b.
\]
Multiplication by $\tau((|Z_\lambda|^2+\eta_1^2)^{-1})^2$ proves
\eqref{eq:regularized-denominator-bound}.

To prove the inequality \eqref{eq:regularized-density-bound}, we apply the calculation in \eqref{eq:CS-estimate} with
\[
h=|Z_\lambda|^2+\eta_1^2,
\qquad
k=|Z_\lambda^*|^2+\eta_1^2,
\]
we obtain
\[
|\tau(Z_\lambda h^{-2})|^2
\le \tau(h^{-1}k^{-1})\,
    \tau(Z_\lambda^*Z_\lambda h^{-2}).
\]
Since $Z_\lambda^*Z_\lambda=h-\eta_1^2$, we have
\[
\tau(Z_\lambda^*Z_\lambda h^{-2})
=\tau(h^{-1})-\eta_1^2\tau(h^{-2}).
\]
Consequently,
\begin{equation}
\label{eq:regularized-CS}
\begin{aligned}
&|\tau(Z_\lambda(|Z_\lambda|^2+\eta_1^2)^{-2})|^2\\
&\quad\le
\tau\!\left(
(|Z_\lambda|^2+\eta_1^2)^{-1}
(|Z_\lambda^*|^2+\eta_1^2)^{-1}
\right)\\
&\qquad\quad\times
\left(
\tau((|Z_\lambda|^2+\eta_1^2)^{-1})
-\eta_1^2\tau((|Z_\lambda|^2+\eta_1^2)^{-2})
\right).
\end{aligned}
\end{equation}
The cyclic
permutations above involve bounded factors and remain
valid when $Z_\lambda$ is unbounded and affiliated.

In the negative-derivative case, the regularized density \eqref{eq:regularized-density} becomes
 \begin{align}    \pi\rho_\varepsilon(\lambda)&=\eta_1^2\tau((|Z_\lambda|^2+\eta_1^2)^{-1}(|Z_\lambda^*|^2+\eta_1^2)^{-1})\nonumber\\
  &\qquad+\frac{|\tau(Z_\lambda (|Z_\lambda|^2+\eta_1^2)^{-2})|^2}
{\frac{\tau((|Z_\lambda|^2+\eta_1^2)^{-1})^2}{2}
\left(H_{\mu_1^{(\lambda)}}'(\iu\eta_1)
-1/H_{\mu_2}'(\iu\eta_2)\right)}\nonumber\\
&=\eta_1^2\tau(h^{-1}k^{-1})
+\frac{|\tau(Z_\lambda h^{-2})|^2}
{\frac{\tau(h^{-1})^2}{2}
\left(H_{\mu_1^{(\lambda)}}'(\iu\eta_1)-1/H_{\mu_2}'(\iu\eta_2)\right)}.
\label{eq:positive-density}
 \end{align}
By \eqref{eq:regularized-denominator-bound}, the denominator in the
second term of \eqref{eq:positive-density} is at least
$\tau(h^{-2})-\tau(h^{-1})^2$. Indeed, the additional term on the
right-hand side of \eqref{eq:regularized-denominator-bound} is
nonnegative, and Cauchy--Schwarz gives
\[
\tau((|X_2|^2+\eta_2^2)^{-1})^2
\le\tau((|X_2|^2+\eta_2^2)^{-2}).
\]
Under the assumption $\tau(h^{-2})-\tau(h^{-1})^2>0$, we can therefore
use \eqref{eq:regularized-CS} to obtain
\begin{align*}
\pi\rho_\varepsilon(\lambda)
&\le\eta_1^2\tau(h^{-1}k^{-1})
+\frac{\tau(h^{-1}k^{-1})
\left(\tau(h^{-1})-\eta_1^2\tau(h^{-2})\right)}
{\tau(h^{-2})-\tau(h^{-1})^2}\\
&\quad=
\frac{\tau(h^{-1}k^{-1})
\left(\tau(h^{-1})-\eta_1^2\tau(h^{-1})^2\right)}
{\tau(h^{-2})-\tau(h^{-1})^2}\\
&\quad=
\frac{\tau(h^{-1}k^{-1})\tau(h^{-1})
\left(1-\eta_1^2\tau(h^{-1})\right)}
{\tau(h^{-2})-\tau(h^{-1})^2}\\
&\quad\le
\frac{\tau(h^{-2})\tau(h^{-1})}
{\tau(h^{-2})-\tau(h^{-1})^2}\\
&\quad=
\tau(h^{-1})\left(
1+\frac{\tau(h^{-1})^2}
{\tau(h^{-2})-\tau(h^{-1})^2}
\right),
\end{align*}
where we used 
\[
0\le1-\eta_1^2\tau(h^{-1})\le1,
\qquad
0\le\tau(h^{-1}k^{-1})\le\tau(h^{-2}).
\]
The second inequality above follows from 
Cauchy--Schwarz and equality of the spectral distributions
of $h$ and $k$. 

If instead
$H_{\mu_2}'(\iu\eta_2)\ge0$, the correction term in
\eqref{eq:regularized-density} is nonpositive, and
\begin{align*}
    0\le\pi\rho_\varepsilon(\lambda)
&\le\eta_1^2\tau((|Z_\lambda|^2+\eta_1^2)^{-1}(|Z_\lambda^*|^2+\eta_1^2)^{-1})=\eta_1^2\tau(h^{-1}k^{-1})\\
&\le\eta_1^2\tau((|Z_\lambda|^2+\eta_1^2)^{-2})\le\tau((|Z_\lambda|^2+\eta_1^2)^{-1}).
\end{align*}
Thus \eqref{eq:regularized-density-bound} holds in both cases.
\end{proof}

We next apply Lemma~\ref{lem:denominator} near a point
$\lambda_0\in F_1\setminus F$. If $|X_2|$ is not
scalar, the last term in \eqref{eq:regularized-denominator-bound}
is bounded below by a positive constant times $\tau(h^{-1})^4$.
\begin{lem}
\label{lem:denominator-F1}
For every $\lambda_0\in F_1\setminus F$, there are a neighborhood $V$
of $\lambda_0$ and $\varepsilon_0>0$ such that
\begin{equation}
H_{\mu_2}'(\iu\eta_2)<0,\qquad
\tau((|Z_\lambda|^2+\eta_1^2)^{-1})\le\frac2{m_2(|X_2|)}
\label{eq:F1-sign-resolvent}
\end{equation}
for all $\lambda\in V$ and $0<\varepsilon<\varepsilon_0$.
If $|X_2|$ is not scalar, then, after shrinking $V$ and
$\varepsilon_0$, there is $c>0$ such that
\begin{align}
&\frac{\tau((|Z_\lambda|^2+\eta_1^2)^{-1})^2(H_{\mu_1^{(\lambda)}}'(\iu\eta_1)H_{\mu_2}'(\iu\eta_2)-1)}{2H_{\mu_2}'(\iu\eta_2)}
\nonumber\\
&\qquad\ge\tau((|Z_\lambda|^2+\eta_1^2)^{-2})-\tau((|Z_\lambda|^2+\eta_1^2)^{-1})^2+c\,\tau((|Z_\lambda|^2+\eta_1^2)^{-1})^4
\label{eq:F1-denominator-bound}
\end{align}

throughout $V\times(0,\varepsilon_0)$.
\end{lem}
\begin{proof}
Since $\lambda_0\in F_1\setminus F$, we have
$0<m_2(|X_2|)<\infty$. By Corollary~\ref{cor:omega-boundary}, we have
\begin{equation}
\label{joint-cont-etas}
\begin{aligned}
&\bigl(\eta_1(\lambda,\varepsilon),\eta_2(\lambda,\varepsilon)\bigr)
\longrightarrow(0,\infty)\\
&\qquad\text{as }(\lambda,\varepsilon)\longrightarrow(\lambda_0,0),
\quad\varepsilon>0.
\end{aligned}
\end{equation}
In particular, for every $a,b>0$, there exists some neighborhood $V$ of
$\lambda_0$ and $\varepsilon_0>0$ such that
\[
\eta_1(\lambda,\varepsilon)<a,\qquad
\eta_2(\lambda,\varepsilon)>b
\]
for all $\lambda\in V$ and $0<\varepsilon<\varepsilon_0$.
By Lemma \ref{lem:increasing function} and \eqref{eqn:deritative-HT} (or a direct calculation),
we can deduce
\begin{equation}
-\iu yH_{\mu_2}(\iu y)\to m_2(|X_2|),\qquad
y^2H_{\mu_2}'(\iu y)\to-m_2(|X_2|)\quad(y\to\infty).
\label{eq:H2-large-height}
\end{equation}
Thus we can choose $b>0$ so that, for $y>b$,
\[
H_{\mu_2}'(\iu y)<0,\qquad
-\iu yH_{\mu_2}(\iu y)\ge\tfrac12m_2(|X_2|).
\]
By joint continuity of subordination functions (Theorem \ref{thm:dw-continuity}), choose $V$ and
$\varepsilon_0$ so that $\eta_1<1$ and $\eta_2>b$ throughout
$V\times(0,\varepsilon_0)$. Subordination gives
\begin{align*}
\eta_1\eta_2&\ge\eta_2(\eta_1-\varepsilon)
=-\iu\eta_2H_{\mu_2}(\iu\eta_2)\ge\tfrac12m_2(|X_2|),\\
\tau((|Z_\lambda|^2+\eta_1^2)^{-1})&=\frac1{\eta_1(\eta_1+\eta_2-\varepsilon)}
\le\frac1{\eta_1\eta_2}\le\frac2{m_2(|X_2|)}.
\end{align*}
This proves \eqref{eq:F1-sign-resolvent}. In particular,
$H_{\mu_2}'(\iu\eta_2)$ is real-valued and tends to zero from below, whereas its scaled
limit in \eqref{eq:H2-large-height} is strictly negative.

Suppose $|X_2|$ is not scalar. Since $\mu_{|X_2|^2}$ has at least two points in its support, there are
intervals $I,J\subset[0,M]$ of positive $\mu_{|X_2|^2}$-mass with
$d=\operatorname{dist}(I,J)>0$. For $s\in I$ and $t\in J$,
\[
\left|\frac{1}{s+y^2}-\frac{1}{t+y^2}\right|
=\frac{|s-t|}{(s+y^2)(t+y^2)}
\ge \frac{d}{(M+y^2)^2}.
\]
By the variance identity, restricting the integral to
$(I\times J)\cup(J\times I)$ gives
\begin{align*}
&\tau((|X_2|^2+y^2)^{-2})
-\tau((|X_2|^2+y^2)^{-1})^2\\
&\quad=\frac12\iint_{[0,\infty)^2}
\left(\frac{1}{s+y^2}-\frac{1}{t+y^2}\right)^2
\di{\mu_{|X_2|^2}(s)}\di{\mu_{|X_2|^2}(t)}
\ge \frac{c}{(M+y^2)^4},
\end{align*}
where $c=\mu_{|X_2|^2}(I)\mu_{|X_2|^2}(J)d^2>0$.
Using \eqref{joint-cont-etas}, shrink $V$ and $\varepsilon_0$ so
that $\eta_2>\max\{1,\sqrt M\}$ and $\varepsilon_0\le1/2$.
Since $\eta_j>\varepsilon$,
\[
\frac{\eta_2(\eta_2-\varepsilon)}{\eta_1(\eta_1-\varepsilon)}
\ge\frac{\eta_2^2}{2\eta_1^2}.
\]
The last term on the right of \eqref{eq:regularized-denominator-bound}
is therefore bounded below, after changing $c$, by
\[
c\frac{\tau((|Z_\lambda|^2+\eta_1^2)^{-1})^2}{(\eta_1\eta_2)^2}
\ge c\,\tau((|Z_\lambda|^2+\eta_1^2)^{-1})^4,
\]
because $\tau((|Z_\lambda|^2+\eta_1^2)^{-1})^{-1}=\eta_1(\eta_1+\eta_2-\varepsilon)
\ge\eta_1\eta_2$. This proves \eqref{eq:F1-denominator-bound}.
\end{proof}


\begin{prop}
\label{prop:local-ac-F1}
For every $\lambda_0\in\partial\Omega\cap(F_1\setminus F)$, there
exist some neighborhood $V$ of $\lambda_0$ and constants $C,\varepsilon_0>0$
such that
\begin{equation}
\rho_\varepsilon(\lambda)\le C,\qquad
\lambda\in V,\quad 0<\varepsilon<\varepsilon_0.
\label{eq:uniform-density-F1}
\end{equation}
Consequently, $\mu_X|_V\ll\di{A}$.
\end{prop}
\begin{proof}
By Lemma~\ref{lem:denominator-F1}, choose $V$ and $\varepsilon_0$
so that $H_{\mu_2}'(\iu\eta_2)<0$ and $\tau((|Z_\lambda|^2+\eta_1^2)^{-1})$ is uniformly
bounded on $V\times(0,\varepsilon_0)$. By
Corollary \ref{cor:omega-boundary}, we may also assume $\eta_1<1$ on $V\times(0,\varepsilon_0)$.

First suppose $|X_1-\lambda_0|$ is not scalar. Choose two compact
intervals $I,J\subset[0,M]$, separated by a distance $d>0$, carrying
positive mass for $\mu_{|X_1-\lambda_0|^2}$. By weak continuity, after
shrinking $V$ there is $p>0$ such that both intervals have mass at least
$p$ for every $\mu_{|X_1-\lambda|^2}$ with $\lambda\in V$. Thus
\begin{align*}
&\tau((|Z_\lambda|^2+\eta_1^2)^{-2})-\tau((|Z_\lambda|^2+\eta_1^2)^{-1})^2\\
&\quad=\frac12\iint_{\mathbb R_+^2}
\left(\frac1{s+\eta_1^2}-\frac1{t+\eta_1^2}\right)^2
\di{\mu_{|X_1-\lambda|^2}(s)}\di{\mu_{|X_1-\lambda|^2}(t)}\\
&\quad\ge\frac{p^2d^2}{(M+\eta_1^2)^4}
\ge\frac{p^2d^2}{(M+1)^4}.
\end{align*}
Together with \eqref{eq:regularized-density-bound} and the uniform bound
on $\tau((|Z_\lambda|^2+\eta_1^2)^{-1})$, this proves \eqref{eq:uniform-density-F1} in this case.

Now suppose $|X_1-\lambda_0|=a1$ for some $a>0$. Then $X_1-\lambda_0$ is then $a$ times a unitary.
After shrinking $V$ and $\varepsilon_0$, we have
\[
a/2\le|X_1-\lambda|\le3a/2,\qquad \eta_1<a/2.
\]
Consequently, for constants $c,C>0$,
\begin{align*}
    &c\le\tau((|Z_\lambda|^2+\eta_1^2)^{-1})\le C,\\
&\tau((|Z_\lambda|^2+\eta_1^2)^{-1}(|Z_\lambda^*|^2+\eta_1^2)^{-1})\le C,\\
&|\tau(Z_\lambda (|Z_\lambda|^2+\eta_1^2)^{-2})|\le C
\end{align*}
uniformly on $V\times(0,\varepsilon_0)$.

 Since $|X_1-\lambda_0|=a1$, we have
$m_{-2}(|X_1-\lambda_0|)=a^{-2}$. 
Thus $\lambda_0\in F_1$ implies $m_2(|X_2|)\le a^2$.
Suppose that this inequality is strict. The operator
$X_1-\lambda_0$ is invertible, so $X_1-\lambda$ remains
invertible for $\lambda$ near $\lambda_0$, and the map
$\lambda\longmapsto
\frac{1}{\tau(|X_1-\lambda|^{-2})}$
is continuous there. Consequently,
\[
m_2(|X_2|)
<
\frac{1}{\tau(|X_1-\lambda|^{-2})}
\]
throughout some neighborhood of $\lambda_0$.
This neighborhood is contained in $F_1$ and hence is
disjoint from $\Omega$, contradicting
$\lambda_0\in\partial\Omega$. Therefore
$m_2(|X_2|)=a^2$. 
Moreover, $|X_2|$ cannot be scalar. Otherwise the preceding
equality would imply $|X_2|=a1$, and hence
\[
m_2(|X_1-\lambda_0|)
=a^2
=\frac{1}{m_{-2}(|X_2|)}.
\]
Thus $\lambda_0\in F_2$ as well as $F_1$, so
$\lambda_0\in F$, contrary to our assumption.

By \eqref{eq:F1-denominator-bound} and the lower bound on
$\tau((|Z_\lambda|^2+\eta_1^2)^{-1})$, there is $c'>0$ such that
\[
\frac{\tau((|Z_\lambda|^2+\eta_1^2)^{-1})^2}{2}
\left(H_{\mu_1^{(\lambda)}}'(\iu\eta_1)
-\frac1{H_{\mu_2}'(\iu\eta_2)}\right)\ge c'
\]
throughout a possibly smaller $V\times(0,\varepsilon_0)$.
The positive-term density formula \eqref{eq:positive-density},
together with the preceding bounds, proves
\eqref{eq:uniform-density-F1}.

Finally, for every nonnegative $f\in C_c(V)$, weak convergence gives
\[
\int f\di{\mu_X}
=\lim_{\varepsilon\downarrow0}\int f\rho_\varepsilon\di{A}
\le C\int f\di{A}.
\]
Hence $\mu_X|_V\le C\di{A}|_V$, and in particular
$\mu_X|_V\ll\di{A}$.
\end{proof}

We next treat $F_2$. The argument is the regularized version of the gradient estimate used above. By \eqref{eq:regularized-gradient} and the polar decomposition of $\lambda-X_1$,
\begin{align}
|\partial_\lambda\Phi_\varepsilon(\lambda)|
&\le\frac12\tau\left(
|\lambda-X_1|(|\lambda-X_1|^2+\eta_1^2)^{-1}
\right)
\le\frac{1}{4\eta_1}.
\label{eq:regularized-gradient-bound}
\end{align}
We shall also use a uniform bound near $F$. Indeed, if $\lambda_0\in F$, Lemma~\ref{lem:desciption of four parts} gives
$|X_1-\lambda_0|=c1$ for some $c>0$. Hence
$|X_1-\lambda|\ge(c/2)1$ whenever $|\lambda-\lambda_0|<c/2$, and therefore
\begin{align}
|\partial_\lambda\Phi_\varepsilon(\lambda)|\le\frac1c
\label{eq:regularized-gradient-near-F}
\end{align}
on this neighborhood, uniformly in $\varepsilon>0$.

\begin{prop}
\label{prop:F2-zero}
We have $\mu_X(F_2)=0$.
\end{prop}
\begin{proof}
We may assume that $F_2\neq\emptyset$. By
Lemma~\ref{lem:desciption of four parts}, we have
$F_2=\{\lambda\in\CC:|\lambda-z_0|\le R\}$ for some $z_0\in\CC$ and
$R\ge0$. For $d>0$, Green's theorem gives
\begin{align}
\mu_X^{(\varepsilon)}(\{|\lambda-z_0|<R+d\})
=\frac{1}{\pi\iu}\int_{|\lambda-z_0|=R+d}
\partial_\lambda\Phi_\varepsilon(\lambda)\di{\lambda}.
\label{eq:regularized-F2-flux}
\end{align}
Since the disk on the left is open, weak convergence gives
\begin{align}
\mu_X(F_2)
\le
\liminf_{\varepsilon\downarrow0}
\mu_X^{(\varepsilon)}(\{|\lambda-z_0|<R+d\}).
\label{eq:F2-portmanteau}
\end{align}

First suppose that $F=F_1\cap F_2=\emptyset$. Given $M>0$, the boundary limits in
Corollary~\ref{cor:omega-boundary}, together with Theorem~\ref{thm:dw-continuity}, imply
$\eta_1(\lambda,\varepsilon)\to\infty$ as
$(\lambda,\varepsilon)\to(z,0)$ for each $z\in F_2$.
By compactness of $F_2$, there exist $d_0,\varepsilon_0>0$ such that
\[
\eta_1(\lambda,\varepsilon)>M
\quad\text{whenever}\quad
\dist(\lambda,F_2)<d_0,\qquad 0<\varepsilon<\varepsilon_0.
\]
Fix $0<d<\min\{d_0,1\}$. The circle $|\lambda-z_0|=R+d$ has length
at most $2\pi(R+1)$. Thus \eqref{eq:regularized-gradient-bound} and
\eqref{eq:regularized-F2-flux} give
\[
\mu_X^{(\varepsilon)}(\{|\lambda-z_0|<R+d\})
\le\frac{R+1}{2M},\qquad 0<\varepsilon<\varepsilon_0.
\]
For this fixed $d$, let $\varepsilon\downarrow0$ and use
\eqref{eq:F2-portmanteau}. We obtain
$\mu_X(F_2)\le(R+1)/(2M)$. Letting $M\to\infty$ proves
$\mu_X(F_2)=0$.

\begin{figure}[htbp]
\centering
\resizebox{0.9\textwidth}{!}{%
\begin{tikzpicture}[
    scale=1,
    >=Stealth,
    every node/.style={font=\small},
    nearF/.style={draw=red!75!black,fill=red!8},
    farF/.style={draw=teal!75!black,line width=1.4pt},
    neararc/.style={draw=red!75!black,line width=1.8pt}
]
\def\Rad{2.2}
\def\Out{2.7}
\def\Del{0.7}

\coordinate (O) at (0,0);
\coordinate (P) at (40:\Rad);
\coordinate (Q) at (155:\Rad);

\filldraw[
    fill=blue!5,
    draw=blue!65!black,
    line width=.8pt
] (O) circle[radius=\Rad];

\draw[farF] (O) circle[radius=\Out];

\filldraw[nearF] (P) circle[radius=\Del];
\filldraw[nearF] (Q) circle[radius=\Del];

\begin{scope}
    \clip (P) circle[radius=\Del];
    \draw[neararc] (O) circle[radius=\Out];
\end{scope}
\begin{scope}
    \clip (Q) circle[radius=\Del];
    \draw[neararc] (O) circle[radius=\Out];
\end{scope}

\fill (O) circle[radius=1.2pt];
\node[below left] at (O) {$z_0$};

\fill[red!75!black] (P) circle[radius=1.5pt];
\fill[red!75!black] (Q) circle[radius=1.5pt];
\node[below right,text=red!75!black] at (P) {$F$};
\node[below right,text=red!75!black] at (Q) {$F$};

\node[text=blue!65!black,font=\large]
    at (0,.65) {$F_2=\overline{D(z_0,R)}$};

\draw[gray]
    (O) -- (-35:\Rad)
    node[midway,below] {$R$};

\draw[
    {Stealth[length=3pt]}-{Stealth[length=3pt]}
]
    (-35:\Rad) -- (-35:\Out)
    node[midway,sloped,above=3pt,inner sep=1pt] {$d$};

\draw[red!75!black]
    (P) -- ++(110:\Del)
    node[midway,left] {$\delta$};

\node[text=teal!75!black,align=center]
    at (0,-3.05)
    {$|\lambda-z_0|=R+d$};

\node[
    anchor=north west,
    text width=6.0cm,
    align=left,
    text=red!75!black
] (neartext) at (3.3,2.55) {%
    \textbf{Near \(F\): short arcs}\\[3pt]
    Total length $\le C\delta$, and
    $|\partial_\lambda\Phi_\varepsilon|\le C$.
};

\node[
    anchor=north west,
    text width=6.0cm,
    align=left,
    text=teal!75!black
] (fartext) at (3.3,.9) {%
    \textbf{Away from \(F\): small gradient}\\[3pt]
    Joint continuity gives $\eta_1>M$, hence
    \[
    |\partial_\lambda\Phi_\varepsilon|
    \le \frac{1}{4\eta_1}
    \le \frac{1}{4M}.
    \]
};

\draw[->,red!75!black]
    (neartext.west) to[bend right=12] (40:\Out);
\draw[->,teal!75!black]
    (fartext.west) to[bend left=12] (0:\Out);

\end{tikzpicture}%
}
\caption{%
The contour decomposition used to prove $\mu_X(F_2)=0$
in the case $R>0$ and $F\ne\emptyset$.
The red arcs lie in the $\delta$-neighborhood of $F$;
the remaining arcs are shown in teal.
}
\label{fig:F2-zero-strategy}
\end{figure}


Now suppose that $F\ne\emptyset$. If $R=0$, then
$F=F_2=\{z_0\}$. By \eqref{eq:regularized-gradient-near-F}, Green's
theorem, and weak convergence, for sufficiently small $\delta>0$,
\begin{align*}
\mu_X(\{z_0\})
\le
\liminf_{\varepsilon\downarrow0}
\mu_X^{(\varepsilon)}(\{|\lambda-z_0|<\delta\})
\le C\delta.
\end{align*}
Letting $\delta\downarrow0$ gives $\mu_X(F_2)=0$.

Assume therefore that $R>0$. Since $F$ is a finite set,
\eqref{eq:regularized-gradient-near-F} gives $\delta_0,C>0$ such that
$|\partial_\lambda\Phi_\varepsilon|\le C$ whenever
$\dist(\lambda,F)<\delta_0$, uniformly in $\varepsilon>0$. Fix
$0<\delta<\min\{\delta_0,R/4\}$. For $0<d<1$ and $0<\delta<R/4$, the portion of the circle
$|\lambda-z_0|=R+d$ lying within distance $\delta$ of $F$
has length at most $C\delta$, where $C$ is independent
of $d$ and $\delta$. Indeed, $F$ is finite, and each
disk of radius $\delta$ cuts out an arc of length
$O(\delta)$ on a circle of radius at least $R$.
Together with \eqref{eq:regularized-gradient-near-F},
this shows that the contribution of these arcs to
the absolute flux in \eqref{eq:regularized-F2-flux}
is at most $C\delta$.

For the remaining arcs, consider the compact set
\[
K=\{z\in F_2:\dist(z,F)\ge\delta/2\}\subset F_2\setminus F.
\]
If $d<\delta/2$, radial projection from $|\lambda-z_0|=R+d$ onto
$|\lambda-z_0|=R$ moves each point a distance $d$. Therefore every
point of the remaining arcs has distance at most $d$ from $K$.
If $K$ is empty, there are no remaining arcs.

Otherwise, given $M>0$, compactness of $K$, the boundary limits in
Corollary~\ref{cor:omega-boundary}, and joint continuity give
$d_0,\varepsilon_0>0$ such that
\[
\eta_1(\lambda,\varepsilon)>M
\quad\text{whenever}\quad
\dist(\lambda,K)<d_0,\qquad 0<\varepsilon<\varepsilon_0.
\]
Fix $0<d<\min\{1,\delta/2,d_0\}$. By
\eqref{eq:regularized-gradient-bound}, the contribution of these
remaining arcs to the absolute flux is at most
\[
\frac1\pi\,2\pi(R+1)\frac1{4M}=\frac{R+1}{2M}.
\]
Consequently, for this fixed $d$ and all sufficiently small
$\varepsilon>0$,
\[
\mu_X^{(\varepsilon)}(\{|\lambda-z_0|<R+d\})
\le\frac{4C|F|}{\pi}\delta+\frac{R+1}{2M}.
\]
Letting $\varepsilon\downarrow0$ and using \eqref{eq:F2-portmanteau},
we obtain
\[
\mu_X(F_2)\le\frac{4C|F|}{\pi}\delta+\frac{R+1}{2M}.
\]
The constants on the right are independent of $\delta$ and $M$.
For fixed $\delta$, let $M\to\infty$, and then let $\delta\downarrow0$.
This proves $\mu_X(F_2)=0$.
\end{proof}

\begin{proof}[Proof of Theorem~\ref{thm:no-boundary-singular}]
By Theorem~\ref{thm:main-sec3}, the measure $\mu_X$ is supported on
$\overline\Omega$ and is absolutely continuous on $\Omega$.
Every point of $\partial\Omega\cap(F_1\setminus F)$ has a neighborhood
on which $\mu_X$ is absolutely continuous, by
Proposition~\ref{prop:local-ac-F1}. These neighborhoods admit a countable subcover of
$\partial\Omega\cap(F_1\setminus F)$. Consequently, $\mu_X$ is
absolutely continuous on their union.

By Proposition~\ref{prop:F2-zero}, $\mu_X(F_2)=0$; in particular,
$\mu_X(F)=0$. The decomposition in
Lemma~\ref{lem:desciption of four parts} therefore shows that the
singular part of $\mu_X$ is supported on $S$. Since $S$ is finite,
this singular part is purely atomic. Thus there exists a nonnegative
$\rho_X\in L^1(\di{A})$ such that
\[
\mu_X=\rho_X\di{A}
+\sum_{\lambda\in S}\mu_X(\{\lambda\})\delta_\lambda.
\]
This proves \eqref{eq:main}.
\end{proof}

\subsection{The atomic part}
We next turn to the atomic part of $\mu_X$. 
We recall some basic results on the generalized singular values of operators \cite{FackKosaki}.
For $t>0$ and $T\in \tilde{\mathcal{A}}$, let
\begin{align*}
\mu_t(T)=\inf\left\{x\ge0:\mu_{|T|}((x,\infty))<t\right\}
\end{align*}
be the generalized singular numbers of $T$. We shall use
\begin{align}
\mu_{s+t}(S+T)&\le\mu_s(S)+\mu_t(T),
\qquad s,t>0,\label{eq:singular-sum}\\|\mu_t(T+B)-\mu_t(T)|&\le\|B\|,
\qquad t>0,\quad B\in\mathcal M.
\label{eq:bounded-singular-perturbation}
\end{align}
Moreover, for every Borel function
$f\colon[0,\infty)\to[-\infty,\infty]$ for which the integrals are
well defined, we have 
\begin{align}
\tau(f(|T|))=\int_0^1 f(\mu_t(T))\,dt.
\label{eq:singular-trace}
\end{align}
For $T\in\mathrm{Log}^+(\tau)$ and $0<q<1$, we denote
\begin{align*}
J_q(T)=\int_0^q\log\mu_t(T)\,dt.
\end{align*}
In particular, we have 
\[
\log\Delta(T-\lambda)=\tau(\log|T-\lambda|)=J_1(T-\lambda).
\]

For $\lambda_0\in S$, we denote 
\[
a=\mu_{|X_1-\lambda_0|}(\{0\}), \quad \text{and}\quad b=\mu_{|X_2|}(\{0\}).
\]
We note that $a+b\geq 1$ by Lemma \ref{lem:desciption of four parts}.
The main result of this section is the following. 
\begin{thm}\label{thm:atom-main}
Let $\lambda_0\in S$, we have  
\begin{align}
\mu_X(\{\lambda_0\})=\lim_{r\downarrow0}\frac{1}{2\pi\log r}\int_0^{2\pi}J_b(X_1-\lambda_0-re^{\iu\theta})\,d\theta.
\label{eq:atom-formula}
\end{align}
Moreover,
\begin{align}
a+b-1\le\mu_X(\{\lambda_0\})\le\frac{a+b-1}{a}\,\mu_{X_1}(\{\lambda_0\}).
\label{eq:atomic-sharp-bounds}
\end{align}
In particular, $\mu_X(\{\lambda_0\})=0$ when $a+b=1$.
If $X_1$ is normal, then 
\[
\mu_X(\{\lambda_0\})=a+b-1.
\]
\end{thm}

One can show that the function $\lambda\mapsto J_q(T-\lambda)$ is subharmonic for any $0<q\leq 1$ (see \cite{RV2005Lelong} for a proof for $T$ bounded). Motivated by theory of Lelong numbers of plurisubharmonic functions (see \cite{Demailly, Rashkovskii} for example), we introduce the following definition. 

\begin{defn}
 For $T\in\mathrm{Log}^+(\tau)$, $\lambda\in\mathbb{C}$ and $0\leq t\leq 1$, the Lelong number of $J_t(T-\lambda)$ at $\lambda$ is defined by
 \begin{equation}\label{eq:alphacircular}
\alpha_T(\lambda,t)
=\lim_{r\downarrow0}\frac1{2\pi\log r}
\int_0^{2\pi}J_t(T-\lambda-re^{i\theta})d\theta.
\end{equation}
\end{defn}

It is known that the Lelong number of a subharmonic function is the atom of its corresponding Riesz measure (see \cite{Berndtsson} for example).
We present this classic result in the context of Brown measure and provide a self-contained proof. 
\begin{lem}\label{lem:circle}
For $T\in\mathrm{Log}^+(\tau)$ and $\lambda\in\CC$, we have
\begin{align}
    \mu_T(\{\lambda \})=\alpha_T(\lambda,1)=\lim_{r\downarrow0}\frac1{2\pi\log r}
\int_0^{2\pi}\tau(\log|T-\lambda-re^{i\theta}|)d\theta.
\end{align}
\end{lem}

\begin{proof}
By Jensen's formula, for every $z\in\CC$ and $r>0$,
\begin{align*}
\frac1{2\pi}\int_0^{2\pi}\log|z-re^{\iu\theta}|\,d\theta=\log\max\{r,|z|\}.
\end{align*}
Thus Fubini's theorem yields
\begin{align*}
\frac1{2\pi}\int_0^{2\pi}\tau(\log|T-\lambda-re^{i\theta}|)\,d\theta&=\frac1{2\pi}\int_0^{2\pi}\int_\CC\log|\zeta-\lambda-re^{\iu\theta}|\,d\mu_T(\zeta)\,d\theta\\
&=\int_\CC\left(\frac1{2\pi}\int_0^{2\pi}\log|\zeta-\lambda-re^{\iu\theta}|\,d\theta\right)d\mu_T(\zeta)\\
&=\int_\CC\log\max\{r,|\zeta-z|\}\,d\mu_T(\zeta).
\end{align*}
Hence, for $0<r<1$ we have 
\begin{align*}
&\frac1{2\pi\log{r}}\int_0^{2\pi}\tau(\log|T-\lambda-re^{i\theta}|)\,d\theta\\
&\qquad=
 \int_{\{|\zeta-\lambda|\le1\}}
 \frac{\log\max\{r,|\zeta-\lambda|\}}{\log r}\,d\mu_T(\zeta)
 +\frac{1}{\log r}
 \int_{\{|\zeta-\lambda|>1\}}\log|\zeta-\lambda|\,d\mu_T(\zeta).
\end{align*}
For $\{|\zeta-w|\le1\}$, we have
\[
0\le
\frac{\log\max\{r,|\zeta-w|\}}{\log r}
\le1.
\]
As $r\downarrow0$, this quotient is identically $1$ at
$\zeta=\lambda$ and tends to $0$ at every $\zeta\ne w$.
Dividing by $\log r$ and applying dominated convergence as $r\downarrow0$, we have
\begin{align*}
    \mu_T(\{\lambda \})=\lim_{r\downarrow0}\frac1{2\pi\log r}
\int_0^{2\pi}\tau(\log|T-\lambda-re^{i\theta}|)d\theta.
\end{align*}
This completes the proof. 
\end{proof}

We next study the Lelong number $\alpha_T(\lambda,t)$ for any $0\leq t\leq 1$. 
\begin{thm}\label{thm:cutoff-profile}
Let $T\in {\rm Log}^{+}(\tau)$, $\lambda\in\mathbb{C}$ and $t\in [0,1]$. 
Put $f(t)=\alpha_T(\lambda,t)$, 
$\kappa=\tau(P_{\ker(T-\lambda)})$ and 
$m=\mu_T(\{\lambda\})$. 
Then the following statements are valid.
\begin{enumerate}
\item The function $f$ is nonnegative, nondecreasing and convex on $[0,1]$;
$f(0)=0$, $f(1)=m$, and $f$ is continuous on $[0,1)$.
\item $0\le \kappa\le m\le1$. If $\kappa=0$, then
\begin{equation}\nonumber
f(t)=0\quad(0\le t<1),\qquad f(1)=m.
\end{equation}
\item If $\kappa>0$, then for every $0\le t\le1$,
\begin{equation}\label{eq:kernelbounds}
(t+\kappa-1)_+\le f(t)
\le\frac{(t+\kappa-1)_+}{a}\,m.
\end{equation}
In particular, $f(t)>0$ if and only if $t+\kappa>1$. Moreover,
\begin{equation}\nonumber
f(t)-f(s)\ge t-s\quad\text{whenever }1-\kappa\le s<t\le1.
\end{equation}
\item The bound in terms of Brown mass is
\begin{equation}\label{eq:brownupper}
0\le f(t)\le (t+m-1)_+\quad(0\le t\le1).
\end{equation}
Thus $f(t)=0$ for $t\le1-m$.
\item If $\kappa=m$, then
\begin{equation}\label{eq:exact}
f(t)=(t+m-1)_+\quad(0\le t\le1).
\end{equation}
This applies in particular when $T$ is normal.
\end{enumerate}
\end{thm}
\begin{proof}
Proof for (1): 
    for $0<s<t\le1$, decreasing singular values give
\[
J_s(B)\ge\frac{s}{t}J_t(B).
\]
Indeed, the average of the decreasing function $v\mapsto\log s_v(B)$
on $(0,s)$ is at least its average on $(0,t)$. Applying the inequality to $\lambda+re^{i\theta}$ and
and divide by $\log r<0$, we obtain 
\[
\frac1{2\pi\log r}\int_0^{2\pi}J_s(T-\lambda-re^{i\theta})d\theta
\le\frac{s}{t}\,
\frac1{2\pi\log r}\int_0^{2\pi}J_t(T-\lambda-re^{i\theta})d\theta.
\]
Letting $r\downarrow0$ gives
\begin{equation}\label{eq:monotonicity}
0\le f(s)\le\frac{s}{t}f(t)\le f(t)
\qquad(0<s<t\le1).
\end{equation}
Together with $f(0)=0$, this proves that $f$ is nondecreasing on
$[0,1]$. Taking $t=1$ also gives $0\le f(s)\le sm$. The decreasing singular values imply that $t\mapsto J_t(w)$
is concave, with $J_0(w)=0$. By the definition of $f(t)=\alpha_T(\lambda,t)$, then $f$ is convex on $[0,1]$.
Finite convex functions are continuous on $(0,1)$. Note that $f$ is also right continuous at $0$ by the inequality $0\leq f(s)\leq sm$. This proves (1).

Proof for (2): 
It is clear that $0\leq \kappa\leq m\leq 1$. If $\kappa=0$, then $s_v(T-\lambda)>0$ for any $0<v<1$. For $|z|<s_v(T-\lambda)/2$, then \eqref{eq:bounded-singular-perturbation} implies that
\[
 s_v(T-\lambda-z)\geq s_v(T-\lambda)-|z|\geq s_t(T-\lambda)-|z|\geq \frac{1}{2}s_t(T-\lambda)
\]
for any $0<v<t$. In addition, we have 
\[
 \log s_v(T-\lambda-z)\leq \log(1+|z|+s_v(T-\lambda)).
\]
It follows that 
\[
    t\log (s_t(T-\lambda)/2) \leq J_t(T-\lambda-z)\leq C<\infty. 
\]
where $C$ can be chosen to be independent of $z$ in some neighborhood of $\lambda$. This implies that 
\[
\left |\frac{J_t(T-\lambda-re^{i\theta}}{\log r}\right|\leq \frac{\max\{t\log (s_t(T-\lambda)/2), C\}}{\log r}\rightarrow 0
\]
as $r\rightarrow 0$. It follows that $f(t)=0$ for any $0\leq t<1$. By Lemma \ref{lem:circle}, we have $f(1)=m=\mu_T(\{\lambda\})$. 

Proof for (3): If $\kappa=0$, then $\kappa\leq m$ by the positivity of $m$; and if $\kappa=1$, then $T=\lambda$ so $m=1$ and $f(t)=t$. We now suppose that $0<\kappa<1$, then $s_v(T-\lambda)>0$ for any $0<v<1-\kappa$. A similar argument as for (2) yields that $f(t)=0$ for any $0\leq t\leq 1-\kappa$. In particular, we have $f(1-\kappa)=0$. For $1-\kappa\leq t\leq 1$, the convexity property implies that 
\[
  f(t)\leq \frac{1-t}{\kappa}f(1-\kappa)+\frac{\kappa+t-1}{\kappa}f(1)=\frac{(t+\kappa-1)m}{\kappa}.
\]
On the other hand, we have $s_v(T-\lambda)=0$ for $v\geq 1-\kappa$. By \eqref{eq:bounded-singular-perturbation} we have 
\[
  s_v(T-\lambda-z)\leq |z|.
\]
Hence for $t\geq 1-\kappa$ (equivalently, $t+\kappa-1\geq 0$) and $|z|=r<1$, by splitting the logarithmic integral at $1-\kappa$, we obtain
\begin{align*}
    J_t(T-\lambda-z)
&\le\int_0^{1-\kappa}\log^+s_v(T-\lambda-z)d v
+(t+\kappa-1)\log |z|\\
&\le C+(t+\kappa-1)\log |z|.
\end{align*}
Taking the average over the circle $\{re^{i\theta}\}$ and division by $\log r<0$, by the definition of $f(t)=\alpha_T(\lambda,t)$, we obtain
\[
   f(t)\geq t+\kappa-1. 
\]
When $t=1$, this proves $\kappa\leq f(1)=m$. 

For $1-\kappa<s<t\leq 1$, we have $f(1-\kappa)=0$ and $f(s)\geq s+\kappa-1$ by the first part of (3). Moreover, the convexity of $f$ implies that 
\[
  \frac{f(t)-f(s)}{t-s}\geq \frac{f(s)-f(1-\kappa)}{s+\kappa-1}\geq 1.
\]
This shows that $f(t)-f(s)\geq t-s$. For $s=1-\kappa$, we also have
\[
 f(t)-f(1-\kappa)=f(t)\geq t+\kappa-1=t-(1-\kappa).
\]
This completes the proof. 

Proof for (4) and (5): If $\kappa=0$, then \eqref{eq:exact} follows from Part (2). If $\kappa>0$ and $t\leq 1-\kappa$, then $f(t)=0$. If $t>1-\kappa$, then $\kappa\leq m$ and hence
\[
 f(t)\leq \frac{(t+\kappa-1)m}{a\leq m-(1-t)}=t+m-1. 
\]
Thus \eqref{eq:exact} holds for any $t\in [0,1]$. If $\kappa=m$, then \eqref{eq:exact} follows from Part (3). If $T$ is normal, then we have $\kappa=m$. This finishes the proof. 
\end{proof}

We recall that $J_1(T-\lambda)=\tau(\log|T-\lambda|))$
and $J_q(T-\lambda)=\int_0^q\log\mu_t(T-\lambda)\,dt$. 
\begin{lem}\label{lem:atomic-potential-comparison}
Fix $\lambda_0\in S$. For every $0<\varepsilon<1$, there exist $C_\varepsilon>0$ such that, for all sufficiently small $z\ne0$,
\begin{align}
&(1-\varepsilon)J_{b}(X_1-\lambda_0-z)+\varepsilon J_1(X_1-\lambda_0-z)-C_\varepsilon\nonumber\\
&\qquad\le J_1(X-\lambda_0-z)
\le \tau(\log(2+|X-\lambda_0|))+J_b(X_1-\lambda_0-z),
\label{eq:atomic-potential-comparison}
\end{align}
where $b=\mu_{|X_2|}(\{0\})$.
\end{lem}
\begin{proof}
We will derive the upper bound uses only singular-value inequalities.
By applying \eqref{eq:bounded-singular-perturbation}, for every $|z|\leq 1$, we have 
\begin{align*}
   \log^+\mu_t(X-\lambda_0-z)
   &\leq \log(1+\mu_t(X-\lambda_0-z))\\
   &\leq \log (2+\mu_t(X-\lambda_0)).
\end{align*}
Hence, by \eqref{eq:singular-trace} we have 
\begin{align}
\sup_{|z|\le1}\int_0^1\log^+\mu_t(X-\lambda_0-z)\,dt\nonumber
 &\leq \int_0^1\log(2+\mu_t(X-\lambda_0))\,dt\\
 &=\tau(\log(2+|X-\lambda_0|)<\infty.
 \label{inequality-in-proof(4.33)}
\end{align}
Since $\mu_s(X_2)=0$ for $s>1-b$, applying \eqref{eq:singular-sum} gives, for each $0<t<b$ and all sufficiently small $\delta>0$,
\begin{align*}
\mu_{t+1-b}(X_1+X_2-\lambda_0-z)&\le\mu_{t-\delta}(X_1-\lambda_0-z)+\mu_{1-b+\delta}(X_2)\\
&=\mu_{t-\delta}(X_1-\lambda_0-z).
\end{align*}
Letting $\delta\downarrow0$ and using left-continuity, we obtain the desired inequality.
By combining \eqref{eq:singular-trace} and recalling that $X=X_1+X_2$, we obtain
\begin{align*}
&J_1(X-\lambda_0-z)=\int_0^1\log\mu_t(X-\lambda_0-z)dt\\
 &=\int_0^{1-b}\log \mu_t(X-\lambda_0-z)dt
    +\int_0^{b}\log\mu_{t+1-\mu_{|X_2|}(\{0\})}(X-\lambda_0-z)dt \\
&\le\int_0^{1-b}\log^+\mu_t(X-\lambda_0-z)\,dt+J_{b}(X_1-\lambda_0-z)\\
&\le \tau(\log(2+|X-\lambda_0|))+J_{b}(X_1-\lambda_0-z),
\end{align*}
withe $b=\mu_{|X_2|}(\{0\})$.
This proves the upper bound. 

For the lower bound, we will use both the subordination formula and some singular value inequalities. 
Lemma~\ref{lem:desciption of four parts} shows that $\lambda_0+z\in\Omega$ for all sufficiently small $z\ne0$. The inequality \eqref{eq:bounded-singular-perturbation} implies $\mu_{|X_1-\lambda_0-z|}\Rightarrow\mu_{|X_1-\lambda_0|}$ weakly.  
For the rest of this proof, we denote $\eta_j=\eta_j(\lambda_0+z)>0$. Hence Theorem~\ref{thm:dw-continuity} and Proposition~\ref{prop:omega(0)}(3) give $\eta_1,\eta_2\to0$ as $z\to0$.
The subordination relation \eqref{eqn:subordination-lambda-index} (by taking imaginary part and letting $z\rightarrow {0}$) yields
\begin{align}
\frac{\eta_2}{\eta_1+\eta_2}=\eta_2^2 \cdot \tau((|X_2|^2+\eta_2^2)^{-1}),
\label{eq:subordination}
\end{align}
where we used the identity
\begin{equation}
 \label{eqn:eta's-ratio}
      -\Im G_{|X_2|}(\omega_2^{(\lambda)}(0))=\eta_2 \cdot \tau((|X_2|^2+\eta_2^2)^{-1}).
\end{equation}
By the dominated convergence theorem, we obtain 
\begin{align*}
\frac{\eta_2}{\eta_1+\eta_2}=\int_0^\infty \frac{\eta_2^2}{u^2+\eta_2^2}d\mu_{|X_2|}(u)\longrightarrow\mu_{|X_2|}(\{0\})
\end{align*}
when $z\rightarrow {0}$ (hence $\eta_2=\eta_2(\lambda_0+z)\rightarrow {0}$). For $0<y<1$, we then have 
\begin{align*}
\frac12\tau(\log(|X_2|^2+y^2))
&\ge(\mu_{|X_2|}(\{0\})+\mu_{|X_2|}((0,d)))\log y+\mu_{|X_2|}([d,\infty))\log d\\
&\ge(\mu_{|X_2|}(\{0\})+\varepsilon(1-\mu_{|X_2|}(\{0\})))\log y+\log d\\
&=(b+\varepsilon(1-b))\log y +\log d,
\end{align*}
where we choose $d\in(0,1)$ such that $\mu_{|X_2|}((0,d))\le\varepsilon(1-\mu_{|X_2|}(\{0\}))$. 
Using this estimate with $y=\eta_2$ and Proposition~\ref{prop:calculating the logarithmic integral}, we obtain
\begin{align}
&\tau\bigl(\log|X_1+X_2-\lambda_0-z|\bigr)\notag\\
 &=\frac{1}{2}\tau(\log(|X_1-\lambda-z|^2+\eta_1^2))
   +\frac{1}{2}\tau(\log(|X_2|^2+\eta_2^2))
 -\log(\eta_1+\eta_2)\notag\\
 &\ge
\frac12\tau\bigl(\log(|X_1-\lambda_0-z|^2+\eta_1^2)\bigr)
+(b+\varepsilon(1-b))\log\eta_2\nonumber\\
&\qquad\qquad-\log(\eta_1+\eta_2)+\log d\nonumber\\
&\ge
\frac12\tau\bigl(\log(|X_1-\lambda_0-z|^2+\eta_1^2)\bigr)
-(1-\varepsilon)(1-b)\log\eta_2
-C_\varepsilon.
\label{eq:atomic-regularized-lower}
\end{align}
where 
\[
  C_\varepsilon=\log\eta_2-\log(\eta_1+\eta_2)-\log {d}.
\]
By \eqref{eqn:eta's-ratio}, we note that $\frac{\eta_2}{\eta_1+\eta_2}\rightarrow 1-b$ as $z\rightarrow 0$. 

Moreover, we have 
\begin{align*}
&\frac12\tau(\log(|X_1-\lambda_0-z|^2+\eta_1^2))\\
&\quad\ge\int_0^{\mu_{|X_2|}(\{0\})}\log\mu_t(X_1-\lambda_0-z)\,dt+(1-\mu_{|X_2|}(\{0\}))\log\eta_1\\
&=J_b(X_1-\lambda_0-z)+(1-b)\log\eta_1.
\end{align*}
We also have
\begin{align*}
\frac12\tau(\log(|X_1-\lambda_0-z|^2+\eta_1^2))\ge\tau(\log|X_1-\lambda_0-z|)=J_1(X_1-\lambda_0-z).
\end{align*}
By taking linear combination of the above two inequalities, 
it follows that for $0<\varepsilon<1$ we obtain
\begin{align*}
    &\frac12\tau(\log(|X_1-\lambda_0-z|^2+\eta_1^2))\\
     &\geq (1-\varepsilon)(J_b(X_1-\lambda_0-z)+(1-b)\log\eta_1)
     +\varepsilon J_1(X_1-\lambda_0-z).
\end{align*}
By substituting it into \eqref{eq:atomic-regularized-lower}, we prove the lower bound.
\end{proof}

We are now ready to prove our main theorem in this section. 
\begin{proof}[Proof of Theorem \ref{thm:atom-main}]
We apply Theorem \ref{thm:cutoff-profile} for $f(t)=\alpha_{X_1}(\lambda_0,t)$.
By taking average over the circle $\{re^{i\theta}\}$ and division by $\log r<0$, the inequality in \eqref{eq:atomic-potential-comparison} implies that 
\[
 (1-\varepsilon)\alpha_{X_1}(\lambda_0, b)+\varepsilon \alpha_{X_1}(\lambda_0, 1)\geq   \alpha_X(\lambda_0, 1)\geq \alpha_{X_1}(\lambda_0,b).
\]
where we used the fact that $\tau(\log(2+|X-\lambda_0|))<\infty$. 
Lemma \ref{lem:circle} implies that $\alpha_{X_1}(\lambda_0, 1)=\mu_{X_1}(\{\lambda_0\})<\infty$. By taking $\varepsilon\rightarrow 0$, we obtain
\[
   \mu_X(\{\lambda_0 \})=  \alpha_X(\lambda_0,1)=\alpha_{X_1}(\lambda_0, b)
\]
with $b=\mu_{|X_2|}(\{0\})$. Hence by Theorem \ref{thm:cutoff-profile}, we obtain
\[
a+b-1\le\mu_X(\{\lambda_0\})\le\frac{a+b-1}{a}\,\mu_{X_1}(\{\lambda_0\})
\]
with $a=\tau(P_{\ker(X_1-\lambda)})=\mu_{|X_1-\lambda_0|}(\{0\})$.
If $X_1$ is normal, then both sides are identical and hence 
$\mu_X(\{\lambda_0\})=a+b-1$.
\end{proof}

The following example shows that both inequalities in \eqref{eq:atomic-sharp-bounds} can be strict.

\begin{example}
\label{exa:atomic-bounds-zero-block}
Let $X_1$ have the $*$-distribution of
\begin{align*}
\begin{pmatrix}
0&1&0\\
0&0&0\\
0&0&0
\end{pmatrix}
\end{align*}
with respect to the normalized trace, and let $X_2$ be $R$-diagonal and $*$-free from $X_1$, with $\mu_{|X_2|}=b\delta_0+(1-b)\delta_1$, where $1/3<b<1$. Set $X=X_1+X_2$. Observe that $0\in S$. For $r=|\lambda|>0$, the singular values of $X_1-\lambda$ are $r,s_+(r),s_i(r)$, where
\begin{align*}s
s_+(r)=\frac{\sqrt{1+4r^2}+1}{2},\quad s_-(r)=\frac{\sqrt{1+4r^2}-1}{2}.
\end{align*}
Since $s_+(r)s_-(r)=r^2$, we obtain
\begin{align*}
J_b(X_1-\lambda)=\begin{cases}
(b-\frac13)\log r+\frac13\log s_+(r),& \frac13<b\le\frac23,\\
(2b-1)\log r+(1-b)\log s_+(r),& \frac23<b<1.
\end{cases}
\end{align*}
The formula \eqref{eq:atom-formula} yields
\begin{align*}
\mu_X(\{0\})=
\begin{cases}
b-\frac13,& \frac13<b\le\frac23,\\
2b-1,& \frac23<b<1.
\end{cases}
\end{align*}
For $b=1/2$, we have
\begin{align*}
a+b-1=\frac16=\mu_X(\{0\})<\frac14=\frac{a+b-1}{a}\,\mu_{X_1}(\{0\}).
\end{align*}
Thus the lower bound is attained and the upper bound is strict. For $b=3/4$, we deduce that
\begin{align*}
a+b-1=\frac5{12}<\mu_X(\{0\})=\frac12<\frac58=\frac{a+b-1}{a}\,\mu_{X_1}(\{0\}).
\end{align*}
Thus both bounds are strict.
\end{example}

\section{Examples\label{sec:Examples}}

Most examples below were already known. We show how the techniquest
developed above may lead to effective calculation. Examples \ref{exa:R-diag+0}
and \ref{exa:Addition to a circular Cauchy operator} are in \cite{HA-SC-brown},
where some of the techniques we use were first introduced. Example
\ref{exa:nilpotent+Haar} illustrates the fact that the addition of
freely independent variables does not lead to a new convolution operation
on probability measures on $\mathbb{C}$. More precisely, the Brown
$\mu_{X_{1}+X_{2}}$ does not depend solely on $\mu_{X_{1}}$ and
$\mu_{X_{2}}$ for general $*$-free variables $X_{1}$ and $X_{2}$,
even when $X_{2}$ is $R$-diagonal. 

The calculations require determining first the set $\Omega$ and setting
up the fixed point equation for $\omega_{1}^{(\lambda)}(0)$. In case
the result is sufficiently explicit (for instance, in Examples \ref{exa:R-diag+0}--\ref{exa:Addition to a circular operator}),
one can calculate the density of $\mu_{X_{1}+X_{2}}$ by further differentiation
of the expression in Lemma \ref{lem:magic simplification}. Otherwise,
we use Lemma \ref{lemma:lasciate ogne speranza voi ch'entrate} and
calculate $\overline{\partial}_{\lambda}\omega_{1}^{(\lambda)}(0)$
in order to determine the second term in the formula provided therein. 
\begin{example}
\label{exa:R-diag+0}\textbf{Brown measure of an} $R$-\textbf{diagonal
operator}. Observe that the measure $\mu_{X_{1}+X_{2}}$ requires
knowledge of $\mu_{X_{2}^{*}X_{2}},$in addition to informatin about
$X_{1}$. The particular case in which $X_{1}=0$ reflects the fact
that the Brown measure of an arbitrary $R$-diagonal operator $X$
is determined by $\mu_{X^{*}X}$. An actual calculation was first
done in \cite{ha-larsen} (in the bounded case) and \cite{HA-SC-brown}
(in general). That calculation can also be derived from our approach,
as we indicate briefly.

Fix an arbitrary $R$-diagonal operator $X$, and write it as $X=X_{1}+X_{2}$,
where $X_{1}=0$ and $X_{2}=X$. According to Definition \ref{def:DSF1F2F},
we have $S=\{0\},$ 
\[
F_{1}=\{\lambda:|\lambda|^2\ge m_{2}(X)\},\quad F_{2}=\{\lambda:|\lambda|^2\le1/m_{-2}(X)\},
\]
and thus 
\[
\Omega=\{\lambda:1/m_{-2}(X)\le|\lambda|^2\le m_{2}(X)\}.
\]
Furthermore, 
\[
\mu_{1}^{(\lambda)}=\frac{1}{2}(\delta_{|\lambda|}+\delta_{-|\lambda|}),\text{ and }\mu_{2}=\widetilde{\mu}_{X},
\]
so 
\[
H_{\mu_{1}^{(\lambda)}}(z)=-\frac{|\lambda|^{2}}{z},
\]
and 
\[
\psi_{\lambda}(z)=H_{\mu_{2}}(H_{\mu_{1}^{(\lambda)}}(z))=\frac{1}{G_{\mu_{2}}(-|\lambda|^{2}/z)}+\frac{|\lambda|^{2}}{z}.
\]
For $\lambda\in\Omega$, $\omega_{1}^{(\lambda)}(0)$ is the fixed
point of the function $\psi_{\lambda}$. Using the notation

\[
M(z)={\it \int_{\mathbb{R}}\frac{zt}{1-zt}\,d\mu_{X^{*}X}(t)}
\]
for the moment generating function of $\mu_{X^{*}X}$, the equation
$\psi_{\lambda}(\omega_{1}^{(\lambda)}(0))=\omega_{1}^{(\lambda)}(0)$
is rewritten a
\begin{align*}
\frac{1}{|\lambda|^{2}-\omega_{1}^{(\lambda)}(0)^{2}}=\frac{1+M(\omega_{1}^{(\lambda)}(0)^{2}/|\lambda|^{4})}{|\lambda|^{2}},
\end{align*}
which implies 
\begin{equation}
M(\omega_{1}^{(\lambda)}(0)^{2}/|\lambda|^{4})=\frac{\omega_{1}^{(\lambda)}(0)^{2}}{|\lambda|^{2}-\omega_{1}^{(\lambda)}(0)^{2}}.\label{eq:example1}
\end{equation}
Next we use Lemma \ref{lem:magic simplification} and (\ref{eq:example1})
to calculate 
\begin{align*}
2\partial_{\lambda}\tau(\log(|X-\lambda)|) & =\frac{\overline{\lambda}}{|\lambda|^{2}-\omega_{1}^{(\lambda)}(0)^{2}}\\
 & =\overline{\lambda}\frac{1+M(\omega_{1}^{(\lambda)}(0)^{2}/|\lambda|^{4})}{|\lambda|^{2}}.
\end{align*}
This equation can be simplified using the $S$-transform of $\mu_{X^{*}X},$that
we denote by $S$, which satisfies 
\[
S(M(z))=z\frac{M(z)+1}{M(z)}.
\]
Combined with (\ref{eq:example1}), this yields 
\[
S(M(\omega_{1}^{(\lambda)}(0)^{2}/|\lambda|^{4}))=\frac{1}{|\lambda|^{2}},
\]
or, equivalently, 
\[
S^{\langle-1\rangle}(1/|\lambda|^{2})=M(\omega_{1}^{(\lambda)}(0)^{2}/|\lambda|^{4}).
\]
Thus, 
\[
2\partial_{\lambda}\tau(\log(|X-\lambda)|)=\frac{1+S^{\langle-1\rangle}(1/|\lambda|^{2})}{\lambda},
\]
and applying $2\partial_{\overline{\lambda}}/\pi$ we obtain the formula
\[
-\frac{(S^{\langle-1\rangle})'(1/|\lambda|^{2})}{\pi|\lambda|^{4}}
\]
for the density of $\mu_{X}$at $\lambda\in\Omega$. This, of course,
agrees with \cite[Theorem 4.17]{HA-SC-brown}. (See \cite{R-diag-revisit}
for a more detailed exposition of this derivation of the formula for
$\mu_{X}$.) 
\end{example}

\begin{example}
\textbf{Nilpotent plus Haar unitary}. \label{exa:nilpotent+Haar}For
this example, we assume that $X_{1}$ is distributed like an $n\times n$
complex matrix $A$ and $X_{2}$ is a Haar unitary. Later, we specialize
to $n=2$ and $A$ nilpotent. (The case in which $n=2$ and $A$ is
selfadjoint with eigenvalues $\pm1$ was considered in \cite{biane-lehner})
As in the preceding example, we have $H_{\mu_{2}}(z)=-1/z$. Suppose
that $\lambda\in\mathbb{C}$, and denote by $P_{\lambda}(z)=\det(z-|\lambda-A|^{2})$,
and by $a_{j}(\lambda)$, $j=1,\dots,n$ the eigenvalues of $|\lambda-A|$.
We have 
\[
P_{\lambda}(z)=\prod_{j=1}^{n}(z-a_{j}(\lambda)^{2})=z^{n}-t(\lambda)z^{n-1}+\cdots+(-1)^{n-1}s(\lambda)z+(-1)^{n}d(\lambda),
\]
where $t(\lambda)$ and $d(\lambda)$ are the trace and determinant
of $|\lambda-A|^{2}$, respectively, while $s(\lambda)$ is the sum
of all products of $n-1$ distinct $a_{j}(\lambda)$. Then
\[
G_{\mu_{1}^{(\lambda)}}(z)=\frac{1}{2n}\sum_{j=1}^{n}\left[\frac{1}{z-a_{j}(\lambda)}+\frac{1}{z+a_{j}(\lambda)}\right]=\frac{z}{n}\sum_{j=1}^{n}\frac{1}{z^{2}-a_{j}(\lambda)^{2}}=\frac{zP'_{\lambda}(z^{2})}{nP_{\lambda}(z^{2})},
\]
therefore 
\[
H_{\mu_{1}^{(\lambda)}}(z)=\frac{nP_{\lambda}(z^{2})-z^{2}P'_{\lambda}(z^{2})}{zP'_{\lambda}(z^{2})},
\]
and
\begin{align*}
H_{\mu_{2}}(H_{\mu_{1}^{(\lambda)}}(z)) & =\frac{zP'_{\lambda}(z^{2})}{z^{2}P'_{\lambda}(z^{2})-nP_{\lambda}(z^{2})}\\
 & =\frac{nz^{2n-1}-(n-1)t(\lambda)z^{2n-3}+\cdots+(-1)^{n-1}s(\lambda)z}{t(\lambda)z^{2n-2}+\cdots+(-1)^{n-1}nd(\lambda)}.
\end{align*}
Observe that $0$ and $\infty$ are fixed points of this function
for every $\lambda$. However, $\infty$ is the Denjoy-Wolff point
only when $t(\lambda)\le n$, demonstrating again that $F_{2}$ is
either empty of a closed disk (possibly of zero radius). Also, $0$
is the Denjoy-Wolff point when $s(\lambda)\le nd(\lambda)$, and the
boundary of $F_{1}$ is described by the real-algebraic equation $s(\lambda)=nd(\lambda)$.
When neither of these inequalities is satisfied, $\omega_{1}^{(\lambda)}(0)$
is the unique solution in $\mathbb{C}^{+}$ of a polynomial equation
of degree $2n-2$. 

In the special case in which $n=2$, the above formula specializes
to
\[
(H_{\mu_{2}}\circ H_{\mu_{1}^{(\lambda)}})(z)=\frac{2z^{3}-t(\lambda)z}{t(\lambda)z^{2}-2d(\lambda)}.
\]
We have 
\[
\omega_{1}^{(\lambda)}(0)=i\sqrt{\frac{t(\lambda)-2d(\lambda)}{t(\lambda)-2}}.
\]
when $t(\lambda)>2$ and $t(\lambda)>2d(\lambda)$. The measure $\mu_{X_{1}+X_{2}}$
can be calculated quite explicitly when $A$ is a nilpotent operator.
Suppose that 
\[
A=\left[\begin{array}{cc}
0 & \varepsilon\\
0 & 0
\end{array}\right]
\]
for some $\varepsilon>0$. (Every nonzero nilpotent in $M_{2}(\mathbb{C})$
is unitarily equivalent to such an operator $A$.) In this case, we
have 
\[
t(\lambda)=2|\lambda|^{2}+\varepsilon^{2},\quad d(\lambda)=|\lambda|^{4},
\]
and thus $\omega_{j}^{(\lambda)}(0)$ belongs to $\mathbb{C}^{+}$
provided that $2|\lambda|^{2}+\varepsilon^{2}>2$ and $2|\lambda|^{2}+\varepsilon^{2}>2|\lambda|^{4}$.
The second condition amounts to 
\[
|\lambda|<\sqrt{\frac{1+\sqrt{1+2\varepsilon^{2}}}{2}},
\]
while the first condition is vacuous if $\varepsilon>\sqrt{2}$ and
it amounts to 
\[
|\lambda|>\sqrt{\frac{2-\varepsilon^{2}}{2}}
\]
for $\varepsilon\le\sqrt{2}$. Without any further calculation, we
may conclude that the support of $\mu_{X_{1}+X_{2}}$ is an annulus
(or disk for $\varepsilon\ge\sqrt{2})$ that depends on $\varepsilon$.
In particular, $\mu_{X_{1}+X_{2}}$ does not depend just on $\mu_{X_{1}}$
and $\mu_{X_{2}}$, in contrast with the case of free selfadjoint
variables $X_{1}$ and $X_{2}$. We also observe that as $\varepsilon\downarrow0$
the support of $\mu_{X_{1}+X_{2}}$ converges, as may be expected,
to the unit circle.

We proceed to calculate the density of $\mu_{X_{1}+X_{2}}$, using
the restatement of \ref{lem:magic simplification} as 
\[
2\partial_{\lambda}\tau(\log(|X-\lambda|))={\rm tr}_{2}((\lambda-A)^{*}(|\lambda-A|^{2}-\omega_{1}^{(\lambda)}(0)^{2})^{-1}),
\]
where ${\rm tr}_{2}$ denotes the normalized trace on $M_{2}(\mathbb{C})$,
and 
\[
\omega_{1}^{(\lambda)}(0)^{2}=-\frac{2|\lambda|^{2}+\varepsilon^{2}-2|\lambda|^{4}}{2|\lambda|^{2}+\varepsilon^{2}-2}
\]
for $\lambda$ in the annulus (or punctured disk) 
\[
\Omega_{\varepsilon}=\left\{ \lambda\in\mathbb{C}:\sqrt{\left(\frac{2-\varepsilon^{2}}{2}\right)_{+}}<|\lambda|<\sqrt{\frac{1+\sqrt{1+2\varepsilon^{2}}}{2}}\right\} .
\]
Since 
\[
|\lambda-A|^{2}-\omega_{1}^{(\lambda)}(0)^{2}=\left[\begin{array}{cc}
|\lambda|^{2}-\omega_{1}^{(\lambda)}(0)^{2} & -\varepsilon\overline{\lambda}\\
-\varepsilon\lambda & |\lambda|^{2}+\varepsilon^{2}-\omega_{1}^{(\lambda)}(0)^{2}
\end{array}\right],
\]
we have 
\[
(|\lambda-A|^{2}-\omega_{1}^{(\lambda)}(0)^{2})^{-1}=\frac{1}{D_{\varepsilon}(\lambda)}\left[\begin{array}{cc}
|\lambda|^{2}+\varepsilon^{2}-\omega_{1}^{(\lambda)}(0)^{2} & \varepsilon\overline{\lambda}\\
\varepsilon\lambda & |\lambda|^{2}-\omega_{1}^{(\lambda)}(0)^{2}
\end{array}\right],
\]
where 
\[
D_{\varepsilon}(\lambda)=(|\lambda|^{2}-\omega_{1}^{(\lambda)}(0)^{2})^{2}-\varepsilon^{2}\omega_{1}^{(\lambda)}(0)^{2}
\]
is the determinant of $|\lambda-A|^{2}-\omega_{1}^{(\lambda)}(0)^{2}$.
Thus, 
\begin{align*}
(\lambda-A)^{*}(|\lambda-A|^{2}-\omega_{1}^{(\lambda)}(0)^{2})^{-1}\\
=\frac{1}{D_{\varepsilon}(\lambda)}\left[\begin{array}{cc}
\overline{\lambda}(|\lambda|^{2}+\varepsilon^{2}-\omega_{1}^{(\lambda)}(0)^{2}) & \varepsilon\overline{\lambda}^{2}\\
-\varepsilon(|\lambda|^{2}+\varepsilon^{2}-\omega_{1}^{(\lambda)}(0)^{2}) & -\varepsilon^{2}\overline{\lambda}+\overline{\lambda}(|\lambda|^{2}-\omega_{1}^{(\lambda)}(0)^{2})
\end{array}\right],
\end{align*}
and we find that 
\[
2\partial_{\lambda}\tau(\log(|X-\lambda|))=\overline{\lambda}f_{\varepsilon}(|\lambda|^{2}),
\]
where $f(z)$ is a rational function of $z$ such that 
\[
f(|\lambda|^{2})=\frac{|\lambda|^{2}-\omega_{1}^{(\lambda)}(0)^{2}}{(|\lambda|^{2}-\omega_{1}^{(\lambda)}(0)^{2})^{2}-\varepsilon^{2}\omega_{1}^{(\lambda)}(0)^{2}}.
\]
It follows that, in the annulus $\Omega_{\varepsilon}$, the measure
$\mu_{X_{1}+X_{2}}$ has density 
\[
\frac{1}{\pi}\overline{\partial}_{\lambda}(\overline{\lambda}f(\lambda\overline{\lambda}))=\frac{1}{\pi}\left[f(|\lambda|^{2})+|\lambda|^{2}f'(|\lambda|^{2})\right],
\]
where $f'$ is the usual derivative of $f$. The integral of this
density relative to area measure is 
\begin{align*}
2\int_{\sqrt{((2-\varepsilon^{2})/2})_{+}}^{\sqrt{(1+\sqrt{1+2\varepsilon^{2}})/2}}(f(r^{2})+r^{2}f'(r^{2}))r\,dr & =\int_{((2-\varepsilon^{2)}/2)_{+}}^{(1+\sqrt{1+2\varepsilon^{2}})/2}(f(\rho)+\rho f'(\rho))\,d\rho\\
 & =\rho f(\rho)|_{((2-\varepsilon^{2)}/2)_{+}}^{(1+\sqrt{1+2\varepsilon^{2}})/2}.
\end{align*}
This quantity is easily evaluated to be $1$ if we observe that $\omega_{1}^{(\lambda)}(0)=0$
when $|\lambda|=\sqrt{(1+\sqrt{1+2\varepsilon^{2}})/2}$. Thus, $\mu_{X}$
is absolutely continuous. 
\end{example}

\begin{example}
\label{exa:X_2=00003DHaar}\textbf{Addition to a free Haar unitary.}
As in the preceding example, $X_{2}$ is a Haar unitary operator,
but $X_{1}$ is an arbitrary variable $*$-free from $X_{2}$. As
seen above, $H_{\mu_{2}}(z)=H_{\widetilde{\mu}_{|X_{2}|}}(z)=-1/z$,
so $\omega_{1}^{(\lambda)}(0)$ is the Denjoy-Wolff point ot the map
$w\mapsto\psi_{\lambda}(w)=-1/H_{\mu_{1}^{(\lambda)}}(w)$. In particular,
for $\lambda\in\Omega$, $\omega_{1}^{(\lambda)}(0)$ is the unique
solution $w\in\mathbb{C}^{+}$ of the equation 
\[
w=-\frac{1}{H_{\mu_{1}^{(\lambda)}}(w)}.
\]
We use this equation (rather than Lemma \ref{lem:delbar psi_lambda})
in order to calculate $\overline{\partial}_{\lambda}\omega_{1}^{(\lambda)}(0)$
and the second term in the formula of Lemma \ref{lemma:lasciate ogne speranza voi ch'entrate}.
Using the notation 
\[
Y(\lambda)=(\omega_{1}^{(\lambda)}(0)^{2}-|X_{1}-\lambda|^{2})^{-1},
\]
it is easily verified that 
\[
G_{\mu_{1}^{(\lambda)}}(\omega_{1}^{(\lambda)}(0))=\omega_{1}^{(\lambda)}(0)\tau(Y(\lambda)).
\]
Recalling that $H_{\mu_{1}^{(\lambda)}}(w)=(1/G_{\mu_{1}^{(\lambda)}}(w))-w$,
the fixed point equation reduces to 
\[
\omega_{1}^{(\lambda)}(0)^{2}=1+\frac{1}{\tau(Y(\lambda))},
\]
and thus 
\begin{align*}
2\omega_{1}^{(\lambda)}(0)\overline{\partial}_{\lambda}\omega_{1}^{(\lambda)}(0) & =-\frac{\overline{\partial}_{\lambda}\tau(Y(\lambda))}{Y}\\
 & =\frac{\tau[(\overline{\partial}_{\lambda}Y(\lambda))Y(\lambda)^{2}]}{\tau(Y(\lambda))}\\
 & =\frac{\tau[(2\omega_{1}^{(\lambda)}(0)\overline{\partial}_{\lambda}\omega_{1}^{(\lambda)}(0)-(\lambda-X_{1}))Y(\lambda)^{2}]}{\tau(Y(\lambda))}.
\end{align*}
Solving for $\overline{\partial}_{\lambda}\omega_{1}^{(\lambda)}(0)$
yields 
\[
\overline{\partial}_{\lambda}\omega_{1}^{(\lambda)}(0)=\frac{\tau((\lambda-X_{1})Y(\lambda)^{2})}{2\omega_{1}^{(\lambda)}(0){\rm Var}(Y(\lambda))},
\]
where ${\rm Var}(Y(\lambda))=\tau(Y(\lambda)^{2})-\tau(Y(\lambda))^{2}$.
The second term in Lemma \ref{lemma:lasciate ogne speranza voi ch'entrate}
becomes 
\[
\frac{1}{\pi}\frac{|\tau((\lambda-X_{1})Y(\lambda)^{2})|^{2}}{{\rm Var}(Y(\lambda))},
\]
and this shows that the density of $\mu_{X_{1}+X_{2}}$ is positive
in the entire open set $\Omega$. 
\end{example}

\begin{example}
\label{exa:Addition to a circular operator}\textbf{Addition to a
circular operator}. Consider now a circular operator $X_{2}$ with
variance $t$, and an arbitrary variable $X_{1}$, $*$-free from
$X_{2}$. In this case, $\tau(\text{ker}(X_{2}))=0$, $m_{2}(|X_{2}|)=t$
and $1/m_{-2}(|X_{2}|)=0$. Hence, 
\[
\Omega=\left\{ \lambda:m_{-2}(|X_{1}-\lambda|)>1/t\right\} .
\]
Since $\mu_{2}$ is the semicircular distribution with variance $t$,
we have 
\begin{equation}
\frac{1}{G_{\mu_{2}}(z)}+tG_{\mu_{2}}(z)=z.\label{eq14-example3}
\end{equation}
Using equations (\ref{eq:G_jcompose with omega_j_}) and (\ref{eq:omega1+omega2}),
we obtain
\begin{align*}
-tG_{\mu_{2}}(\omega_{1}^{(\lambda)}(0)) & =-tG_{\mu_{2}}(\omega_{2}^{(\lambda)}(0))\\
 & =\frac{1}{G_{\mu_{2}}(\omega_{2}^{(\lambda)}(0))}-\omega_{2}^{(\lambda)}(0)=\omega_{1}^{(\lambda)}(0),
\end{align*}
which can be rewritten as 
\begin{equation}
\tau(\omega_{1}^{(\lambda)}(0)^{2}-|X_{1}-\lambda|^{2})^{-1})=-1/t.\label{eqn:fixEqnCircular}
\end{equation}
Denote $Y(\lambda)=(\omega_{1}^{(\lambda)}(0)^{2}-|X_{1}-\lambda|^{2})^{-1}$.
Applying implicit differentiation, we have 
\[
\overline{\partial}_{\lambda}\omega_{1}^{(\lambda)}(0)=\frac{\tau((\lambda-X_{1})Y(\lambda)^{2})}{2\omega_{1}^{(\lambda)}(0)\cdot\tau(Y(\lambda)^2)}.
\]
The second term in Lemma \ref{lemma:lasciate ogne speranza voi ch'entrate}
becomes 
\[
\frac{1}{\pi}\frac{|\tau((\lambda-X_{1})Y(\lambda)^{2})|^{2}}{\tau(Y(\lambda)^2)}.
\]
This formula was proved in \cite[Theorem 4.2]{pzh-ajm} (see also
\cite{bordenave-caputo-chafai} for a special case and in \cite{byzh-ADV}
for the unbounded case). It is shown in \cite[Section 7]{byzh-ADV}
that $\mu_{X_{1}+X_{2}}$ is absolutely continuous in $\mathbb{C}$
and the density is bounded by $1/{\pi}t$. 
\end{example}

\begin{example}
\label{exa:Addition to a circular Cauchy operator}\textbf{Addition
to a circular Cauchy operator}. Suppose that $c_{1}$ and $c_{2}$
are {*}-free standard circular operators, set $X_{2}=c_{1}c_{2}^{-1}$,
and let $X_{1}$ be an arbitrary variable that is {*}-free from $X_{2}$.
It is shown in \cite{HA-SC-inv} that $\widetilde{\mu}_{|X_{2}|}$
is the standard Cauchy distribution, and hence $H_{\mu_{2}}(z)=i$
is constant. The fixed point equation shows that $\omega_{1}^{(\lambda)}(0)=i$
for every $\lambda\in\mathbb{C}$, in paricular, $\overline{\partial}_{\lambda}\omega_{1}^{(\lambda)}(0)=0$.
Lemma \ref{lemma:lasciate ogne speranza voi ch'entrate} shows now
that the density of $\mu_{X_{1}+X_{2}}$ is equal to 
\[
\frac{1}{\pi}\tau[(|\lambda-X_{1}|^{2}+1)^{-1}(|\overline{\lambda}-X_{1}^{*}|^{2}+1)^{-1}],\quad\lambda\in\mathbb{C}.
\]
This formula was first proved in \cite[Corollary 4.6]{HA-SC-inv}. 
\end{example}

\section{Application to free multiplicative Brownian motions}
\label{bzrw:section}
{
After an earlier version of this paper was posted on arXiv,
Driver--\allowbreak Hall--\allowbreak Ho--\allowbreak Kemp--\allowbreak Nemish--\allowbreak Nikitopoulos--\allowbreak Parraud
\cite[Section 6]{Seven} applied our subordination results to
multiplicative free random walks. They proved, among other results,
locally uniform convergence of the Brown measure densities in the
limiting interior and weak convergence of the Brown measures
\cite[Theorem 7.16 and Corollary 7.17]{Seven}.
In this section, we use our Fuglede--Kadison determinant formula
and their subordination convergence theorem to calculate the limit
of the logarithmic potentials. The convergence is locally uniform
in $(t,z)$ on the set $\mathcal I$ defined below. We identify the
limit with the potential calculated by Driver--Hall--Kemp
\cite{DHK} using the Hamilton--Jacobi method. We consider only
the case $u_0=I$, but a similar calculation works for arbitrary unitary $u_0$ which can recover potential formula obtained in \cite{Ho-Pi-jems}. 
}

Let $a_1,a_2,\ldots$ be freely independent copies of a bounded $R$-diagonal
operator with $\tau(a_1^*a_1)=1$. Put $B=a_1^*a_1$ and
\begin{equation}\label{bzrw:walk}
 \begin{gathered}
 B_k(t)=\prod_{j=1}^k(1+\sqrt{t/k}\,a_j),\qquad t>0,\\
 L_k(t,z)=\tau(\log|B_k(t)-z|),\quad S_k(t,z)=2L_k(t,z).
 \end{gathered}
\end{equation}
{
The factors are ordered in increasing $j$. Let $b_t$ be the free
multiplicative Brownian motion satisfying $db_t=b_t\,dc_t$, $b_0=I$,
where $c_t$ is a circular Brownian motion with $\tau(c_t^*c_t)=t$.
}

{
For $k\geq3$, consider the tracial probability space
$(\mathcal A_k,\tau_k)=(\mathcal A\overline\otimes M_k(\mathbb C),
\tau\otimes\mathrm{tr}_k)$, where $\mathrm{tr}_k$ is the normalized
matrix trace. Define
}
\begin{equation}\label{bzrw:linearization}
 Z_k=\sum_{j=1}^k1\otimes E_{j,j+1},\qquad
 A_k=\sum_{j=1}^ka_j\otimes E_{j,j+1},\qquad
 X_{k,t}=Z_k+\sqrt{t/k}\,A_k,
\end{equation}
where indices are taken modulo $k$.

\begin{lem}\label{bzrw:freeness}
The operators $Z_k$ and $A_k$ are $*$-free with respect to the
scalar-valued trace $\tau_k$, and $A_k$ has the same $*$-distribution
as $a_1$; in particular it is $R$-diagonal.
\end{lem}
{
This is the case for $u_0=I$ in \cite[Corollary 5.19]{Seven}. It can also be deduced by a
direct matrix-valued cumulant computation and the restriction and
factorization criteria of Nica--Shlyakhtenko--Speicher
\cite[Theorems 3.1 and 3.5]{NSS}.
}

For $\lambda\in\mathbb C$, write, as in Section 3,
\[
 \mu_{1,k}^{(\lambda)}=\widetilde\mu_{|Z_k-\lambda|},\qquad
 \mu_{2,k,t}=\widetilde\mu_{\sqrt{t/k}|a_1|},
\]
and denote the associated subordination functions by
$\omega_{1,k,t}^{(\lambda)},\omega_{2,k,t}^{(\lambda)}$.
The tilde denotes symmetrization. By Lemma~\ref{bzrw:freeness}
and the sum formula \eqref{eq:definition of mu1lambda and mu2},
\[
 \widetilde\mu_{|X_{k,t}-\lambda|}
   =\mu_{1,k}^{(\lambda)}\boxplus\mu_{2,k,t}.
\]
Recall that $m_p(T)=\tau(|T|^p)$, with the trace changed to $\tau_k$ when
$T\in\mathcal A_k$; negative moments are $+\infty$ when the kernel
is nonzero. Set
\[
 p=\tau(\ker a_1),\qquad \beta=m_{-2}(a_1)\in[1,+\infty],\qquad
 d_k(\lambda)=\frac1k\sum_{\zeta^k=1}|\zeta-\lambda|^{-2}.
\]
Here $d_k(\lambda)=+\infty$ at a $k$th root of unity, and $1/(+\infty)=0$.
{
We apply Definition \ref{def:DSF1F2F} to $X_1=Z_k$ and
$X_2=\sqrt{t/k}A_k$. Following \cite[Definition 6.1]{Seven}, we denote
the sets $F_1,F_2,F$ by $E_{1,k}(t),E_{2,k}(t),E_k(t)$,
respectively, and write $\Omega_k^{\mathrm{reg}}(t)$ for the set
$\Omega$ of Section~3. We use $\mathcal S_k$ for the exceptional set
to distinguish it from the potential $S_k(t,z)$. Thus
}
\begin{align}\label{bzrw:finite-sets}
 {\mathcal S_k}&=\{\lambda:\tau_k(\ker(Z_k-\lambda))+p\geq1\},\notag\\
 E_{1,k}(t)&=\{\lambda:t/k\leq1/d_k(\lambda)\},\notag\\
 E_{2,k}(t)&=\{\lambda:1+|\lambda|^2\leq t/(k\beta)\},\notag\\
 E_k(t)&=E_{1,k}(t)\cap E_{2,k}(t),\notag\\
 \Omega_k(t)&=\mathbb C\setminus(E_{1,k}(t)\cup E_{2,k}(t)),\qquad
 \Omega_k^{\mathrm{reg}}(t)=\Omega_k(t)\setminus {\mathcal S_k}.
\end{align}
{
Since $p<1$ and $\tau_k(\ker(Z_k-\lambda))$ equals $1/k$ at
a $k$th root of unity and zero elsewhere, $\mathcal S_k$ is empty
if $p+1/k<1$ and otherwise consists of all $k$th roots of unity.
}

\begin{lem}\label{bzrw:finite-domain}
The set $\Omega_k^{\mathrm{reg}}(t)$ is open. For $\lambda\ne0$,
\begin{equation}\label{bzrw:strict-domain}
 \lambda\in\Omega_k^{\mathrm{reg}}(t)\quad\Longleftrightarrow\quad
 \lambda\notin {\mathcal S_k},\quad
 d_k(\lambda)>k/t,\quad 1+|\lambda|^2>t/(k\beta).
\end{equation}
{
For $\lambda\ne0$, these conditions are also equivalent to
}
$\omega_{1,k,t}^{(\lambda)}(0),\omega_{2,k,t}^{(\lambda)}(0)
\in i(0,+\infty)$.
{
The conditions defining $\Omega_k^{\mathrm{reg}}(t)$ are unchanged
when $\lambda$ is multiplied by a $k$th root of unity.
}
\end{lem}
{
The lemma follows from parts (2) and (5) of
Lemma \ref{lem:desciption of four parts}, applied to $X_1=Z_k$ and
$X_2=\sqrt{t/k}A_k$. Here we use the fact that $E_k(t)\subset\{0\}$
for $k\geq3$.
}

For $R>0$ define
\begin{equation}\label{bzrw:K}
 K_\theta(R)=\frac{\sinh R}{2R(\cosh R-\cos\theta)}.
\end{equation}
The limit $K_\theta(0+)$ equals $[2(1-\cos\theta)]^{-1}$,
with value $+\infty$ at $\theta=0$ modulo $2\pi$.
If $K_\theta(0+)>1/t$, let $R_t(\theta)>0$ solve
$K_\theta(R_t(\theta))=1/t$; otherwise set $R_t(\theta)=0$.
Existence and uniqueness follow from Zhong \cite[Lemma 3.1]{Zhong}:
with $r=e^{-R}$,
$K_\theta(R)=T_{\cos\theta}(r)/2$ in that lemma's notation.
Thus $K_\theta$ strictly decreases from $K_\theta(0+)$ to zero.
{
The same monotonicity follows from the periodization identity
used in the proof of \cite[Lemma 7.14]{Seven}:
}
\begin{equation}\label{bzrw:period}
 K_\theta(R)=\sum_{n\in\mathbb Z}\frac1{R^2+(\theta+2\pi n)^2},
\end{equation}
and hence
\[
 K_\theta'(R)=-2R\sum_{n\in\mathbb Z}
 \frac1{[R^2+(\theta+2\pi n)^2]^2}<0.
\]
The series may be differentiated term by term, since it and its
$R$-derivative converge locally uniformly for $R>0$.
The implicit function theorem therefore shows that $R_t(\theta)$
is smooth in $(t,\theta)$ wherever $R_t(\theta)>0$.

We will use the sets
\begin{equation}\label{bzrw:domain}
 \Sigma_t=\{e^{\rho+i\theta}:|\rho|<R_t(\theta)\},\qquad
 \mathcal I=\{(t,z):t>0,\ z\in\Sigma_t\}.
\end{equation}
{
The set $\Sigma_t$ agrees with the domain in \cite{DHK}:
its positive boundary parameter $R=R_t(\theta)$ satisfies
$t=2R(\cosh R-\cos\theta)/\sinh R$.
}

{
We also consider the open set
}
\begin{equation}\label{bzrw:product-domain}
 \mathcal D_{k,t}=\{\lambda^k:\lambda\in\Omega_k^{\mathrm{reg}}(t),\ \lambda\ne0\}.
\end{equation}
We describe $\mathcal D_{k,t}$ in logarithmic polar coordinates.
For $z=re^{i\theta}=e^{\rho+i\theta}\ne0$ define
\begin{equation}\label{bzrw:Tk-polar}
 T_k(r,\theta)
 =k\frac{r^{2/k}-1}{r^2-1}(r^2-2r\cos\theta+1),\qquad
 T_k(1,\theta)=2(1-\cos\theta).
\end{equation}
Using the continuous extension at $r=1$, \eqref{bzrw:product-domain}
becomes
\begin{equation}\label{bzrw:polar-domain}
 \begin{split}
 \mathcal D_{k,t}=\{re^{i\theta}:{}&r>0,\quad T_k(r,\theta)<t,\quad
                  1+r^{2/k}>t/(k\beta)\}\setminus\mathcal E_k,\\
 \mathcal E_k={}&
 \begin{cases}\{1\},&p+1/k\geq1,\\ \varnothing,&p+1/k<1.
 \end{cases}
 \end{split}
\end{equation}
Here $k/d_k(\lambda)=T_k(r,\theta)$ when $\lambda^k=re^{i\theta}$,
by \cite[Lemma 6.6(1)]{Seven}. The additional exclusions disappear
when $k>t$ and $p+1/k<1$. In the linearization variable
$\lambda=se^{i\varphi}$, substitute $r=s^k$ and $\theta=k\varphi$.

To identify the limiting domain, put
\begin{equation}\label{bzrw:T-limit}
 T(e^\rho,\theta)=\frac{2\rho(\cosh\rho-\cos\theta)}{\sinh\rho},
 \qquad T(1,\theta)=2(1-\cos\theta).
\end{equation}
Then $\Sigma_t=\{re^{i\theta}:r>0,\ T(r,\theta)<t\}$, and
\begin{equation}\label{bzrw:Tk-factor}
 T_k(e^\rho,\theta)=a_k(\rho)T(e^\rho,\theta),\qquad
 a_k(\rho)=\frac{e^{2\rho/k}-1}{2\rho/k},\quad a_k(0)=1.
\end{equation}
Thus $T_k\to T$ locally uniformly on $\mathbb C\setminus\{0\}$.
For fixed $t$, $\mathcal D_{k,t}=\Sigma_k(1,t)$ of \cite{Seven}
for all large $k$. By \cite[Theorem 7.2 and Corollary 7.3]{Seven},
interior compact sets are eventually included, exterior compact
sets are eventually disjoint from the closures, and
$\overline{\mathcal D_{k,t}}\to\overline{\Sigma_t}$ in Hausdorff distance.
For the theorem below, we will verify the required inclusion of
compact subsets of $\mathcal I$ directly. 

\begin{thm}\label{bzrw:limit}
Every compact set $\mathcal K\Subset\mathcal I$ is contained in
$\{(t,z):z\in\mathcal D_{k,t}\}$ for all sufficiently large $k$.
Under the preceding hypotheses, locally uniformly for $(t,z)\in\mathcal I$,
with $z=e^{\rho+i\theta}$ and $R=R_t(\theta)$,
\begin{equation}\label{bzrw:result}
 S_k(t,z)\longrightarrow
 F_t(z):=\rho+\log[2(\cosh R-\cos\theta)]
                  -\frac{R^2-\rho^2}{t}
       =2\log\Delta_\tau(b_t-z).
\end{equation}
Equivalently, $L_k(t,z)\to F_t(z)/2=\tau(\log|b_t-z|)$.
Consequently $\Delta_\tau(B_k(t)-z)\to\Delta_\tau(b_t-z)$ locally
uniformly on $\mathcal I$.
\end{thm}

\begin{proof}
Fix $\mathcal K\Subset\mathcal I$ and write $z=e^{\rho+i\theta}$,
$R=R_t(\theta)$. All subsequent limits are uniform on $\mathcal K$.
{
The continuous functions $t,R$, and $R^2-\rho^2$ are strictly
positive on $\mathcal K$, and hence are bounded above and away from
zero there. Also, $|\rho|$ is bounded on $\mathcal K$.
}
Since $t-T(z)$ has a positive minimum, \eqref{bzrw:Tk-factor} and
\eqref{bzrw:polar-domain} imply $z\in\mathcal D_{k,t}$ on $\mathcal K$
for all large $k$; the additional exclusions disappear since $p<1$
and $t$ is bounded. Choose $\lambda_k=e^{(\rho+i\theta)/k}$ on local logarithmic charts.
The choice of root is immaterial.

Put $q_k=-i\omega_{1,k,t}^{(\lambda_k)}(0)$ and
$r_k=-i\omega_{2,k,t}^{(\lambda_k)}(0)$. By
Lemma~\ref{bzrw:finite-domain}, both numbers are finite and positive.
In the notation of \cite[Definition 6.3]{Seven},
\begin{equation}\label{bzrw:dictionary}
 \psi_k=iq_k,\qquad \eta_k=q_k^2,\qquad Q_k=kq_k.
\end{equation}
{
By \cite[Theorem 7.8 and Corollary 7.9]{Seven}, the following
convergence is uniform on $\mathcal K$, jointly in $t$ and $z$:
}
\begin{equation}\label{bzrw:input}
 Q_k\longrightarrow Q>0,\qquad
 \sum_{n\in\mathbb Z}\frac1{\rho^2+(\theta+2\pi n)^2+Q^2}=\frac1t.
\end{equation}
{
By \eqref{bzrw:period} and the uniqueness of $R$, we have
$Q^2=R^2-\rho^2$; see also \cite[Lemma 7.14]{Seven}.
} In particular $Q_k$ is bounded above and away
from zero on $\mathcal K$.

The cyclic determinant identity, with normalized matrix trace, is
\begin{equation}\label{bzrw:cyclic}
 S_k(t,z)=2\tau(\log|B_k(t)-z|)=k\tau_k(\log|X_{k,t}-\lambda_k|^2).
\end{equation}
Indeed, the blocks of $X_{k,t}^k$ are cyclic rotations of $B_k(t)$
and have equal determinants after subtracting $z\ne0$, by
$\Delta(1-CD)=\Delta(1-DC)$. Factoring $X_{k,t}^k-z$ proves
\eqref{bzrw:cyclic}, since conjugation by a diagonal unitary
multiplies $X_{k,t}$ by any $k$th root of unity.

Define $R_k>0$ by
\begin{equation}\label{bzrw:Rk}
 \cosh(R_k/k)=\cosh(\rho/k)+\frac{q_k^2}{2e^{\rho/k}}.
\end{equation}
The inequality $\cosh u-1\geq u^2/2$ gives $R_k=O(1)$, and
Taylor expansion yields $R_k^2=\rho^2+Q_k^2+O(k^{-1})$.
Thus $R_k\to R$. The identity
\[
 |\zeta-\lambda_k|^2+q_k^2
 =e^{(\rho-R_k)/k}|\zeta-e^{(R_k+i\theta)/k}|^2,
 \qquad \zeta^k=1,
\]
together with $\prod_{\zeta^k=1}(w-\zeta)=w^k-1$, gives
\begin{equation}\label{bzrw:root}
 \begin{aligned}
     I_k:=k\tau_k\log(|Z_k-\lambda_k|^2+q_k^2)
 &=\rho+\log[2(\cosh R_k-\cos\theta)]\\
 &\longrightarrow\rho+\log[2(\cosh R-\cos\theta)].
 \end{aligned}
\end{equation}
Uniformity follows since $R_k$ is bounded away from zero.
{
To calculate the resolvent trace, replace $q_k$ by a variable
$q>0$ in \eqref{bzrw:Rk} and in the equality in \eqref{bzrw:root}.
Keeping $k,\rho,\theta$ fixed, differentiate with respect to $k^2q^2$
and then set $q=q_k$. We obtain
}
\begin{equation}\label{bzrw:P}
 P_k:=\frac1k\tau_k[(|Z_k-\lambda_k|^2+q_k^2)^{-1}]
 =\frac{e^{-\rho/k}\sinh R_k}
 {2k\sinh(R_k/k)(\cosh R_k-\cos\theta)}\longrightarrow\frac1t.
\end{equation}
For $G_\mu(w)=\int(w-x)^{-1}d\mu(x)$, the subordination relation
at zero gives $q_k kP_k=1/(q_k+r_k)$, and hence
\begin{equation}\label{bzrw:r}
 r_k=\frac1{Q_kP_k}-\frac{Q_k}{k}\longrightarrow\frac tQ.
\end{equation}
Thus $r_k$ and its reciprocal are uniformly bounded on $\mathcal K$.

{
Since $\lambda_k\in\Omega_k^{\mathrm{reg}}(t)$, we may apply item (3) of Proposition \ref{prop:calculating the logarithmic integral} to the free operators in
Lemma~\ref{bzrw:freeness}. Combining that formula with
\eqref{bzrw:cyclic} and cancelling $2k\log r_k$, we obtain
}
\begin{equation}\label{bzrw:decomp}
 S_k=I_k+k\tau\log\left(1+\frac{tB}{kr_k^2}\right)
             -2k\log\left(1+\frac{q_k}{r_k}\right).
\end{equation}
In particular $S_k$ is finite. Since $\tau B=1$ and
$0\leq x-\log(1+x)\leq x^2/2$ for $x\geq0$,
\begin{align}\label{bzrw:terms}
 0\leq\frac{t}{r_k^2}
       -k\tau\log\left(1+\frac{tB}{kr_k^2}\right)
 &\leq\frac{t^2\tau(B^2)}{2kr_k^4},\notag\\
 0\leq\frac{2Q_k}{r_k}
       -2k\log\left(1+\frac{q_k}{r_k}\right)
 &\leq\frac{Q_k^2}{kr_k^2}.
\end{align}
Both bounds tend to zero uniformly on $\mathcal K$. Since
$r_k\to t/Q$, the second and third terms in \eqref{bzrw:decomp}
converge to $Q^2/t$ and $-2Q^2/t$, respectively. Together with
\eqref{bzrw:root} and $Q^2=R^2-\rho^2$, this proves
$S_k\to F_t$ uniformly on $\mathcal K$.

It remains to show that $F_t(z)=\tau\log|b_t-z|^2$.
For this purpose, we use the Hamilton--Jacobi formulas of
\cite{DHK}. Set
\begin{equation}\label{bzrw:coordinates}
 \begin{gathered}
 C=\cosh R-\frac{\sinh R}{R}\rho,\qquad
 \rho_0=-\log C,\qquad \lambda_0=e^{\rho_0+i\theta},\\
 \varepsilon_0=2e^{\rho_0}(\cosh R-\cosh\rho_0),\qquad
 p_0=\frac{e^{-\rho_0}}{2(\cosh R-\cos\theta)}.
 \end{gathered}
\end{equation}
Since $|\rho|<R$, $e^{-R}<C<e^R$. Thus $|\rho_0|<R$,
$\lambda_0\in\Sigma_t$, and $\varepsilon_0>0$.
For DHK's initial potential
$S_0(\lambda,\varepsilon)=\log(|1-\lambda|^2+\varepsilon)$,
direct substitution gives
\begin{align}\label{bzrw:initial}
 S_0(\lambda_0,\varepsilon_0)
 &=\rho_0+\log[2(\cosh R-\cos\theta)],\quad
 p_0=\partial_\varepsilon S_0(\lambda_0,\varepsilon_0),\notag\\
 t^2\varepsilon_0p_0^2
 &=\frac{R^2(2C\cosh R-C^2-1)}{\sinh^2R}
   =R^2-\rho^2.
\end{align}
Therefore
\begin{equation}\label{bzrw:action}
 F_t(z)=S_0(\lambda_0,\varepsilon_0)-t\varepsilon_0p_0^2
                         +\log\frac{|z|}{|\lambda_0|}.
\end{equation}
The characteristic with these initial data has lifetime $t$ and
endpoint $z$. Indeed, with $r_0=e^{\rho_0}$, the lifetime and
terminal-radius formulas in \cite{DHK} give
\[
 \begin{gathered}
 \delta=(r_0^2+1+\varepsilon_0)/r_0=2\cosh R,\qquad
 t_* =\frac{R}{p_0r_0\sinh R}=t,\\
 \frac t2p_0(r_0^2-1+\varepsilon_0)
 =\frac R{\sinh R}(\cosh R-C)=\rho;
 \end{gathered}
\]
the angle remains $\theta$; see
\cite[Eq.~(6.48), Propositions 6.4 and 6.11, Theorem 6.7]{DHK}.
{
We use DHK's convention in which $\varepsilon_0$ is added to
$|1-\lambda_0|^2$.
}
By the action formula and continuity at zero regularization in
\cite[Theorem 6.2 and Corollary 7.5]{DHK}, the right hand side of
\eqref{bzrw:action} is $\tau\log|b_t-z|^2$. This proves
\eqref{bzrw:result}. Finally, $F_t(z)$ is bounded on $\mathcal K$,
so exponentiating the locally uniform convergence of $S_k/2$
gives the asserted convergence of the determinants.
\end{proof}

For $z=e^{\rho+i\theta}$ and $R=R_t(\theta)$, our formula gives
\[
F_t(z)=\tau\log|b_t-z|^2
=\rho+\log[2(\cosh R-\cos\theta)]
-\frac{R^2-\rho^2}{t}.
\]
Since $R$ is independent of $\rho$, we have
\[
\partial_\rho^2F_t=\frac2t.
\]
By the chain rule,
\begin{align*}
\partial_\theta F_t
&=
\frac{\sinh R_t(\theta)\,R_t'(\theta)+\sin\theta}
{\cosh R_t(\theta)-\cos\theta}
-\frac{2R_t(\theta)R_t'(\theta)}{t}\\
&=
\frac{\sin\theta}{\cosh R_t(\theta)-\cos\theta}
+
R_t'(\theta)
\left[
\frac{\sinh R_t(\theta)}
{\cosh R_t(\theta)-\cos\theta}
-\frac{2R_t(\theta)}{t}
\right].
\end{align*}
The expression in brackets vanishes by the defining equation
$K_\theta(R_t(\theta))=1/t$.
Consequently,
\[
\partial_\theta F_t
=\frac{\sin\theta}{\cosh R_t(\theta)-\cos\theta}.
\]

Using
$\Delta_z=e^{-2\rho}(\partial_\rho^2+\partial_\theta^2)$,
we obtain, for $z=re^{i\theta}\in\Sigma_t$,
\[
\frac{d\mu_{b_t}}{dA}(re^{i\theta})
=\frac{1}{4\pi}\Delta_zF_t
=\frac{1}{4\pi r^2}
\left[
\frac2t+
\frac{d}{d\theta}
\left(
\frac{\sin\theta}{\cosh R_t(\theta)-\cos\theta}
\right)
\right]
=\frac{w_t(\theta)}{r^2}.
\]
This is the density formula obtained in \cite{DHK}.
Since $dA=r^2\,d\rho\,d\theta$, the density in logarithmic
polar coordinates is $w_t(\theta)$, independent of $\rho$.

\section*{Acknowledgements}
P.~Zhong is supported in part by NSF grant LEAPS-MPS-2516951 and NSF CAREER Award DMS-2516950. He also wants to thank Brian Hall and Ching-Wei Ho for useful discussions on free random walk in April 2026. 

The main ideas and proof strategies in this paper are due to the authors. Most results were obtained in an earlier version of the manuscript, before the availability of the AI tools used in this work. We used ChatGPT-6 Astra while working on Sections 4 and 6. For the absolute continuity result in Section 4.1, the authors supplied the proof strategy and key arguments. For the analysis of atoms, the authors proposed studying the local behavior of the logarithmic potential. The AI suggested a formulation of Definition 4.7 for calculating the mass of an atom, which the authors recognized as the Lelong number of a subharmonic function and used in the subsequent analysis. Assistance on computational details on Section 6 was provided by the AI. The authors independently verified all mathematical arguments and take full responsibility for the content of the paper.

\end{document}